\documentclass[11pt]{article}
\RequirePackage[OT1]{fontenc}
\RequirePackage{amsthm,amsmath}
\RequirePackage{hyperref}
\usepackage{amsfonts}
\usepackage[mathscr]{eucal}
\usepackage{amssymb}
\usepackage{amstext}
\usepackage{parskip}
\usepackage{fullpage}
\usepackage{graphicx}
\usepackage{subcaption}
\usepackage{lscape}

\RequirePackage[numbers]{natbib}
\usepackage{amsfonts,amssymb,amsthm,amsmath,enumerate}

\RequirePackage[OT1]{fontenc}
\RequirePackage{amsthm,amsmath}

\theoremstyle{plain}
\newtheorem{theorem}{Theorem}[section]

\newtheorem{lemma}[theorem]{Lemma}
\newtheorem{proposition}[theorem]{Proposition}

\newtheorem{definition}[theorem]{section}

\newtheorem{remark}[theorem]{Remark}

{\end{eqnarray}\end{subequations}\hskip-4.0pt}%

\newenvironment{proofof}[1]{\hspace*{20pt}{\it Proof}{ of #1}.\hskip10pt}{\qed\vskip5pt}
\newenvironment{proofof2}{}{\qed\vskip5pt}

\newcommand{\R}{\ensuremath{\mathbb R}}
\newcommand{\Q}{\ensuremath{\mathcal Q}}

\newcommand{\Sc}{\ensuremath{S^c}}

\newcommand{\prob}[1]{\ensuremath{{\mathbb P}\left(#1\right)}}

\newcommand{\size}[1]{\ensuremath{\left|#1\right|}}

\newcommand{\argmin}{\operatorname{argmin}}

\newcommand{\up}{\upsilon}
\newcommand{\e}{\epsilon}
\newcommand{\ve}{\varepsilon}

\newcommand{\silent}[1]{}

\newcommand{\ip}[1]{\;\langle{\,#1\,}\rangle\;}

\newcommand{\Cone}{{\mathcal C}}
\newcommand{\T}{{\mathcal T}}

\newcommand{\RE}{\textnormal{\textsf{RE}}}

\newcommand{\trace}{\textsf{Tr}}
\newcommand{\basepen}{\ensuremath{\lambda_{\sigma,a,p}}}

\newcommand{\lars}{\textsf{LARS}}

\newcommand{\init}{\text{\rm init}}
\newcommand{\drop}{{\ensuremath{\mathcal{D}}}}
\newcommand{\dropS}{\ensuremath{{S_{\drop}}}}

\newcommand{\inv}[1]{\frac{1}{#1}}
\newcommand{\abs}[1]{\left\lvert#1\right\rvert}
\newcommand{\twonorm}[1]{\left\lVert#1\right\rVert_2}

\newcommand{\shtwonorm}[1]{\lVert#1\rVert_2}
\newcommand{\shonenorm}[1]{\lVert#1\rVert_1}
\newcommand{\norm}[1]{\left\lVert#1\right\rVert}

\newcommand{\bens}{\begin{eqnarray*}}
\newcommand{\eens}{\end{eqnarray*}}
\newcommand{\ben}{\begin{eqnarray}}
\newcommand{\een}{\end{eqnarray}}
\newcommand{\bit}{\begin{enumerate}}
\newcommand{\eit}{\end{enumerate}}

\newcommand{\beq}{\begin{equation}}
\newcommand{\eeq}{\end{equation}}

\newcommand{\E}{{\textbf{E}}}
\newcommand{\supp}{\mathop{\text{\rm supp}\kern.2ex}}
\newcommand{\OLS}{\mathop{\text{\rm ols}\kern.2ex}}
\newcommand{\ext}{\mathop{\natural}\kern.2ex}

\begin{document}

\title{Thresholded Lasso for high dimensional variable selection}

\author{Shuheng Zhou \\
  University of California, Riverside}

\date{}

\maketitle

\begin{abstract}
Given $n$ noisy samples with $p$ dimensions, where $n \ll p$, we show 
that the multi-step thresholding procedure based on the 
Lasso -- we call it the {\it Thresholded Lasso}, can accurately estimate a 
sparse vector $\beta \in \R^p$ in a linear model $Y = X \beta + \epsilon$, 
where $X_{n \times p}$ is a design matrix normalized to have column 
$\ell_2$-norm $\sqrt{n}$ and $\epsilon \sim N(0, \sigma^2 I_n)$.
Here $I_n$ denotes the identity matrix.
We show that under the restricted eigenvalue (RE) condition,
it is possible to achieve the $\ell_2$ loss within a logarithmic factor of the ideal mean square error one would achieve with an {\em oracle }  while selecting a sufficiently sparse model -- 
hence achieving {\it sparse oracle inequalities}; the oracle would supply perfect 
information about which coordinates are non-zero and which are above the 
noise level. We also show for the Gauss-Dantzig selector (Cand\`{e}s-Tao 07), if $X$ 
obeys a uniform uncertainty principle,
one will achieve the sparse oracle inequalities as above, while
allowing at most $s_0$ irrelevant variables in the model  in the worst case,
where $0< s_0 \leq s$ is the smallest integer such that for $\lambda = \sqrt{2 \log
p/n}$, $\sum_{i=1}^p \min(\beta_i^2, \lambda^2 \sigma^2) \leq s_0 \lambda^2
\sigma^2$.  Our simulation results on the Thresholded Lasso match our
theoretical analysis excellently.
\end{abstract}

\section{Introduction}
\label{sec:introduction}
In a typical high dimensional setting, the number of variables $p$ is much 
larger than the number of observations $n$. This challenging setting appears
in linear regression, signal recovery, covariance selection in graphical 
modeling, and sparse approximation.
In this paper, we consider recovering a vector $\beta \in \R^p$ in the following 
linear model:
\beq
\label{eq::linear-model}
Y = X \beta + \epsilon.
\eeq
Here $X$ is an $n \times p$ design matrix, $Y$ is a vector
of noisy observations, and $\epsilon$ is the noise term. 
We assume throughout this paper that $p \ge n$ (i.e. high-dimensional)
and $\e$ is a vector of i.i.d. normal $N(0, \sigma^2)$ random
variables. Denote by $[p] = \{1, \ldots, p\}$. 
The notation $\twonorm{\beta}= (\sum_{j} \beta_j^2)^{1/2}$
stands for the $\ell_2$ norm of $\beta$.
Given such a linear
model, two key tasks are:
(1) to select the relevant set of variables and (2) to estimate
$\beta$ with bounded $\ell_2$ loss.
In particular, recovery of the sparsity pattern 
$S = \supp(\beta) := \left\{j \,:\, \beta_j \neq 0\right\}$, 
also known as variable (model) selection,  refers to the task of
correctly identifying the support set, 
or a subset of ``significant'' coefficients in $\beta$, based on the
noisy observations.
Here and in the sequel, we assume that each 
column vector $X_j \in \R^n, j \in [p]$ of the fixed design matrix $X$ has $\ell_2$-length
$\twonorm{X_j} = \sqrt{n}$.

Even in the noiseless case, recovering $\beta$ (or its support) 
from $(X, Y)$ seems impossible when $n \ll p$ given that we have more
variables than observations.
Over the past two decades, a line of research shows that 
when $\beta$ is sparse, that is, when it has a relatively small number of 
nonzero coefficients, and when design matrix $X$ behaves sufficiently 
nicely in a sense that it satisfies certain incoherence conditions, it becomes possible to reconstruct
$\beta$~\citep{Donoho:cs,CT05,CT07}.

Throughout this paper, we refer to a vector $\beta \in \R^{p}$ with at 
most $s$ non-zero entries, where $0 \le s \leq p$, as an
{\bf $s$-sparse} vector.
Consider now the linear regression model in~\eqref{eq::linear-model}.
For a chosen penalization parameter $\lambda_n \geq 0$, regularized 
estimation with the $\ell_1$-norm penalty, 
also known as the Lasso \citep{Tib96}
refers to the following convex optimization problem  
\begin{eqnarray}
\label{eq::origin} \; \; 
\hat \beta = \arg\min_{\beta \in \R^p} \frac{1}{2n}\|Y-X\beta\|_2^2 + 
\lambda_n \|\beta\|_1,
\end{eqnarray}
where the scaling factor $1/(2n)$ is chosen by convenience and
$\|\beta\|_1 = \sum_{i=1}^p \abs{\beta_i}$.
In the present work, we explore model selection beyond focusing on the
notion of exact recovery of the support set $\supp(\beta)$
which crucially depends on the so called $\beta_{\min}$ condition as well as the {\it
  Neighborhood stability} or {\it irrepresentability
  condition}~\citep{MB06,ZY06,Wai09,Wai09b}.
One can not hope that such incoherence conditions always hold in reality.
As pointed out by~\cite{GBZ11}, the {\it irrepresentability condition}
which is essentially a necessary condition for exact recovery of the
non-zero coefficients (for which a $\beta_{\min}$ condition needs to
hold) by the Lasso, is much too restrictive in comparison to the
Restricted Eigenvalue condition~\citep{BRT09}; cf.~\eqref{eq::admissible}.
For some integer $s \in [p]$ and a positive number $k_0 >0$, 
we say $\RE(s, k_0, X)$ holds with $K(s, k_0)$ if for all $\upsilon \not=0$,
\beq
\label{eq::admissible}
\inv{K(s, k_0)} \stackrel{\triangle}{=}
\min_{J \subseteq [p], |J| \leq s} \min_{\norm{\upsilon_{J^c}}_1 \leq k_0 \norm{\upsilon_{J}}_1}
\; \;  \frac{\norm{X \upsilon}_2}{\sqrt{n}\norm{\upsilon_{J}}_2} > 0
\eeq
where $\upsilon_{J}$ represents the subvector of $\upsilon \in \R^p$
confined to a subset $J$ of $[p]$.   It is clear that as $k_0$ and
$s_0$ become smaller, this condition is easier to satisfy.
To be clear, $\RE$ conditions alone are not sufficient for the Lasso
to recover the model $S$ exactly. 
Moreover, to ensure variable selection consistency, a
\ben
\label{eq::betamin}
\text{ $\beta_{\min}$ condition:} \quad \min_{j \in S} \abs{\beta_j} \ge C \sigma \sqrt{{2 \log p}/{n} },
\een
is imposed for some constant $C > 1/2$, and shown to be crucial
to recover the support of $\beta$ in the information theoretic limit
by~\cite{Wai09b} and~\cite{Zhang10}.
Such a $\beta_{\min}$ condition and the corresponding
signal-to-noise ratio (SNR) defined as $\beta_{\min}^2/\sigma^2$, rather than the typical $\shtwonorm{\beta}^2/\sigma^2$, is shown to be
the key quantity that controls subset selection~\citep{Wai09b}.

Ideally, we aim to remove or relax the $\beta_{\min}$ condition, which is
rather unnatural for many applications.
Toward this end,
we define {\bf sparse oracle inequalities} 
as the new criteria for {\it model selection consistency} when the
{\it irrepresentability condition} or related mutual incoherence
conditions are replaced with the Restricted Eigenvalue ($\RE$) type of
conditions.
Roughly speaking, the new criteria ask one to identify a sparse model
such that the corresponding least-squares (OLS) estimator based on the
selected model achieves an oracle inequality in terms of the $\ell_2$
loss while keeping the selection set small. We deem the bound on the
$\ell_2$-loss as a natural criterion for evaluating a sparse model especially when it is not exactly $S$.
We achieve this goal by controlling the false positive selection through thresholding
initial estimates of $\beta$  obtained via the Lasso (or the Dantzig selector) at the critical threshold level.

\noindent{\bf Contributions.}
Our contributions in this work are twofold. 
From a methodological point of view, we propose to study the {\bf Thresholded
  Lasso} estimator with the aforementioned goals in mind:
\begin{itemize}
\item[Step 1]
First, we obtain an initial estimator $\beta_{\init}$ using the 
Lasso~\eqref{eq::origin} with $\lambda_n = d_0 \sigma \sqrt{2\log p/n}$,
for some constant $d_0 > 0$, which is allowed to depend on sparse and restricted eigenvalue parameters;
\item[Step 2]
\label{first-step}
Threshold the initial Lasso estimator $\beta_{\init}$ with $t_0$,
with the general goal such that, we get a set $I$ with cardinality at 
most $2s$; in general, we also have $|I \cup S| \leq 2s$, where 
$I =  \left\{j \in [p]: \abs{\beta_{j, \init}} \geq t_0 \right\}$ 
for some $t_0 \asymp \sigma \sqrt{2 \log p/n}$ with hidden constant to
be specified;
\item[Step 3]
Feed $(Y, X_{I})$ to the ordinary least squares (OLS) estimator to
obtain $\hat{\beta}$, where we set
$\hat\beta_{I} = (X_I^T X_{I})^{-1} X_{I}^T Y$ and the other
coordinates to zero.
\end{itemize}
In Theorem~\ref{thm::RE-oracle-main}, we show that the critical
threshold level for estimating a high dimensional sparse vector $\beta
\in \R^p$ should be set at the level $t_0 \asymp \lambda \sigma$ for
$\lambda = \sqrt{2 \log  p/n}$, to retain signals in $\beta_{\init}$ 
at or above the level of  $\lambda \sigma$, where 
$\beta_{\init}$ is the solution to the Lasso
estimator~\eqref{eq::origin} obtained in Step 1.
Moreover, we show that $d_0$ and $t_0$ are allowed to depend 
on sparse and restricted eigenvalue parameters of the design 
matrix; cf. Section~\ref{sec::relaxed}.
From a theoretical point of view, the framework for our analysis is 
set upon the Restricted Eigenvalue  type of condition and an upper 
sparse eigenvalue (USE) condition, namely,
\ben 
\label{eq::eigen-max}
\text{(USE)} \quad
\Lambda_{\max}(2s) \; \stackrel{\triangle}{=} \; 
\max_{\upsilon \not=0; 2s-\text{sparse}} \; 
\twonorm{X \upsilon}^2/(n\twonorm{\upsilon}^2) < \infty. 
\een
These are among the most general assumptions on the design 
matrix, guaranteeing sparse recovery and oracle inequalities in the $\ell_2$ loss for the 
Lasso estimator as well as the Dantzig selector: $\RE(s, k_0, X)$ is shown to be a relaxation of the restricted isometry property (RIP) under suitable choices of 
parameters involved in each condition~\citep{CT07,BRT09}.
Part of this work was presented in a conference paper 
by~\cite{Zhou09th}. Importantly,  we present significant and novel
extensions in both theory and numerical simulations, with regards to the Thresholded Lasso and
the Lasso under the $\RE$ and sparse eigenvalue conditions.

While the crucial theoretical and methodological ideas presented here 
originate from~\cite{Zhou09th}, the current work significantly expands
the original ideas and show new results on the sparse oracle inequalities in
Theorems~\ref{thm::RE-oracle-main} and~\ref{thm::RE-oracle}.
Compared with the original paper~\citep{Zhou09th}, we further 
study the behavior of the Thresholded Lasso in several challenging
situations in Sections~\ref{sec:type-II-intro},~\ref{sec::numeric},
and the supplementary Section~\ref{sec:experiments}.
We find that our estimator is robust and adaptive to the overwhelming presence of weak signals in
$\supp(\beta)$ in both theoretical and practical senses.
We show that the Thresholded Lasso tradeoffs false positives and false negatives nicely in this case: its
advantage in terms of model selection over the Lasso and adaptive
Lasso~\citep{Zou06,ZGB09} is clearly evident by examining their ROC
curves empirically.  Our numerical simulations in
Section~\ref{sec:experiments} show that the rates for exact recovery
of the support rise sharply for a few types of random matrices once the number of samples passes a certain threshold, using the Thresholded
Lasso estimator.

\noindent{\bf Notation and definitions.}
Let $T \subseteq [p]$ be a fixed subset of indices.
Let $X_T$ be the $n \times  \abs{T}$ submatrix obtained by extracting columns of $X$ indexed by $T$.
We use $\upsilon_{T}$ to represent the subvector of
$\upsilon \in \R^p$ confined to a subset $T$ of $[p]$.
Let $\shtwonorm{\upsilon}^2 = \sum_{j=1}^p \upsilon_j^2$.
Occasionally, we use $\upsilon_T \in \R^{\abs{T}}$ to also represent
its $0$-extended version $\upsilon' \in \R^p$ such that
$\upsilon'_{T^c} = 0$ and $\upsilon'_{T} = \upsilon_T$.
We will use $\beta^{\ext}(T) \in \R^p$ to represent the $0$-extended
version of $\beta_T \in \R^{\abs{T}}$ such that $\beta_{T}^{\ext}({T}) = \beta_{T}$
and $\beta_{T^c}^{\ext}({T}) = 0$.
For a matrix $A$, let $\Lambda_{\min}(A)$ and $\Lambda_{\max}(A)$
denote the smallest and the largest eigenvalues respectively.
Let {\bf $s = |S|$}. We assume the lower sparse eigenvalue (LSE)
condition, namely,
\beq
\label{eq::eigen-admissible-s}
\text{(LSE)} \quad \Lambda_{\min}(2s) \;
\stackrel{\triangle}{=} \;
\min_{\upsilon \not= 0; 2s-\text{sparse}} \; \;
\twonorm{X \upsilon}^2/(n \twonorm{\upsilon}^2) > 0,
\eeq
where $n \geq 2s$ is necessary, as any submatrix with more than $n$ 
columns must be singular.
Also relevant is the $(s, s')$-{\em restricted orthogonality constant}
$\theta_{s, s'}$~\citep{CT07}, which is defined to be the smallest quantity such that for all disjoint sets $T,T' \subseteq [p]$ of cardinality $|T| \leq s$ and $|T'| \leq s'$:
\ben
\label{label:correlation-coefficient}
&& {\abs{{\ip{X_T \upsilon, X_{T'} {\upsilon'}}}}}/{n} \leq 
\theta_{s,  s'} \twonorm{\upsilon} \twonorm{\upsilon'}, \; \text{
  where } s + s' \leq p, \\
  &&\text{ and } \;  \theta_{s,  s'} \leq (\Lambda_{\max}(s) \Lambda_{\max}(s'))^{1/2} \; 
\text{ by the Cauchy Schwarz inequality.}
\een
Note that 
small values of $\theta_{s, s'}$ indicate that disjoint subsets of
covariates in $X_T$ and $X_{T'}$ span nearly orthogonal subspaces.
Moreover, we have
$\theta_{s, s'} \leq {(\Lambda_{\max}(s+s')
  -\Lambda_{\min}(s+s'))}/{2}$; cf. Lemma~\ref{lemma:parallel}.
Technically speaking, each of the entities defined above, namely,
$1/\Lambda_{\min}(2s)$, $\Lambda_{\max}(2s)$, $\theta_{s, s'}$, and $K(s, k_0)$ as introduced in~\eqref{eq::admissible}, 
is a non-decreasing function of $s$, $s'$, and $k_0$.
Nonetheless, we crudely consider these as constants following how they
are typically treated in the literature as it is to be understood that
they grow very slowly with $s$ and $s'$; see for
example~\cite{CT05,CT07},~\cite{MY09}, and~\cite{BRT09}.
For two numbers $a, b$, $a \wedge b:= \min(a,b)$, and $a \vee b:= \max(a,b)$.
We write $a \asymp b$ if $ca \le b \le Ca$ for some positive absolute
constants $c,C$ which are independent of $n, p$, and $\gamma$.
We write $f = O(h)$ or $f \ll h$ if $\abs{f} \le C h$ for some absolute constant
$C< \infty$ and $f=\Omega(h)$ or $f \gg h$ if $h=O(f)$.
We write $f = o(h)$ if $f/h \to 0$ as $n \to \infty$.

\subsection{Sparse oracle inequalities}
\label{sec:oracle-intro}
In this section, we define {\bf sparse oracle inequalities} as the new 
criteria for model  selection consistency when some of the signals in
$\beta$ are relatively weak, for example, well below the information
theoretic  detection limit~\eqref{eq::betamin} for high dimensional sparse 
recovery. While the idea of thresholding and 
refitting is widely used in statistical theory and applications in
various contexts, we quantify the threshold level based on the oracle
$\ell_2$ loss for the Lasso (and Dantzig selector respectively) in the
present work, with the following goals.

Specifically, (a) we wish to obtain $\hat{\beta}$ such that $|\supp(\hat{\beta}) \setminus S|$ 
(and sometimes the set difference between $S$ and $\supp(\hat{\beta})$
denoted by $|S \triangle \supp(\hat{\beta})|$ also) is small, 
with high probability; (b) while at the same time, we wish to bound
$\shtwonorm{\hat{\beta} - \beta}^2$,
within logarithmic factor of the ideal mean squared error one would achieve with an oracle that would
supply perfect information about which coordinates are non-zero and
which are above the noise level (hence achieving the {\it oracle inequality} as 
studied by~\cite{Donoho:94} and~\cite{CT07}).
Here we denote by $\beta_I$ the restriction of $\beta$ to the set $I$,
where $I \subset [p]$, and $\beta^{\ext}({I})$ its $0$-extended version.

Formally, we evaluate the selection set through the following
criterion. Consider the least squares estimators 
$\hat{\beta}_I = (X_I^T X_{I})^{-1} X_{I}^T Y$, where $I \subset [p]$
and $ |I| \ \leq s$.
Here and in the sequel, let $\hat{\beta}_I^{\OLS}(I) = \hat{\beta}_I$ and $\hat{\beta}_{I^c}^{\OLS}(I) =0$.
Consider the {\it ideal} least-squares estimator $\beta^{\diamond}$
based on a subset $I$ of size at most $s$, which minimizes the mean squared error:
\begin{eqnarray}
\label{eq::beta-diamond}
\beta^{\diamond} = \argmin_{I \subseteq [p],\; |I| \leq  s}  \E
  \twonorm{\beta- \hat{\beta}^{\OLS}(I)}^2.
\end{eqnarray}
It follows from the analysis by~\cite{CT07} that for
$\Lambda_{\max}(s) < \infty$,
\begin{eqnarray}
   \label{eq::ME-diamond}
\E \twonorm{\beta- \beta^{\diamond}}^2 & \geq & 
\min \left(1, {1}/{\Lambda_{\max}(s)}
                                                \right)\sum_{i=1}^p\min(\beta_i^2,\sigma^2/n),
                                                \text{ where } \\
  \nonumber
\sum_{i=1}^p \min(\beta_i^2, \sigma^2/n) & = &  \min_{I \subseteq [p]}
\twonorm{\beta - \beta^{\ext}({I})}^2 + {|I| \sigma^2}/{n} 
\end{eqnarray}
represents the squared bias and variance.
Now we check if~\eqref{eq::log-MSE}
\begin{eqnarray}
\label{eq::log-MSE}
\shtwonorm{\hat{\beta} - \beta}^2 & = &  
O(\lambda^2 \sigma^2 + \sum_{i=1}^p \min(\beta_i^2, \lambda^2
                                        \sigma^2)) 
\end{eqnarray}
holds with high probability;
If so, we claim the following holds:
\begin{eqnarray}
\label{eq::log-MSE-E}
\shtwonorm{\hat{\beta} - \beta}^2 =
O_P(\log p \max(1, \Lambda_{\max}(s)) 
  \E\twonorm{\beta^{\diamond}- \beta}^2),
\end{eqnarray}
in view of~\eqref{eq::ME-diamond}; cf. the supplementary Section~\ref{section:append-diamond}.
Here the $\ell_2$ loss in \eqref{eq::log-MSE} is optimal up to a $\log p$
factor. We note that~\eqref{eq::log-MSE-E} is not the tightest upper bound 
that we could derive due to a relaxation we have on the lower bound
as stated in~\eqref{eq::ME-diamond}.
Nevertheless, we  use it for its simplicity.

\noindent{\bf Essential sparsity and objectives.}
The current paper answers the following question: 
Is there a good thresholding rule that enables us to obtain a 
sufficiently {\it sparse} estimator $\hat{\beta}$ that 
satisfies an {\it oracle inequality} in the sense of~\eqref{eq::log-MSE}, 
when some components of $\beta_S$ are well below $\sigma/\sqrt{n}$? 
Such oracle results are accomplished without any knowledge of the 
significant coordinates or parameter values of $\beta$. 
Both Theorem~\ref{thm::RE-oracle-main} and the supplementary 
Theorem~\ref{thm:ideal-MSE-prelude} 
answer this question positively, where we elaborate upon the sparse 
recovery properties of the Lasso and Dantzig selector in combination 
with thresholding and refitting.

For a given pair of $(n, p)$ values, the essential sparsity parameter $s_0$ characterizes more 
accurately than $s$, the number of significant coefficients of 
$\beta$ with respect to the noise level $\sigma$
that we should (could) try to recover. 
Denote by $s_0$ the smallest integer such that
the following holds~\citep{CT07}:
\ben
\label{eq::define-s0}
\sum_{i=1}^p \min(\beta_i^2, \lambda^2 \sigma^2) \leq 
s_0 \lambda^2 \sigma^2, \text{ where } \lambda = \sqrt{ 2 \log p/n}.
\een
The parameter $s_0$ is relevant especially when
we do not wish to impose any lower bound on $\beta_{\min}$.
This is the focus of Theorem~\ref{thm::RE-oracle-main}.
To make this statement precise, we have as a consequence of the 
definition in~\eqref{eq::define-s0},
\ben
\label{eq::beta-2-small-intro}
|\beta_j| < \lambda \sigma\; \; \; \text{ for all } j > s_0, \; \text{ if we order } \; |\beta_1| \geq |\beta_2| ... \geq |\beta_p|;
\een
cf. Remark~\ref{rem::sparsity}.
For simplicity of presentation, we set $|I| < 2 s_0$ as our first
goal while achieving the oracle inequality as  in~\eqref{eq::log-MSE}.
One could aim to bound $|I| < c s_0$ for some other constant $c > 0$.
Moreover, to put the bound of $\size{I} \leq 2 s_0$ in perspective, 
we show in Proposition~\ref{PROP:COUNTING-S0} (by setting $c' = 1$) 
that the number of variables in $\beta$ that are larger than or equal to 
$\sigma \sqrt{\log p/n}$ in magnitude is bounded by $2s_0$.
Roughly speaking, we wish to include most of them by taking $2s_0$
as the upper bound on the model size $I$.

The Thresholded Lasso algorithm is constructive in that it relies neither on the 
unknown parameters $\abs{S}$ or {\bf $\beta_{\min} := \min_{j \in S}
  |\beta_{j}|$}, nor the exact knowledge of those that characterize
the incoherence conditions on $X$. Instead, our choices of
penalty $\lambda_n$ (Step 1) and threshold $t_0$ (Step 2) only depend on
$\sigma, n$, and $p$, and some crude estimation of certain sparse eigenvalue parameters;
cf. Section~\ref{sec::relaxed}. In practical settings, one can choose
$\lambda_n$ using cross-validation; See for example the subsequent work by~\cite{ZRXB11}, where we use
cross-validation to choose both penalty and threshold parameters in
the context of covariance selection based on Gaussian graphical
models.

In Section~\ref{sec:type-II-intro}, we briefly discuss possibilities of 
recovering a subset of strong signals via thresholding,
despite the existence of (or influence from) other relatively weaker
signals. Let $T_0$ denote the largest $s_0$ coordinates of $\beta$ in absolute 
values. As a consequence of the definition in~\eqref{eq::define-s0}, 
we have $T_0 = \{1, \ldots, s_0\}$ and $|\beta_j| < \lambda \sigma$
for all $j \in T_0^c$; cf.~\eqref{eq::beta-2-small-intro}.
More precisely, we decompose
$T_0 = \{1, \ldots, s_0\}$ into two sets: $A_0$ and $T_0 \setminus
A_0$, where $A_0$ contains the set of coefficients of $\beta$ strictly
above $\lambda \sigma$, for which we define a constant
$\beta_{\min,A_0}$:
\ben 
\label{eq::betaminA0}
\beta_{\min,A_0}:= \min_{j \in  A_0} |\beta_{j}| >
\lambda \sigma & \text{ where} &A_0 = \{j: |\beta_j| > \lambda \sigma
\}.
\een
The goal
is to demonstrate the remarkable properties
of the Lasso and the Thresholded Lasso estimators:
while exact recovery of all non-zero variables requires very stringent
incoherence and $\beta_{\min}$ conditions, we can significantly relax
both conditions when we only require a subset $A_0$ of active
variables to be included in our selection set.
We loosely refer to $A_0$ or its superset $T_0 \supseteq A_0$ which we
aim to identify as an {\emph active set} throughout this work.
When $\beta_{\min,A_0}$~\eqref{eq::betaminA0} is sufficiently large,
we have $A_0 \subset I$ while achieving the sparse oracle
inequalities in the sense of~\eqref{eq::log-MSE};
cf.~\eqref{eq::ideal-t0} and Theorem~\ref{thm::threshA0}.

One of the reviewers brought to our attention that~\cite{ZH08} provide bounds
similar to the Dantzig selector under the upper and lower sparse Riesz
condition (SRC), which are similar to the upper and lower
bounds in~\eqref{eq::eigen-cond} in the present work.
Both models allow potentially many small coefficients in the true
$\beta$.
Specifically, we design a new set of experiments to evaluate the
impact of sparsity $s$ and  the $\beta_{\min, A_0}$ (e.g., $C_a \lambda
\sigma$ in \eqref{eq::tigermodel}) condition on the recovery of the
first $s_0$ components in $\beta$, where we set $\beta_{T_0^c}$ (with
support size $s - s_0$) to have a fixed $\ell_2$ norm but potentially
many small coordinates with magnitude $\asymp \sigma/\sqrt{n}$.
See Section~\ref{sec::numeric} for the setup of numerical 
simulations.

Not included in the present work are Theorem 2.1 and the Iterative Procedure by~\cite{Zhou09th},
where we show conditions under which one can recover a sparse subset of strong 
signals when $\beta_{\min} := \min_{j \in S} \abs{\beta_j} \ge C
\sigma \sqrt{{2 s \log p}/{n} }$, where $S = \supp(\beta)$, $s =
\abs{S}$ and $C$ depends on the restricted eigenvalue parameter;
cf. Theorems 3.1 in~\cite{Zhou10}.
When $\beta_{\min}$ is sufficiently large, the range of thresholding parameters is even more flexible, 
which we elaborate in~\cite{Zhou09th} and~\cite{Zhou10}, cf. Theorem
3.1, and hence details are omitted from the current paper.
We do show numerical examples for which the 
Thresholded Lasso recovers the support $S$ {\it exactly} with high probability, 
using {\emph a small number of samples} per non-zero component in $\beta$, 
for which the Lasso would certainly have failed, as predicted by the 
work of~\cite{Wai09,Wai09b}.
These results have been presented in part in an earlier version of the present 
paper~\citep{Zhou10} and the conference paper by~\cite{Zhou09th}.

Finally, we define a quantity $\basepen$, 
which bounds the maximum correlation between the noise and 
covariates of $X$; For each $a \geq 0$, let
\ben
\label{eq::low-noise}
{\T_a} := 
\biggl \{\e: \norm{{X^T \e}/{n}}_{\infty} \leq \basepen, 
\text{ where }
\basepen = \sigma \sqrt{1 + a} \sqrt{{2\log p}/{n}}\biggr \}.
\een
Then, we have $\prob{\T_a}  \geq 1 - (\sqrt{\pi \log p} p^a)^{-1}$
when $X$ has column $\ell_2$ norms bounded by $\sqrt{n}$.

\noindent{\bf Organization of the paper.}
The rest of the paper is organized as follows.
In Section~\ref{sec::theory}, we describe a thresholding framework
for the general setting, and highlight the role thresholding plays in
terms of recovering the best subset of variables; we present the main
Theorem~\ref{thm::RE-oracle-main} and oracle results for the Lasso
estimator.
We discuss related work in Section~\ref{sec:linear-sparsity}.
Section~\ref{sec::proofsketchRE} provides the proof sketch of
Theorems~\ref{thm::RE-oracle-main} and~\ref{thm::RE-oracle},
and the proof of Theorem~\ref{thm::RE-oracle-main} appears in Section~\ref{sec:proof-TH-main}.
Section~\ref{sec:type-II-intro} discusses Type II errors and the
$\ell_1$- and $\ell_2$-loss. Section~\ref{sec::numeric} and the
supplementary Section~\ref{sec:experiments} include simulation
results. We prove
Lemmas~\ref{prop:MSE-missing},~\ref{prop:MSE-missing-orig},~\ref{lemma:threshold-general},
and Lemma~\ref{lemma:threshold-general-II} in
Sections~\ref{sec::OLSproof-orig} to~\ref{sec::proofofA0}.
We give a sketch proof of Theorem 3~\citep{ZH08} in the
supplementary Section \ref{sec::ZH08proof}.
We conclude in Section~\ref{sec:conclude}. 
Additional technical proofs are included in the supplement.

\section{The Thresholded Lasso estimator}
\label{sec::theory}
Theorem~\ref{thm::RE-oracle-main} states
that sparse oracle inequalities
as elaborated in Section~\ref{sec:oracle-intro} hold for the Thresholded Lasso
under no restriction on $\beta_{\min}$.
We do not optimize constants in this paper.
Theorem~\ref{thm::RE-oracle-main} is the key contribution of this 
paper and is proved in 
Section~\ref{sec:proof-TH-main}. 
\begin{theorem}\textnormal{(\bf{Ideal model selection for the Thresholded Lasso})}
\label{thm::RE-oracle-main}
Suppose $\beta \in \R^p$ is {\bf $s$-sparse}.
Let $s_0$ be as in~\eqref{eq::define-s0}.
Let $Y = X \beta + \e$, where $\e =(\e_1, \ldots, \e_n)^T$ is a vector
containing independent and identically distributed (i.i.d.) noise with
$\e_i \sim N(0, \sigma^2)$ for all $i \in [n]$.
Suppose the columns of $X$ are normalized to have $\ell_2$ norm $\sqrt{n}$.
Suppose $\RE(s_0, 4, X)$ holds with $K(s_0, 4)$, and
the sparse eigenvalue conditions~\eqref{eq::eigen-max} and~\eqref{eq::eigen-admissible-s} hold. 
Let $\beta_{\init}$ be an optimal solution to the Lasso~\eqref{eq::origin} with 
$\lambda_{n} = d_0 \sqrt{2 \log p/n} \sigma$ $\geq 2 \basepen$, 
where $d_0 \geq 2 \sqrt{1 + a}$ for $a \geq 0$.
Set $t_0 = C_4 \lambda \sigma \ge 2 \sqrt{1+a} \lambda \sigma$, for
some constant $C_4 \geq D_1$ for $D_1$ as in~\eqref{eq::D1-define}.
Set $I =  \left\{j \in [p]: \abs{\beta_{j, \init}} \geq t_0 \right\}$. 
Set $\hat\beta_{I} = (X_I^T X_{I})^{-1} X_{I}^T Y$ and
$\hat\beta_{I^c} =0$. Then with probability at least
$1-\prob{\T_a^c}$, we obtain
\ben
\nonumber
&& |I| \leq   s_0 (1 + D_1/C_4)   < 2 s_0, \quad|I \cup S|   \leq  s
+ s_0 \; \text { and } \\
\label{eq::MSE}
&& \shtwonorm{\hat\beta - \beta}^2 \leq   D_4^2 s_0 \lambda^2
\sigma^2, \text{ where for $D'_0$ be as in~\eqref{eq::D0-define-orig}, }  \\
\label{eq::D4-constant}
&& D_4^2  \leq ((D'_0 + C_4)^2 + 1)
\left(\frac{3}{2} + \frac{(\Lambda_{\max}(2s) - 
\Lambda_{\min}(2s))^2}{2\Lambda_{\min}^2(2s_0)} \right).
\een
\end{theorem}
Theorem~\ref{thm::RE-oracle-main} relies on new oracle results that 
we prove for the Lasso estimator under $\RE(s_0, 4, X)$ and
the upper sparse eigenvalue condition~\eqref{eq::eigen-max} in Theorem~\ref{thm::RE-oracle}.
The Lasso estimator achieves essentially the same bound in terms of
$\ell_2$ loss as stated in~\eqref{eq::log-MSE}, which adapts nearly
ideally not only to the uncertainty in the support set $S$ but also
to the size of the ``significant'' set; cf. Section~\ref{sec::relaxed}. 
It is clear from our analysis in~\cite{Zhou09th} that, showing an oracle inequality 
as in~\eqref{eq::log-MSE} for the initial Lasso estimator $\beta_{\init}$, cf. 
Theorem~\ref{thm::RE-oracle}, as well as applying new
techniques for analyzing algorithms involving thresholding followed by 
OLS refitting will be crucial in proving sparse oracle inequalities for the 
Thresholded Lasso. We discuss the initial estimator $\beta_{\init}$
and the set $I$ in Section~\ref{sec::thresholdrules}.
Moreover, \eqref{eq::MSE} implies that \eqref{eq::log-MSE} holds for the 
Thresholded Lasso estimator.

\begin{remark}
  \label{rem::sparsity}
To see \eqref{eq::beta-2-small-intro}, we have by definition of $s_0$, 
where $0 \leq s_0 \leq s < p$,
\bens
\label{eq::s0-upper-bound}
s_0 \lambda^2 \sigma^2 & \leq & 
 \lambda^2 \sigma^2 + \sum_{i=1}^p \min(\beta_i^2, \lambda^2 \sigma^2) \\
& \leq &  2 \log p \big({\sigma^2}/{n} + \sum_{i=1}^p 
\min\big(\beta_i^2, {\sigma^2}/{n}\big)\big) \\
\label{eq::s0-lower-bound}
\text{ and} \; \; 
s_0 \lambda^2 \sigma^2  & \geq &  \sum_{j=1}^{s_0 + 1} 
\min(\beta_j^2, \lambda^2 \sigma^2)
\geq (s_0 + 1) \min(\beta_{s_0 +1}^2, \lambda^2 \sigma^2),
\eens
which immediately implies that
$\min(\beta_{s_0 +1}^2, \lambda^2 \sigma^2) < \lambda^2 \sigma^2$ and 
hence~\eqref{eq::beta-2-small-intro} holds.
\end{remark}

\subsection{The thresholding rules}
\label{sec::thresholdrules}
Consider the linear regression model in~\eqref{eq::linear-model}.
Let $T_0$ denote the largest $s_0$ coordinates of $\beta$ in absolute 
values, for $s_0$ as in~\eqref{eq::define-s0}.
Suppose we aim to select the set of variables of size at least 
$\sigma \sqrt{2 \log p/n}$ while controlling the bias.
Lemma~\ref{lemma:threshold-RE}
is a deterministic result toward this goal.
We prove Lemma \ref{lemma:threshold-RE} in Section~\ref{sec::proofLemmaTRE}.
\begin{lemma}({\bf A deterministic result.})
 \label{lemma:threshold-RE}
Let $\beta_{\init}$ be an initial estimator of the {\bf $s$-sparse} $\beta$
in the linear model~\eqref{eq::linear-model}, where $\e \sim N_n(0, \sigma^2 I)$
and $\twonorm{X_j} = \sqrt{n}, j \in [p]$.
Let $T_0$ denote the largest $s_0$ coordinates of $\beta$ in absolute 
values.
Let $\beta_{T_0}$ be the restriction of $\beta$ to the set  
$T_0$. Let $\beta^{\ext}({T_0}) \in \R^p$ denote its $0$-extended
version such that $\beta_{T_0^c}^{\ext}({T_0}) 
= 0$ and $\beta_{T_0}^{\ext}({T_0}) = \beta_{T_0}$. 
Let $h = \beta_{\init} - \beta^{\ext}({T_0})$.
Suppose for  $\lambda := \sqrt{{2 \log p}/{n}}$,
\ben 
\label{eq::part-norms}
\twonorm{h_{T_{0}}} &\leq &D'_0 \lambda \sigma \sqrt{s_0} \; \text{
  and } \; \norm{h_{T_0^c}}_1 \leq  D_1 \lambda \sigma s_0.
\een
Set $t_0 = C_4 \lambda \sigma$ for some positive constant $C_4$.
Let $I = \{j: \abs{\beta_{j, \init} }\geq t_0\}$ and $\drop :=
[p]\setminus I$. Then the set $I$ satisfies
\begin{eqnarray}
\label{eq::modelsize}
  |I| & \leq &  s_0 (1 + D_1/C_4),\quad  |I \cup S|
             \leq  s + s_0 D_1 /C_4, \quad \text{ and } \\
\label{eq::off-beta-norm-bound} 
\twonorm{\beta_{\drop}} & \leq & \sqrt{(D'_0 + C_4)^2 + 1} 
                                 \lambda \sigma \sqrt{s_0}, \; \; \text{ for $D'_0, D_1$ as in~\eqref{eq::part-norms}.}                 
\end{eqnarray}
\end{lemma}
Then, a tighter bound on 
$\norm{h_{T_{0}^c}}_1 := \norm{\beta_{\init, T_0^c}}_1$ or
$\shtwonorm{\beta_{\init, T_{0}^c}}$ will decrease the threshold $t_0$ while a tighter bound on 
$\shtwonorm{h_{T_{0}}} :=\shtwonorm{(\beta-\beta_{\init})_{T_0} }
\le D_0' \lambda \sigma \sqrt{s_0}$ will tighten 
the upper bound in~\eqref{eq::off-beta-norm-bound} on the bias 
component through the triangle inequality.
In general, we allow $t_0$ to be chosen from a reasonably 
wide range, where we tradeoff the width of the range with the 
tightness of the upper bound on the $\ell_2$ loss for $\hat{\beta}$ in
Step 3. This saves us from having to estimate incoherence parameters 
in a refined (and tedious) manner.
Theorem~\ref{thm::RE-oracle} may be of independent 
interest. We give a proof sketch in Section~\ref{sec::relaxed}. 
\begin{theorem}\textnormal{\bf (Oracle inequalities of the Lasso)}
  \label{thm::RE-oracle}
Suppose $\beta \in \R^p$ is {\bf $s$-sparse}.   Let $s_0$ be as in~\eqref{eq::define-s0}.
Suppose $\RE(s_0, 4, X)$ and the upper sparse eigenvalue condition~\eqref{eq::eigen-max} hold.
Under the settings of Lemma~\ref{lemma:threshold-RE},
let $\beta_{\init}$ be an optimal solution to the Lasso~\eqref{eq::origin} with $\lambda_{n} = d_0 \lambda \sigma \geq 2 \basepen$, where $a \geq 0$ and $d_0 \geq 2 \sqrt{1 + a}$. 
Let $h = \beta_{\init} - \beta^{\ext}({T_0})$ be defined as in Lemma~\ref{lemma:threshold-RE}.
Then on $\T_a$,
\bens
\twonorm{\beta_{\init} - \beta}
& \leq &  \lambda \sigma \sqrt{s_0} (\sqrt{D_0^2 + D_1^2} + 1), \\
\norm{h_{T_0}}_1 + \norm{h_{T_0^c}}_1
& =& 
\norm{\beta_{T_0}-\beta_{\init,T_0}}_1 + \norm{\beta_{\init, T_0^c}}_1  \leq
D_2 \lambda \sigma  s_0, \quad \text{and} \\
\twonorm{ X \beta_{\init} - X \beta }/\sqrt{n} & \leq & D_3 \lambda \sigma \sqrt{s_0},
\eens
where $D_0, \ldots, D_3$ are defined  in~\eqref{eq::D0-define-orig} to~\eqref{eq::D3-define}.
Moreover, for any subset $I_0 \subset S$, 
by assuming that $\RE(|I_0|, 4, X)$ holds with $K(|I_0|, 4)$, we have
\ben
\label{eq::pred-error-gen}
\twonorm{ X \beta_{\init} - X \beta}/\sqrt{n} \leq 
\twonorm{X \beta - X_{I_0} \beta_{I_0}}/\sqrt{n} + \frac{3}{2}  K(|I_0|, 4) \lambda_n \sqrt{|I_0|}.
\een
\end{theorem}
We compare Theorem~\ref{thm::RE-oracle} with the $\ell_p$ error bounds
by~\cite{BRT09} (cf. Theorem 7.2) in Section~\ref{sec::relaxed}.
The sparse oracle properties of the Thresholded Lasso in terms of variable selection, $\ell_2$ loss, and
prediction error then follow from Theorem~\ref{thm::RE-oracle},
Lemmas~\ref{lemma:threshold-RE} and~\ref{prop:MSE-missing-orig};
cf. Section~\ref{sec:proof-TH-main}.
The full proof of Theorem~\ref{thm::RE-oracle} appears in 
the supplementary Section~\ref{sec::LassoProof}, which yields the
following:  cf. \eqref{eq::part-norms},
\ben 
\label{eq::D0-define-orig}
D_0 & = & D \vee  \sqrt{2} \big(d_0 K^2(s_0, 4) +
  K(s_0, 3) \sqrt{\Lambda_{\max}(s- s_0)}  +  2 d_0 K^2(s_0, 3) \big), \\
\label{eq::D0-prime}
D'_0 & = & D \vee (d_0 K^2(s_0,  4)) \vee (K(s_0, 3)
  \sqrt{\Lambda_{\max}(s- s_0)}  +  3 d_0 K^2(s_0, 3)),\\
\label{eq::defineD}
&& \text{ where } \; D = 
\sqrt{\frac{\Lambda_{\max}(s -
    s_0)}{\Lambda_{\min}(2s_0)}}
\big(1+ \frac{3 \ell(s_0) \sqrt{\Lambda_{\max}(s -s_0)}}{d_0}
\big) \;\text{ and  } \\
\label{eq::definell}
&& \ell(s_0) = ({\theta_{s_0, 2s_0} }/{\sqrt{\Lambda_{\min}(2s_0)}})
\wedge \sqrt{\Lambda_{\max}(s_0) }, \\
\label{eq::D1-define}
D_1 & = &
d_0 [\big(\frac{\Lambda_{\max}(s- s_0)}{d^2_0} + \frac{9}{4}  K^2(s_0, 
 3)\big) \vee 4 K(s_0, 4)^2 \vee \frac{3 \Lambda_{\max}(s- s_0)
 }{d^2_0} ], \\
\label{eq::D2-define}
D_2 & = &
d_0 [\big(\frac{\Lambda_{\max}(s- s_0)}{d^2_0} + 4 K^2(s_0,
  3)\big) \vee 5 K(s_0, 4)^2\vee \frac{4 \Lambda_{\max}(s- s_0)
  }{d^2_0} ], \\
\label{eq::D3-define}
&& \text{ and } \quad
D_3   =  \sqrt{\Lambda_{\max}(s-s_0)}  + d_0 K(s_0, 4)/2 + d_0 K(s_0, 3).
\een
Hence
\bens 
D'_2  & = & (\inv{d_0}\Lambda_{\max}(s- s_0) + 4 K(s_0, 3)^2 d_0) \vee (\frac{4}{d_0} \Lambda_{\max}(s-s_0))  \vee 5 d_0 K(s_0, 4)^2 \\
& \leq & 5 \big(\Lambda_{\max}(s- s_0)/d_0 \vee (d_0 K(s_0, 4)^2)\big).
\eens
The bounds in~\eqref{eq::orthocauchy}
and~\eqref{eq::admissible2} ensure that $ \theta_{s_0, 2s_0} < \infty$
under~\eqref{eq::eigen-max} and $1/\Lambda_{\min}(2s_0) < 2 K^2(s_0,  1) \le  2   K^2(s_0,  4) < \infty$ under  $\RE(s_0, 4, X)$.
In fact, we can obtain an upper bound on $\theta_{s_0, 2s_0}$ in two ways.
Let disjoint sets $J, J' \subset [p]$ satisfy
$\abs{J} \le s_0$ and $\abs{J'} \le 2s_0$ and vectors $\upsilon, \upsilon'$ satisfy
$\twonorm{\upsilon} =\twonorm{\upsilon'} =1$.
We have by the Cauchy–Schwarz inequality and the parallelogram identity:
$${\abs{{\ip{X_{J} \upsilon, X_{J'} {\upsilon'}}}}} \le \twonorm{X_{J} \upsilon}
\twonorm{X_{J'} {\upsilon'} } \le n\sqrt{\Lambda_{\max}(s_0)
  \Lambda_{\max}(2s_0) } \le n \Lambda_{\max}(2s_0)$$
and moreover, cf. the proof of Lemma~\ref{lemma:parallel},
\ben
&& {\abs{{\ip{X_{J} \upsilon, X_{J'} {\upsilon'}}}}}/{n} 
\nonumber
\le  (\Lambda_{\max}(3s_0) -\Lambda_{\min}(3s_0) )/2,
\text{ and  hence},\\
\label{eq::orthocauchy}
&& \theta_{s_0, 2s_0} \le  \sqrt{\Lambda_{\max}(s_0) 
\Lambda_{\max}(2s_0) } \wedge \big(\Lambda_{\max}(3s_0) -
\Lambda_{\min}(3s_0) \big)/2 < \infty.
\een
To help build intuition, we first state in Lemma~\ref{prop:MSE-missing} 
a general result on  the $\ell_2$ loss for the OLS estimator,
when a subset of relevant variables is missing from the fixed model $I$.
Lemma~\ref{prop:MSE-missing} is also an important technical 
contribution of this paper,  which may 
be of independent interest.
The assumption on $I$ being fixed is then 
relaxed in Lemma~\ref{prop:MSE-missing-orig}.
Given $\hat\beta_{I} = (X_I^T X_{I})^{-1} X_{I}^T Y$, let $\hat\beta^{\OLS}(I)$ be the 
$0$-extended version of $\hat\beta_{I}$ such that $\hat\beta^{\OLS}_{I^c} = 0$ and 
$\hat\beta^{\OLS}_I = \hat\beta_{I}$ in both
lemmas.
Let $c_0, c >0$ be some absolute constants.
\begin{lemma}{\textnormal{({\bf OLS estimator with missing variables})}}
\label{prop:MSE-missing}
Suppose sparse eigenvalue conditions~\eqref{eq::eigen-max} and~\eqref{eq::eigen-admissible-s} hold. 
Given a deterministic set $I \subset [p]$, set $\drop := [p] \setminus I$ and $S_\drop = \drop \cap S =
S \setminus I$. Let $\size{I} = m \le (c_0 s_0) \wedge s$.
Suppose $|I \cup S| \leq 2s$.
Let $\hat\beta_{I} = (X_I^T X_{I})^{-1} X_{I}^T Y$.
Then, with probability at least $1 - 2\exp(-3m/64)$, 
\bens
\twonorm{\hat{\beta}^{\OLS}(I) - \beta}^2 \leq 
\big(\frac{2 \theta^2_{|I|, |\dropS|}}{\Lambda_{\min}^2(m)} +
  1\big) \twonorm{\beta_{\drop}}^2 + \frac{3 m \sigma^2}{n \Lambda_{\min}(m)}.
  \eens
\end{lemma}

\begin{remark}
Lemmas~\ref{prop:MSE-missing} and~\ref{prop:MSE-missing-orig} apply to $X$ so long as sparse 
eigenvalue conditions~\eqref{eq::eigen-max} and~\eqref{eq::eigen-admissible-s} hold.
As a consequence,
for any subset $I$ such that $|I| \leq 2s$, 
\begin{eqnarray}
\nonumber
\infty > 
\Lambda_{\max}(2s) & \geq &
\Lambda_{\max}(|I|)  \geq        \Lambda_{\max}  \left({X_I^T X_I}/{n}\right) \\
\label{eq::eigen-cond}
                   & \geq & \Lambda_{\min}  \left({X_I^T X_I}/{n}\right)
                            \geq \Lambda_{\min}(|I|)
                            \geq \Lambda_{\min}(2s) >0.                          
\end{eqnarray}
For disjoint sets $I$ and $S_{\drop} = S \setminus I$ as in
Lemmas~\ref{prop:MSE-missing} and~\ref{prop:MSE-missing-orig}, we have
\ben
\label{eq::model2s}
|I| + |S_{\drop}| = |I \cup S| \leq 2s, \; \text{and}\;
\theta_{|I|, |S_{\drop}|} \le {(\Lambda_{\max}(2s) -
  \Lambda_{\min}(2s))}/{2},
\een
where it is understood that $\Lambda_{\min}(2s) = 0$ is also 
permitted; cf. Lemma~\ref{lemma:parallel} and its proof
in the supplementary Section~\ref{sec::proofparallel}.
\end{remark}
Lemma~\ref{prop:MSE-missing}
implies that even if we miss some 
columns of $X$ in $S$, we can still hope to get the $\ell_2$ loss
bounded as in Theorem~\ref{thm::RE-oracle-main}
so long as $\twonorm{\beta_{\drop}}$ is bounded by
$O_P(\lambda \sigma \sqrt{s_0})$ while $\abs{I}$ is sufficiently small.
Both conditions are guaranteed to hold by our choices of the thresholding parameters
as shown in Lemma~\ref{lemma:threshold-RE}.
Although the tight analysis of Lemma~\ref{prop:MSE-missing} depends on
the fact that the selection set $I$ is deterministic, a simple
variation of the statement makes it work well with the thresholded
estimators as considered in the present paper, with $\ell_2$ error
bounded essentially at the same order of magnitude as
in~\eqref{eq::log-MSE} so long as $\size{I} =O(s_0)$ and $|I \cup S| \leq 2s$.
Lemma~\ref{prop:MSE-missing-orig} is presented by~\cite{Zhou09th}.
\begin{lemma}~\citep{Zhou09th}
  \label{prop:MSE-missing-orig}
Suppose sparse eigenvalue conditions~\eqref{eq::eigen-max} and~\eqref{eq::eigen-admissible-s}
hold. Given an arbitrary set $I \subset [p]$, possibly random, set
$\drop := [p] \setminus I$ and $S_\drop = \drop \cap S$.
Suppose on event $\T_a$, $\size{I} \le (c_0 s_0) \wedge s$ and $|I \cup S| \leq 2s$. Then,
for $\hat\beta_{I} = (X_I^T X_{I})^{-1} X_{I}^T Y$ and
$\hat\beta^{\OLS}(I)$ as defined in Lemma~\ref{prop:MSE-missing}, it holds on
$\T_a$ that
\ben
\twonorm{\hat{\beta}^{\OLS}(I) - \beta}^2
& \leq & 
\label{eq::plossintro}
\big(\frac{2 \theta^2_{|I|, |\dropS|}}{\Lambda^2_{\min}(|I|) } + 1\big)\twonorm{\beta_{\drop}}^2
+ \frac{2|I| (1+a) \sigma^2 \lambda^2}{\Lambda^2_{\min}(|I|)}.
\een
\end{lemma}

The proofs for Lemmas~\ref{prop:MSE-missing}
and~\ref{prop:MSE-missing-orig} are deferred to 
Sections~\ref{sec::OLSproof-orig} and~\ref{sec::OLSmissing2}. 
We now state Lemma~\ref{lemma:parallel}, which follows 
from~\cite{CT05}.
\begin{lemma}\textnormal{[Lemma 1.2~\cite{CT05}]}
  \label{lemma:parallel}
Suppose sparse eigenvalue conditions~\eqref{eq::eigen-max} and~\eqref{eq::eigen-admissible-s}
hold.
The following statements hold:
(a) For all disjoint sets $I,S_{\drop} \subseteq [p]$ 
of cardinality $|S_{\drop}| < s$ and $|I| + |S_{\drop}| \leq 2s$, 
$$\theta_{|I|, |S_{\drop}|} \leq (\Lambda_{\max}(2s) -
\Lambda_{\min}(2s))/{2};$$
(b) Without the lower sparse eigenvalue (LSE)
condition~\eqref{eq::eigen-admissible-s}, we still obtain
\ben
\label{eq::sparseupper}
\theta_{|I|, |S_{\drop}|} \leq \Lambda_{\max}(s) \wedge 
(\Lambda_{\max}(2s)/2) < \infty \; \text{ if } \; \abs{I} \vee |S_{\drop}| \le
s,
\een
under the USE condition~\eqref{eq::eigen-max}, following the same arguments leading to~\eqref{eq::orthocauchy}.
\end{lemma}

\begin{remark}
  \label{rem::RESPA}
(I.) In Lemma~\ref{prop:MSE-missing-orig}, we have
$\abs{I} \vee |S_{\drop}| \le s$. Hence, the LSE condition~\eqref{eq::eigen-admissible-s} can be replaced with the
following relaxed lower sparse  
eigenvalue condition:
\ben 
\label{eq::bound-on-I}
\Lambda_{\min}((cs_0) \wedge s) > 0, \text{ in case }
\;  |I| \leq (c s_0) \wedge s \text{ for some $c \ge 2$}
\een
so that \eqref{eq::plossintro} holds on $\T_a$, with $\theta_{|I|, |\dropS|} \leq \Lambda_{\max}(s) \wedge 
(\Lambda_{\max}(2s)/2)$ as in~\eqref{eq::sparseupper}. \\
(II.) We note that if $\RE(s_0, k_0, X)$ as defined
in~\eqref{eq::admissible} is satisfied with $k_0 \geq 1$ and $1\le s_0 < p$
then
$1/{\Lambda_{\min}(2s_0) }\le 2 K^2(s_0, 1)$. Consider $2s_0$ sparse
vector $\up$. Let $T_0$ denote  the locations of the  $s_0$ largest
coefficients of $\up$ in absolute values.  Then $\twonorm{\up}^2 \le 2
\twonorm{\upsilon_{T_0}}^2$, since $\twonorm{\up_{T_0^c}}  \leq \twonorm{\up_{T_0}}.$
Thus we have for any $2s_0$-sparse vector $\up$, by $\RE(s_0, 4, X)$~\eqref{eq::admissible},
\bens
\frac{\twonorm{X \upsilon}^2}{n \twonorm{\upsilon}^2}  \ge 
\frac{\twonorm{X \upsilon}^2}{2 n \twonorm{\upsilon_{T_0}}^2}  \ge
1/{(2K^2(s_0, 1))}, \text{ where }  \norm{\up_{T_0^c}}_1  \leq \norm{\up_{T_0}}_1.
\eens
Hence \eqref{eq::bound-on-I} holds for $c =2$, when $2s_0 < s$, since
\ben
\label{eq::admissible2}
\Lambda_{\min}(2s_0) \; \stackrel{\triangle}{=} \;
\min_{\upsilon \not=0; 2s_0-\text{sparse}} \; 
{\twonorm{X \upsilon}^2}/{(n \twonorm{\upsilon}^2)}
\ge  1/(2K^2(s_0, 1)).
\een
\end{remark}
\begin{remark}
One notion of the incoherence condition which has been formulated in the sparse 
reconstruction literature bears the name of restricted isometry 
property (RIP)~\citep{CT05,CT07}.
For each integer $s = 1, 2, \ldots$ such that $s < p$, 
the $s$-{\em restricted isometry constant}  $\delta_s$ of $X$ is defined to be the smallest quantity such that
\ben
\label{eq::RIP-define}
(1 - \delta_s) \twonorm{\upsilon}^2 \leq \twonorm{X_T \upsilon}^2/n \leq 
(1 + \delta_s) \twonorm{\upsilon}^2
\een
for all $T \subset [p]$ with $|T| \leq s$ and coefficients 
sequences $(\upsilon_j)_{j \in T}$~\citep{CT05}; Hence the upper and lower $s$-sparse eigenvalues 
of design matrix $X$ satisfy $1+\delta_s \ge \Lambda_{\max}(s)
\geq \Lambda_{\min}(s) > 1-\delta_s$.
\end{remark}

\subsection{Discussions}
\label{sec::threshold}
The tight analysis in Theorem~\ref{thm::RE-oracle-main} for the Thresholded Lasso estimator 
is motivated by the sparse oracle inequalities on the Gauss-Dantzig 
selector under the UUP, which is originally shown to hold in a conference 
paper by~\cite{Zhou09th}.
Previously, it has been shown that~\eqref{eq::log-MSE} holds with high 
probability for the Dantzig selector under the condition of a Uniform Uncertainty Principle (UUP), 
where the UUP states that for all $s$-sparse sets $J$, the columns of 
$X$ corresponding to $J$ are almost  orthogonal~\citep{CT07}.
\begin{definition}
\textnormal{(\bf{A Uniform Uncertainly Principle})}
\label{def:CT-cond}
For some integer $1 \leq s < n/3$, assume $\delta_{2s} + \theta_{s,
  2s} < 1-\tau$ for some $\tau>0$.
\end{definition}
Then under the settings of Lemma~\ref{prop:MSE-missing-orig}, we have 
$\theta_{|I|, |\dropS|} \leq \delta_{2s} < 1$, where $|I| +
|S_{\drop}| \leq 2s$; cf. \eqref{eq::model2s}. 
The sparse eigenvalue conditions in the present work are considerably relaxed from the
incoherence condition (UUP) in Definition~\ref{def:CT-cond}:
The UUP condition ensures that the
$\RE(s, 1, X)$ condition as in \eqref{eq::admissible} holds with
$$K(s, 1) =
{\sqrt{\Lambda_{\min}(2s)}}/{(\Lambda_{\min}(2s) - \theta_{s, 2s})}
\leq {\sqrt{\Lambda_{\min}(2s)}}/{\tau};$$
cf.~\cite{BRT09}.
Besides the sparse eigenvalue conditions,
Theorem~\ref{thm::RE-oracle-main} requires $\RE(s_0, 4, X)$ to be satisfied, which
depends on the essential sparsity $s_0$ rather than $s = \abs{\supp(\beta)}$;
cf. Section~\ref{sec::proofsketchRE} for detailed discussions.
In the supplementary  Theorem~\ref{thm:ideal-MSE-prelude}, we show the 
corresponding result for the Gauss-Dantzig selector under the UUP for 
completeness.

In summary, Lemmas~\ref{lemma:threshold-RE} and~\ref{prop:MSE-missing-orig}
ensure that the general thresholding rules with threshold level at about 
$\lambda \sigma$ achieve the following property:
although we cannot guarantee the presence of variables indexed by
$S_R = \{j: |\beta_j| < \sigma \sqrt{\log p/n} \}$ to be 
included due to their lack of strength, we will include in $I$ most variables in
$S \setminus S_R$ such that the OLS estimator based on model $I$ 
achieves the oracle bound~\eqref{eq::log-MSE}.
This goal is accomplished despite some variables from the support set 
$S$ are missing from the model $I$, since their overall $\ell_2$-norm 
$\twonorm{\beta_{\drop}}$ is bounded in
\eqref{eq::off-beta-norm-bound}.

As mentioned, Proposition~\ref{PROP:COUNTING-S0} (by setting $c' = 1$)
shows that the number of variables in $\beta$ that are larger than or equal to 
$\sqrt{\log p/n} \sigma$ in magnitude is bounded by $2s_0$.
In hindsight, it is clear that we wish to retain most of them by
keeping $2s_0$ variables in the model $I$.
Indeed, suppose $D_1 s_0 <( c_0 s_0) \wedge (s-s_0)$.
Then, by choosing $t_0 \asymp \lambda \sigma$ on event $\T_a$, we are
guaranteed to obtain
under the settings of Theorem~\ref{thm::RE-oracle-main} (and Lemma~\ref{lemma:threshold-RE}), 
\bens
|I| & \leq  & s_0 (1 + D_1)   < s, \; \; |I \cup S|   \leq  2 s \; \text { and } 
\shtwonorm{\hat\beta - \beta}^2 \leq D_4^2 s_0 \lambda^2 \sigma^2.
\eens
Here, we assume that $D_1$ as in~\eqref{eq::D1-define} will grow only
mildly with the parameters $s_0, s$, under
conditions~\eqref{eq::admissible} and~\eqref{eq::eigen-max},
and it is not necessary to set $C_4 > D_1$.
The set of missing variables in $\drop$ is the price we pay in order
to obtain a sparse model when some coordinates in the support
$\supp(\beta)$ are well
below $\sigma \sqrt{\log  p/n}$. Note that when we allow the
model size to increase by lowering $t_0$,
the variance term $\propto |I| \lambda^2/\Lambda_{\min}(|I|)$ becomes
correspondingly larger.
Since the larger model $I$ may not include more true variables,
the size of $\dropS = S \setminus I$ may remain invariant;
If so, the overall  interaction term $\theta
\twonorm{\beta_{\drop}}/\Lambda_{\min}(|I|)$ can still increase due to
the increased orthogonality coefficient $\theta:=\theta_{|I|, |S_{\drop}|}$, even though
$\twonorm{\beta_{\drop}}$ is a non-increasing function of the model
size.

This argument favors the selection of a small (yet sufficient)
model as stated in Theorem~\ref{thm::RE-oracle-main}, rather than
blindly including extraneous variables in the set $I$.
We mention in passing that by setting an upper bound on the desired
model size $\abs{I} \le 2 s_0$, we are able to make some interesting
connections between the thresholded estimators as studied in the
present paper and the $\ell_0$ penalized least squares estimators. In
particular, we show that the prediction error,  $\shtwonorm{X \beta_{\init} - X \beta}$, and a complexity-based
penalty term $\sigma \sqrt{\abs{I}\log p}$ on the chosen model $I$ are
both bounded by $O_P(\sigma \sqrt{s_0 \log p})$ in case $|I|
\asymp 2 s_0$ for the thresholded estimators by~\cite{Zhou10}.

\subsection{Background and related work} 
\label{sec:linear-sparsity} 
In this section, we briefly discuss related work.
Lasso and the Dantzig selector are both computational efficient and shown with 
provable nice statistical properties;
see for example~\cite{MB06,GR04,Wai09,ZY06,CT07,BTW07c,CP09,Kol09,Kol07,ZH08,MY09,BRT09}. We
refer to the books for a comprehensive survey of related
results~\cite{BG11,Wain19}.

Prior to our work, a similar two-step procedure,  namely, 
the Gauss-Dantzig selector, has been proposed and empirically 
studied by~\cite{CT07}.
This paper builds upon the methodology originally developed in a 
conference paper by the present author~\citep{Zhou09th}.
Inspired by~\cite{CT07},~\cite{Zhou09th} obtains oracle bounds in the spirit of
Theorem~\ref{thm::RE-oracle-main} for the Gauss-Dantzig selector,
under the stronger restricted isometry type of conditions as originally proposed
by~\cite{CT07}; cf. the supplementary Theorem~\ref{thm:ideal-MSE-prelude}.

Under variants of the RIP conditions, the exact recovery or approximate reconstruction of a
sparse $\beta$ using the basis pursuit program~\citep{Chen:Dono:Saun:1998} has been shown in a series of
results~\citep{Donoho:cs,CT05,CT07,RV08}.
We refer to the book
by~\cite{Wain19}  and the paper by~\cite{GB09} for a complete exposition.
The sparse recovery problem under arbitrary noise is also well studied, 
see for example~\cite{NV09} and~\cite{NT08}, where 
they require $s$ to be part of  the input.
See~\cite{FG94},\cite{Dono:John:1995},~\cite{ZRXB11},~\cite{Horns19}
for further references and applications of the essential sparsity.

For the Lasso,~\cite{MY09} has also shown in theoretical analysis that 
thresholding is effective in obtaining a two-step estimator $\hat\beta$ 
that is consistent in its support with $\beta$ when $\beta_{\min}$ is 
sufficiently large.
A weakening of 
the incoherence condition by~\cite{MY09} is still sufficient 
for~\eqref{eq::admissible} to hold~\citep{BRT09}.
See also~\cite{Mei07},~\cite{ZouLi08} and~\cite{Zhang09}.
A more general framework on multi-step variable selection was 
explored by~\cite{WR08}. They control the probability of false positives at 
the price of false negatives, similar to what we aim for here; their analysis is 
constrained to the case when $s$ is a constant.

\noindent{\bf Subsequent development.}
This choice of the threshold parameter identified in~\cite{Zhou09th,Zhou10}
and the current paper has deep connection with the classic and current
literature on model selection~\citep{FG94,Dono:John:1995,WR08,Zhang10,Wain19}.
~\cite{Zhang10} proves the minimax concave penalty (MCP) procedure is selection consistent 
under a sparse Riesz condition and an information requirement in the 
sense of~\eqref{eq::betamin}.
Nonconvexity of the minimization problem causes computational and
analytical difficulties; Compared to the elegant yet complex method in MCP, 
the Thresholded Lasso procedures~\citep{Zhou09th,Zhou10} provide 
a much simpler framework, which is desirable from the practical point of 
view, with overall good performance. This is confirmed in a subsequent 
study by~\cite{Wang13}.

While the focus of the present paper 
is on variable selection and oracle inequalities for the $\ell_2$ loss,  
prediction errors are also explicitly derived by~\cite{Zhou10}. 
~\cite{GBZ11} revisit the adaptive Lasso~\citep{Zou06,HMZ08,ZGB09} as
well as the Thresholded Lasso with refitting~\citep{Zhou09th,Zhou10}, 
in a high-dimensional linear model, and study prediction error and
bound the  number of false positive selections.
We refer to~\cite{FG94},~\cite{Barr:Birg:Mass:1999},~\cite{BM01} and~\cite{Shen12} for 
related work on complexity regularization criteria.
In a subsequent work, \cite{ZRXB11} develop
error bounds based on an earlier version of the present paper and
applied these to obtain fast rates of convergence for covariance
estimation based on a multivariate Gaussian graphical model. There we
show comprehensive numerical results involving
cross-validation to choose penalty $\lambda_n$ and thresholding
$t_0$ parameters.
We mention that a series of recent papers~\citep{RWY10,RZ13,LW12,RZ17,Zhou22RE} show that $\RE$ condition holds for a broader class of random matrices with
complex row/columnwise (or both) dependencies
once the sample size is sufficiently large.

\section{Proof sketch for the main result}
\label{sec::proofsketchRE}
We will now describe the main ideas of our 
analysis in this section. Combining Theorem~\ref{thm::RE-oracle} with
Lemmas~\ref{prop:MSE-missing} and~\ref{lemma:threshold-RE} allows us
to prove Theorem~\ref{thm::RE-oracle-main}, which we will elaborate in more details in Section~\ref{sec:proof-TH-main}.
For now, we highlight the important differences between our results 
and a previous result on the Lasso by~\cite{BRT09} (cf. Theorem 7.2), 
which we refer to as the BRT results.
While a bound of $O_P(\lambda \sigma \sqrt{s})$ on the $\ell_2$ loss 
as obtained by~\cite{BRT09} makes sense when all signals are strong, 
significant improvements are needed for the general case where
$\beta_{\min}$ is not bounded from below.

In the present work, the goal is to investigate sufficient
conditions under which we could achieve a bound of
$O_P(\lambda \sigma \sqrt{s_0})$ on the $\ell_2$ loss for both the
Lasso and the Thresholded Lasso. Given such error bound for the Lasso, thresholding of an initial estimator
$\beta_{\init}$ at the level of  $\asymp \sigma \sqrt{2 \log p/n}$
will select nearly the best subset  of variables in the spirit of
Theorem~\ref{thm::RE-oracle-main}.
Some  more comments.

\noindent{\bf (a)}
As stated in Theorem~\ref{thm::RE-oracle-main}, we use $\RE(s_0, 4, X)$, 
for which we fix the {\it sparsity level}  at $s_0$ and $k_0 = 4$, and
sparse eigenvalue conditions~\eqref{eq::eigen-admissible-s}
and~\eqref{eq::eigen-max}. 
While the constants in association with the BRT results depend on  $K^2(s, 3)$,
the constants in association with the Lasso and the Thresholded Lasso 
crucially depend on $K^2(s_0, 4)$, $\Lambda_{\max}(2s)$,
$\Lambda_{\min}(2s_0)$,  and $\theta_{s_0, 2s_0}$
(cf.~\eqref{eq::orthocauchy}).

\noindent{\bf (b)}
We note that the lower sparse eigenvalue condition
$\Lambda_{\min}(2s) >0$~\eqref{eq::eigen-admissible-s} is implied by,
and hence is weaker than the $\RE(s, 3, X)$ condition.
Moreover, it is possible to prove Theorem~\ref{thm::RE-oracle-main}
even if we leave condition~\eqref{eq::eigen-admissible-s} out.
In particular, we note that so long as $|I| \leq 2 s_0$, then
$\RE(s_0, 4, X)$ already implies that~\eqref{eq::bound-on-I}
holds for $c = 2$~\citep{BRT09};
cf. Remark~\ref{rem::RESPA}.
We will not pursue such optimizations in the present work.

\noindent{\bf (c)}
We note that in the $\RE(s, 3, X)$ condition as required by~\cite{BRT09} to 
achieve the $\ell_2$ loss of  $O_P(\lambda \sigma \sqrt{s})$: while $k_0 = 3$ is chosen, 
they fix sparsity at $s$ instead of $s_0$, which is not ideal when $s_0$ is much smaller than $s$.
We emphasize that in the $\RE(s_0, 4, X)$ condition that we impose, 
$k_0 = 4$ is rather arbitrarily chosen; in principle, 
it can be replaced by any number that is strictly larger than $3$.
In the context of compressed sensing, $\RE$ conditions can also be 
taken as a way to guarantee recovery for anisotropic 
measurements~\citep{RZ13}.
Results by~\cite{RZ13} reveal that for $\RE$ conditions with a smaller
$s_0$, we need correspondingly smaller sample size $n$ in order for
the random design matrix $X$ of dimension $n \times p$ to satisfy such
a condition, when the independent row vectors of $X_i, i=1,\ldots, n$ have covariance
$\Sigma(X_i) = \E X_i \otimes X_i = \E X_i X_i^T$ satisfying
$\RE(s_0, (1+ \ve) k_0, \Sigma^{1/2})$ condition in the sense that
for any $\ve > 0$, for all $\upsilon \not=0$,
$$\inv{K(s_0, (1+ \ve) k_0, \Sigma^{1/2})} \stackrel{\triangle}{=}
\min_{\stackrel{J_0 \subseteq [p]}{|J_0| \leq s_0}}
\min_{\norm{\upsilon_{J_0^c}}_1 \leq k_0 \norm{\upsilon_{J_0}}_1}
\; \;  {\norm{\Sigma^{1/2} \upsilon}_2}/{\norm{\upsilon_{J_0}}_2} > 0.$$

\noindent{\bf (d)}
We impose an explicit upper bound on 
$\Lambda_{\max}(2s)$, which is absent from the paper by~\cite{BRT09}, 
in order to obtain the tighter bounds in the present work for both the
Lasso and the Thresholded Lasso.
This condition is required by our OLS refitting procedure as stated in 
Lemmas~\ref{prop:MSE-missing-orig} and~\ref{lemma:parallel}.
This is consistent with the fact that known oracle inequalities for
the Dantzig and Gauss-Dantzig selectors are proved under the UUP
which impose tighter upper and lower sparse eigenvalue bounds in the
sense that $\theta_{s, 2s}+\delta_{2s}<1$. See also~\cite{ZH08} and
the exposition in Section~\ref{sec::tightness}.

\subsection{Proof sketch of Theorems~\ref{thm::RE-oracle-main} and~\ref{thm::RE-oracle}}
\label{sec::relaxed}
We specify the initial Lasso estimator $\beta_{\init}$, and the parameters $\lambda_n$, 
$t_0$ in Theorem~\ref{thm::RE-oracle} and Lemma~\ref{lemma:threshold-RE}  respectively. 
Throughout this section, let $h = \beta_{\init} - \beta^{\ext}({T_0})$ be as in
Lemma~\ref{lemma:threshold-RE},   where $T_0$  denotes the largest
$s_0$ coordinates of $\beta$ in absolute values.
Let $T_1$ be the $s_0$ largest positions of $h$ outside of $T_0$. 
 Denote by $T_{01} = T_0 \cup T_1$. 
In this section, we elaborate on  the $\ell_p, p=1, 2$, loss on $h_{T_0}$ and $h_{T_0^c}$,
and their implications on variable selection using thresholding.
As mentioned, improving the bounds on each component of $h$ will
result in a tighter upper bound on controlling the bias;
cf. Lemma~\ref{lemma:threshold-RE}.
Specifically, we will show in the supplementary
Section~\ref{sec::LassoProof} that on $\T_a$,
\eqref{eq::part-norms} holds for $h = \beta_{\init} -
\beta^{\ext}({T_0})$, where $D'_0$ and $D_1$ are as in~\eqref{eq::D0-prime}
and~\eqref{eq::D1-define}. 

These bounds ensure that under the $\RE$ and  
sparse eigenvalue conditions, both the $\ell_2$ loss on the set $T_0$
of significant coefficients, namely, $\norm{h_{T_0}}_2 =\norm{(\beta_{\init} - \beta)_{T_0}}_{2}$,
and the $\ell_1$ (and the  $\ell_2$) norm of the estimated
coefficients on $T_0^c$, namely,  $\norm{h_{T_0^c}}_1 =
\norm{\beta_{\init,T_0^c}}_1$  are tightly bounded with respect to
$\abs{T_0}$ for the Lasso; hence we obtain an oracle inequality on the
$\ell_2$ loss in the sense of~\eqref{eq::log-MSE}.

First, by definition of $h$ and $T_0$, we aim to keep variables in 
$T_0$, while for variables not in $T_{0}$, we may trim these 
off. It is also clear by Lemma~\ref{lemma:threshold-RE} that we cannot
cut too many ``significant'' variables in $T_0$ by following the
thresholding rules in our proposal; for example, for those that are
$\geq \lambda \sigma \sqrt{s_0}$,  we can cut at most a constant
number of them.
Now, what do we do with those in $T_1$? The fate of these variables is
pretty much up to the choice of the threshold for a given
$\beta_{\init}$, knowing these are the largest in magnitude in $h$
(and $\beta_{\init}$)  outside of $T_0$ and hence most likely to be
included in model $I$.
Moreover, even if we were able to retain all variables in $T_0$, we will
include at least some variables in $T_1$ when the selection set has
size $\abs{I} > s_0$; in fact, we will include nearly all variables in $T_1$
when the set $I$ has size close to $2s_0$. In
particular,
\bens 
\size{I \cap T^c_0} & \le &
{\norm{h_{T_0^c}}_1}/{t_0}  \le 
D_1 \lambda_n s_0/(f_0 \lambda_n) \le( D_1/f_0) s_0 \quad \text{in 
  case} \quad t_0 =  f_0 \lambda_n
\eens
for some $f_0 > 0$ by \eqref{eq::part-norms}.
Before we continue, we state Lemma~\ref{lemma::h01-bound-CT}, which is
the same (up to normalization) as Lemma 3.1~\citep{CT07}, to illuminate
the roles of sets $T_0, T_1$ in the overall bounds on $\twonorm{h}$.
We note that in their original statement, the UUP condition is assumed; 
a careful examination of their proof shows that it is a sufficient 
but not necessary  condition; indeed we only need to assume that 
sparse eigenvalues are bounded, namely, 
$\Lambda_{\min}(2s_0) > 0$ and $\Lambda_{\max}(2s_0) < \infty$. 
We also state Lemma~\ref{lemma::cone}.
\begin{lemma}
  \label{lemma::h01-bound-CT}
Suppose $\Lambda_{\min}(2s_0) > 0$ and $\Lambda_{\max}(2s_0) < \infty$.
Then
\ben
\label{eq::simpleT01orig}
\twonorm{h_{T_{01}}} & \leq &  
\inv{\sqrt{\Lambda_{\min}(2s_0)}}
\twonorm{X h}/ \sqrt{n} + \frac{\theta_{s_0, 2s_0}}{\Lambda_{\min}(2s_0)}
\norm{h_{T_0^c}}_1/ \sqrt{s_0}, \\
\label{eq::RE-T0c-2-bound}
\twonorm{h_{T_{01}^c}}^2 & \leq &
\norm{h_{T_0^c}}_1^2 \sum_{k \geq s_0 + 1} 1/k^2  \; \leq \;
\norm{h_{T_0^c}}_1^2/s_0 \; \text { and thus} \\
\nonumber
\twonorm{h}^2 & \leq & \twonorm{h_{T_{01}}}^2 + s_0^{-1}
\norm{h_{T_0^c}}^2_1.
\een
\end{lemma}

\begin{lemma}{\textnormal{~\citep{RZ13}}}
\label{lemma::cone} 
Denote by
\beq
\label{eq::cone}
\Cone(s_0, k_0) :=
\left\{x \in \R^p, \exists {J \subseteq [p], |J| = s_0} \text{ s.t.}
  \norm{x_{J^c}}_1 \leq k_0 \norm{x_{J}}_1 \right\}.
\eeq
Let $T_0$ denote the locations of the largest coefficients of $x$ in
absolute values.
Then $\twonorm{x} \le \sqrt{1+ k_0}\twonorm{x_{T_0}}$ for $x \in \Cone(s_0, k_0)$.
\end{lemma}
Now we provide a proof sketch for Theorem~\ref{thm::RE-oracle}.
For $u =(u_1, \ldots, u_n) \in \R^n$, define the empirical norm of $u$ by 
$\norm{u}_n^2 = \twonorm{u}^2/n.$
The proof of the original Theorem~5.1 in~\cite{Zhou10} draws upon
techniques from a concurrent work by~\cite{GBZ11} and uses
an elementary inequality for the Lasso; cf. the supplementary~\eqref{eq::missingone}.
Denote by $\beta_0 :=\beta^{\ext}({T_0})$.
The current proof replaces the supplementary~\eqref{eq::missingone}
for the Lasso with the following updated inequality:
\ben
\nonumber
\lefteqn{ \norm{X (\beta_{\init} - \beta)}^2_n  + \norm{X (\beta_{\init} -\beta_0 )}^2_n
  \le   \norm{X (\beta -\beta_0)}^2_n} \\
\label{eq::oracleD}
&& + \frac{2}{n} \e^T  X(\beta_{\init} - \beta_0) + 2 \lambda_n
(\norm{\beta_0}_1 -\shonenorm{\beta_{\init}}).
\een
See Eq(20)~\citep{DHL17}, where we set $\bar\beta = \beta_0
=\beta^{\ext}({T_0})$ and $\bar{\delta} = \beta_{\init} -\beta_0  = h$.

Denote by $\delta = \beta_{\init} - \beta$.
Similar to the proof in~\cite{Zhou10}, we have a 
deterministic proof on event $\T_a$, except that now on $\T_a$,
\ben
\label{eq::new-oracle}
\norm{X \delta}^2_{n} + \norm{X h}^2_{n}
+ \lambda_n \norm{h_{T_0^c}}_1 
\leq \norm{X (\beta-\beta_0)}^2_{n} + 3 \lambda_n  \norm{h_{T_0}}_1.
\een
We differentiate between three cases under event $\T_a$.
\begin{enumerate}
\item
In the first case, suppose 
\ben
\label{eq::DHDomi}
&& \norm{X\delta}^2_n  + \norm{X h}^2_{n}  \ge
\norm{X (\beta -  \beta_0)}_n^2; \\
\nonumber
\text{Then,} &&
\norm{h_{T_0^c}}_1   \le  3 \norm{h_{T_0}}_1, \text{ and hence }  \; \; h \in \Cone(s_0, 3). 
\een
We will show that on event $\T_a$, for $\lambda_n  = d_0  \lambda \sigma$,
\ben 
\label{eq::HT0case1a}
\twonorm{h_{T_0}} & \le & K(s_0, 3) \norm{X \beta - X \beta_0}_{n}  +
3K^2(s_0, 3) \lambda_n \sqrt{s_0};
\een
Moreover, the bounds on $\norm{X \delta}_n$ and $\norm{h_{T_0^c}}_1$
(and $\norm{h}_1$)  follow from \eqref{eq::newboundsim}:
 \ben 
  \label{eq::newboundsim}
  \norm{X \delta}^2_n
 + \lambda_n \norm{h_{T_0^c}}_1 
 & \leq &\norm{X (\beta-\beta_0)}_n^2
 + (3 K(s_0, 3) \lambda_n \sqrt{s_0}/2)^2,
\een
where $\norm{X (\beta-\beta_0)}_n \le \lambda \sigma \sqrt{s_0
  \Lambda_{\max}(s-s_0)}$.
Then we have
\ben
\nonumber
\norm{h_{T_0^c}}_1
& \le & \norm{ X \beta - X \beta_0 }^2_{n}/( \lambda_n) 
+ (3K/2)^2 \lambda_n  s_0  \le D_{1, a} \lambda \sigma s_0,\\
\label{eq::TailboundMain}
\text{ where}  \; \; D_{1, a} & := &
\Lambda_{\max}(s- s_0)/d_0 + 9 d_0 K(s_0, 3)^2  /{4}.
\een
\item
In the second case, suppose 
\ben 
\label{eq::DHhead}
 \norm{X\delta}^2_n  + \norm{X h}^2_{n}
  \le \lambda_n \norm{h_{T_0}}_1. 
  \een 
  We will show that on event $\T_a$,  $h \in \Cone(s_0, 4)$ and moreover,
\ben 
\label{eq::Xdelta2}
&&  (2 \norm{X \delta}_{n}) \vee  \norm{X h}_{n}
 \leq \lambda_n  \sqrt{s_0} K(s_0, 4),\\
 \nonumber
 && \twonorm{h_{T_0}}
 \leq  K(s_0, 4) \norm{X h}_{n} \le \lambda_n  \sqrt{s_0} K(s_0, 4)^2, \; \text{ and } \;\\
\label{eq::cone4}
&& \norm{h_{T_0^c}}_1  \le  4 \norm{h_{T_0}}_1 \le 4 \lambda_n K(s_0,
4)^2  s_0.
\een
\item
Finally, we consider
\ben 
\label{eq::DHmiddlemain}
\lambda_n
\norm{h_{T_0}}_1 \le \norm{X\delta}^2_n  + \norm{X h}^2_{n}
  \le \norm{X \beta -  X \beta_0}_n^2. 
  \een
Then clearly \eqref{eq::pred-error-gen} also holds with no need to
invoke $\RE(s_0, c_0, X)$ condition (for any $c_0$) at all: 
\bens 
\norm{X\delta}_n^2  \vee \norm{ X h}_{n}^2 
&\leq &
\norm{X \beta - X \beta_0}_n^2 \le 
\Lambda_{\max}(s- s_0)\lambda^2 \sigma^2 s_0 \text{ and } \\
\norm{h_{T_0}}_1 &\leq &
\norm{X \beta - X \beta_0}_n^2/\lambda_n \asymp \lambda_n  s_0.
\eens 
Moreover, we have $\norm{h}_1 \leq 
4 \Lambda_{\max}(s- s_0) \lambda \sigma s_0/(d_0)$ and
$\norm{h_{T_0^c}}_1 \le D_{1, c} \lambda \sigma s_0$,
where $D_{1, c} =3 \Lambda_{\max}(s- s_0)/d_0$, since
\ben
 \label{eq::HT01Case3}
\norm{h_{T_0^c}}_1 
& \leq &
3 \norm{X \beta - X \beta_0}_n^2/\lambda_n  \leq 
3 \Lambda_{\max}(s- s_0) \lambda \sigma s_0/(d_0).
\een
Hence although the vector $h$ may not satisfy the cone constraint, 
both components of $h$ have bounded $\ell_1$ norm of order $\asymp 
\lambda \sigma s_0$.
\end{enumerate}

Using~\eqref{eq::newboundsim},~\eqref{eq::DHhead}, \eqref{eq::Xdelta2}, and 
\eqref{eq::DHmiddlemain}, we conclude that \eqref{eq::pred-error-gen}
holds.
Although the new proof is introduced to get rid of the constant factor $2$ in 
front of  $\norm{X \beta - X_{I_0} \beta_{I_0}}_{n} =\norm{X (\beta -\beta_{0})}_{n}$ in~\eqref{eq::pred-error-gen} in~\cite{Zhou10}, we are not optimizing
the constants in this paper; cf. the proof of the  original
Theorem~5.1 in~\cite{Zhou10}.
For example, the constant $3$ in~\eqref{eq::new-oracle} can be further 
reduced following~\cite{DHL17} using \eqref{eq::oracleD}. 
Moreover, combining
\eqref{eq::HT0case1a},~\eqref{eq::TailboundMain},~\eqref{eq::cone4},
and \eqref{eq::HT01Case3}, we have
the expression of~\eqref{eq::D1-define} for $\norm{h_{T_0^c}}_1$,
namely, $\norm{h_{T_0^c}}_1  \leq  D_1 \lambda \sigma s_0$.
  We demonstrate the tightness of these bounds in 
  Fig.~\ref{fig:lasso-esimate-L1-b2-b0-l2-b1-error}.

Next we bound $\twonorm{h}$ in view of
Lemma~\ref{lemma::h01-bound-CT}.
Lemma~\ref{lemma::HT01case3} also provides the upper bounds on 
$\twonorm{h_{T_{01}}}$, leading to the expression of  $D_0$ as 
in~\eqref{eq::D0-define-orig}.
We combine the upper bounds in Lemma~\ref{lemma::h01-bound-CT}
and the supplementary Lemma~\ref{lemma::T01simple} in the proof of 
Lemma~\ref{lemma::HT01case3}.
We prove Lemmas~\ref{lemma::h01-bound-CT} and~\ref{lemma::HT01case3},
and Theorem~\ref{thm::RE-oracle} in the supplementary
Sections~\ref{sec::proofofT01},~\ref{sec::suppHT01proof},
and~\ref{sec::LassoProof} respectively.
\begin{lemma}
 \label{lemma::HT01case3}
 Under the settings of  Theorem~\ref{thm::RE-oracle}, we have under
 event $\T_a$,
  \ben
  \nonumber
\text{\bf Case 1}: &&
\twonorm{h} \leq 
2 \twonorm{h_{T_{01}}}  \leq
2 D_{0,a} \lambda \sigma \sqrt{s_0} \;  \text{where} \\
\label{eq::TaDa}
D_{0,a}  & \le & \sqrt{2}\big(\sqrt{\Lambda_{\max}(s-s_0)} K(s_0, 3)
+ 3 d_0 K^2(s_0, 3) \big); \\
  \nonumber
\text{\bf Case 2}: && \twonorm{h} \leq 
\sqrt{5} \twonorm{h_{T_{01}}}  \leq 
\sqrt{10} K^2(s_0, 4) d_0 \lambda \sigma \sqrt{s_0};\\
\nonumber 
\text{\bf Case 3}: &&
\twonorm{h_{T_{01}}} 
\leq \inv{\sqrt{\Lambda_{\min}(2s_0)}}
\big(\norm{X h}_{n} +  \frac{\ell(s_0)}{\sqrt{s_0}}\norm{h_{T_0^c}}_1\big)
 \le D \lambda \sigma \sqrt{s_0}, \\
&&
\label{eq::Dlocal}
\text{where} \; D  :=
  \sqrt{\frac{\Lambda_{\max}(s - s_0)}{\Lambda_{\min}(2s_0)}}
  \big(1 + \frac{3 \ell(s_0)}{d_0}  \sqrt{\Lambda_{\max}(s-s_0)}\big)
  \een
  and $\ell(s_0)$ is as defined in~\eqref{eq::definell}.
  Moreover, under event $\T_a$, we have
  for $D_0$ in~\eqref{eq::D0-define-orig} and $D_1$ in~\eqref{eq::D1-define}: $\norm{h_{T_0^c}}_1  \leq D_1 \lambda \sigma s_0$,
\bens
\twonorm{h_{T_{01}}}
& \le & D_0  \lambda \sigma  \sqrt{s_0}, \quad
\twonorm{h_{T_{01}^c}} \le  \norm{h_{T_0^c}}_1  /\sqrt{s_0} \le D_1
\lambda \sigma \sqrt{s_0}, \\
\twonorm{h}^2
& = & \twonorm{h_{T_{01}}}^2 +\twonorm{h_{T_{01}^c}}^2
\le (D_0^2 + D_1^2) \lambda^2 \sigma^2  s_0, \text{ and }\; \\
\twonorm{\beta_{\init}-\beta} & \le &
\twonorm{h} + \lambda \sigma \sqrt{s_0} \le
[\sqrt{D_0^2 + D_1^2} + 1]
\lambda \sigma \sqrt{s_0}.
\eens
\end{lemma}



\begin{remark}
One can also use $\RE(s_0, 6, X)$ in
  {\bf Case 2} to further tighten the constants in
  \eqref{eq::HT01Case3} and \eqref{eq::Dlocal} for {\bf Case 3},
  for example, by considering 
\bens
\label{eq::DHhead2}
\text{ in Case 2}, \quad
 \norm{X\delta}^2_n  + \norm{X h}^2_{n} \le 
 3 \lambda_n \norm{h_{T_0}}_1,
  \eens 
instead of~\eqref{eq::DHhead}, and in {\bf Case 3}, instead of~\eqref{eq::DHmiddlemain}, 
  \bens
 3 \lambda_n \norm{h_{T_0}}_1 & \le &\norm{X\delta}^2_n  + \norm{X h}^2_{n}
  \le \norm{X \beta -  X \beta_0}_n^2. 
  \eens
Now clearly $\frac{\Lambda_{\max}(s- s_0)}{\Lambda_{\min}(2s_0)} \ge
(\frac{\Lambda_{\max}(s - s_0)}{\Lambda_{\min}(2s_0)})^{1/2}$ for $s \ge 2 s_0$ since
\bens
\Lambda_{\min}(2s_0) \le \Lambda_{\min}(s_0)
\le \Lambda_{\max}(s- s_0), \; \text{ and }  1/\sqrt{\Lambda_{\min}(2s_0)}\le  \sqrt{2} K(s_0, 1).
\eens
Compared with {\bf Case 1} and {\bf Case  2}, where $\RE(s_0, k_0, X)$
holds for $k_0 =3, 4$ respectively, we have an extra term on the 
bound for $\twonorm{h_{T_{01}}}$ in~\eqref{eq::Dlocal}, namely, 
$\ell(s_0) \norm{h_{T_0^c}}_1 /\sqrt{s_0}
\asymp \lambda \sigma \sqrt{s_0}$; 
cf.~\eqref{eq::HT01Case3} and \eqref{eq::definell}.
\end{remark}

\begin{remark}
Under the UUP condition as in Definition~\ref{def:CT-cond}:
\ben
\label{eq::tighttheta}
\forall s_0 \le s, \quad \theta_{s_0, 2s_0} + \delta_{2s_0} \le
\theta_{s, 2s} + \delta_{2s} < 1
\een
since $\theta_{s, 2s}$ and
$\delta_s$ are nondecreasing in $s$. Thus we have
$\theta_{s_0, 2s_0} < 1 -  \delta_{2s_0} \le \Lambda_{\min}(2s_0)$.
When such tight bounds are not available, we can still use the bounds
in~\eqref{eq::orthocauchy} to control $\theta_{s_0, 2s_0}$.
Moreover, suppose $\theta_{s_0, 2s_0} < \Lambda_{\min}(2s_0)$, which holds under 
UUP~\eqref{eq::tighttheta}, one can get rid of the factor 
$\sqrt{\Lambda_{\max}(s_0)}$ in \eqref{eq::definell},
where recall $\ell(s_0)  = ({\theta_{s_0, 
    2s_0} }/{\sqrt{\Lambda_{\min}(2s_0)}}) \wedge \sqrt{\Lambda_{\max}(s_0)}$,
by bounding $\twonorm{h_{T_{01}}}$ for {\bf Case   3}  following
\eqref{eq::simpleT01orig} instead of the supplementary
Lemma~\ref{lemma::T01simple}.
\end{remark}

\section{On Type II errors and $\ell_2$-loss optimality}
\label{sec:type-II-intro}
So far, we have focused on controlling Type I errors,
which would be meaningless if {\it significant variables} are all
missing.
Our goal in this section is to show when
$\beta_{\min,A_0}$~\eqref{eq::betaminA0} is sufficiently large,  we
have $A_0 \subset I$ while achieving the sparse oracle inequalities.
Under the $\RE$ and sparse eigenvalue conditions, 
this result is shown in Theorem~6.3~\citep{Zhou10}, 
which is a corollary of Lemma~\ref{lemma:threshold-general-II}. 
We include it here for self-containment.
First, we state Proposition~\ref{PROP:COUNTING-S0}.
\begin{proposition}
\label{PROP:COUNTING-S0}
Let $A_0 := \{j: |\beta_j| >  \sqrt{2 \log p/n} \sigma\}$.
Let $T_0$ denote 
positions of the $s_0$ largest coefficients of $\beta$ in absolute values,
where $s_0$ is defined in~\eqref{eq::define-s0}.
Let $a_0 = \size{A_0}$ denote the cardinality of $A_0$ 
(see also~\eqref{eq::betaminA0}).
Then $\forall c' > 1/2$, we have
$\size{ \{j \in T_0^c: |\beta_j| \geq \sqrt{\log p/(c' n)} \sigma \}}
\leq  (2c' -1) (s_0  - a_0)$.
\end{proposition}
Again order $\beta_j$'s in decreasing order of magnitude: 
$|\beta_1| \geq |\beta_2| \geq... \geq |\beta_p|$.  Let  $T_0 = \{1, 
\ldots, s_0\}$.
One could choose another target set: for example 
$\{j: |\beta_j| \geq  \sqrt{\log p/(c' n)} \sigma\}$, for some $\log 
p/2 > c' > 1/2$.
Moreover, we consider the consequence of setting 
$t_0 \in [\sigma \sqrt{2/(n)},\sigma \sqrt{2 \log p/(n)}]$.
In particular, when we set $c' = \log p/2$, we have 
\ben
\label{eq::tailcount}
\size{ \{j \in T_0^c: |\beta_j| \geq \sigma \sqrt{2/(n)} \}}
\leq  (\log p -1) (s_0  - a_0).
\een
We prove Proposition~\ref{PROP:COUNTING-S0} in Section~\ref{sec::proofprop1}. 
We first show in Lemma~\ref{lemma:threshold-general} that 
under no restriction on $\beta_{\min}$,  we achieve an oracle bound on the 
$\ell_2$ loss, which depends only on 
the $\ell_2$ loss of the initial estimator on the set $T_0$. 
Bounds in Lemma~\ref{lemma:threshold-RE} are special cases
of~\eqref{eq::off-beta-norm-bound-2}.
In Lemma~\ref{lemma:threshold-general-II},
we impose a lower bound on $\beta_{\min, A_0}$ as
in~\eqref{eq::betaA-min-cond} in order to recover the subset of
variables in $A_0$, while achieving the  nearly ideal $\ell_2$ loss
with a sparse model $I$ in view of Lemma~\ref{prop:MSE-missing-orig}.
We prove Lemma~\ref{lemma:threshold-general} and~\ref{lemma:threshold-general-II} in 
Sections~\ref{sec::proofofBias} and~\ref{sec::proofofA0} respectively.
\begin{lemma}{\textnormal({\bf A deterministic result on the bias.})}
\label{lemma:threshold-general}
Let $\beta_{\init}$ be an initial estimator. 
Let $\lambda := \sqrt{2 \log p/n}$ and $h \in \R^p$ be as defined in Lemma~\ref{lemma:threshold-RE}.
Suppose we choose a thresholding parameter $t_0$ and set
$$I = \{j: \size{\beta_{j, \init}} \geq t_0\}.$$
Then for $\drop := [p] \setminus I$, we have for 
$\drop_{11} := \drop \cap A_0$ and $a_0 = \size{A_0}$,
\begin{eqnarray}
\label{eq::off-beta-norm-bound-2}
\twonorm{\beta_{\drop}}^2 & \leq  & (s_0 - a_0) \lambda^2 \sigma^2 + 
(t_0 \sqrt{a_0} + \twonorm{h_{\drop_{11}}})^2.
\end{eqnarray}
Suppose $t_0 < \beta_{\min,A_0}$ as defined in~\eqref{eq::betaminA0}.
Then~\eqref{eq::off-beta-norm-bound-2} can be replaced by
\begin{eqnarray}
\label{eq::off-beta-norm-bound-alt}
  \twonorm{\beta_{\drop}}^2 
  & \leq  & (s_0 - a_0) \lambda^2 \sigma^2 + 
\twonorm{h_{\drop_{11}}}^2 
\left({\beta_{\min,A_0}}/{(\beta_{\min,A_0} - t_0)}\right)^2.
\end{eqnarray}
\end{lemma}

\begin{lemma}
{\textnormal({\bf A deterministic result on the bias and model selection.})}
\label{lemma:threshold-general-II}
Let $\beta_{\init}$ be an initial estimator.
Let $\lambda := \sqrt{2 \log p/n}$ and
$h \in \R^p$
be as defined in Lemma~\ref{lemma:threshold-RE}.
Suppose
for $\beta_{\min,A_0}$ as defined in~\eqref{eq::betaminA0}, it holds that
\ben
\label{eq::betaA-min-cond}
\beta_{\min, A_0}  \geq \norm{h_{A_0}}_{\infty} + 
\min\left\{(s_0)^{-1/2} \twonorm{h_{T_0^c}}, \;
(s_0)^{-1} \norm{h_{T_0^c}}_1 \right\}.
\een
Now we choose a thresholding parameter $t_0$ such that
\ben
\label{eq::ideal-t0}
\lefteqn{\text{for some } \; \breve{s}_0 \in [s_0, s), \quad
\beta_{\min, A_0}  -  \norm{h_{A_0}}_{\infty} \geq t_0}\\
& \geq &
\nonumber
\min\left\{(\breve{s}_0)^{-1/2} \twonorm{\beta_{\init, T_0^c}},
(\breve{s}_0)^{-1} \norm{\beta_{\init, T_0^c}}_1\right\}
\een
holds and set  $I = \{j: \size{\beta_{j, \init}} \geq t_0\}$.
Then we have
\ben
\label{eq::lemma-last1}
A_0 \subset I && \text{ and } \quad |I \cap T_0^c| \leq \breve{s}_0, \quad
\text{ and hence} \\
\label{eq::lemma-last2}
|I| \leq s_0 + \breve{s}_0  && \text{ and  } \quad
\twonorm{\beta_{\drop}}^2 \leq (s_0 - a_0) \lambda^2 \sigma^2.
\een
\end{lemma}

Choosing the set $A_0$ as in Proposition~\ref{PROP:COUNTING-S0} is
rather arbitrary; one could for example,  consider the set of
variables that are strictly above $\lambda \sigma/2$.
In Theorem~\ref{thm::threshA0},
we show that one can indeed recover a subset 
$A_0$ of variables accurately, for $A_0$ as defined in
Proposition~\ref{PROP:COUNTING-S0},
when  $\beta_{\min,A_0} := \min_{j \in A_0} |\beta_{j}|$ is large enough 
(relative to the $\ell_2$ loss of an initial estimator $\beta_{\init}$ under the $\RE$
condition on the set $A_0$); in addition, a small number of extra 
variables from $T_1 \subset T_0^c := [p] \setminus T_0$ are also
possibly included in the model $I$.
We mention in passing that changing the coefficients of $\beta_{A_0}$
will not change the values of $s_0$ or $a_0$, so long as their
absolute values stay strictly above $\lambda \sigma$; See
Section~\ref{sec::newmodel}.
Compared with the {\it almost exact} sparse recovery result in
Theorem~1.1~\citep{Zhou09th}, we have relaxed the restriction on
$\beta_{\min}$: rather than requiring all non-zero entries to be large
in absolute values, we only require those in a subset $A_0$ to be
recovered to be large.

\begin{theorem}
  \label{thm::threshA0}
  Under the settings of Theorem~\ref{thm::RE-oracle-main},
  let $\beta_{\init}$ be the Lasso estimator as in~\eqref{eq::origin}.
  Let $T_0$ denote the largest $s_0$ coordinates of $\beta$ in absolute 
values.  Let $T_1$ be the largest $s_0$ positions of $\beta_{\init}$
outside of $T_0$. Let $T_{01} = T_0 \cup T_1$. 
Suppose~\eqref{eq::eigen-max},~\eqref{eq::eigen-admissible-s}, and
$\RE(s_0, 4, X)$ condition hold.
Let $D_0, D_1$ be as in~\eqref{eq::D0-define-orig}
and~\eqref{eq::D1-define}. 
Suppose for some constants $\breve{D}_1 > D_1$, $\breve{D}_0 > D_0$,
\ben
\label{eq::betaA-min-RE}
\beta_{\min, A_0} := \min_{j \in  A_0} |\beta_{j}| \geq  D_0 \lambda \sigma \sqrt{s_0} 
+ (\breve{D}_1 \wedge \breve{D}_0)
\lambda \sigma
\een
where $\lambda := \sqrt{2 \log p/n}$.
Choose a thresholding parameter $t_0$ and set 
\bens
I = \{j: \size{\beta_{j, \init}} \geq t_0\}, \; \text{ where }
t_0 = (\breve{D}_1 \wedge \breve{D}_0)\lambda \sigma.
\eens
Let $\hat\beta_{I} = (X_I^T X_{I})^{-1} X_{I}^T Y$ and $\hat\beta^{\OLS}(I)$ be the $0$-extended version of 
$\hat\beta_{I}$ such that $\hat\beta^{\OLS}_{I^c} = 0$ and
$\hat\beta^{\OLS}_I = \hat\beta_{I}$.
Then on event $\T_a$, we have $\abs{I} \le  2s_0$, 
\ben
\nonumber
&& A_0 \subset I, \quad \; I \cap T_{01}^c = \emptyset, \;
\twonorm{\beta_{\drop}} \le \lambda \sigma\sqrt{s_0-a_0}, \text{ and} \\
\label{eq::threshold-general-II}
&& \twonorm{\hat{\beta}^{\OLS}(I) - \beta}^2 \leq 
\big(1+ \frac{(\Lambda_{\max}(2s) -\Lambda_{\min}(2s))^2 + 8(1+a)}{2 
  \Lambda_{\min}^2(|I|)}  \big) s_0 \lambda^2 \sigma^2.
\een
\end{theorem}

\begin{proof}
By definition of $T_0$, we have $T_0 \subset S$ and recall $\dropS = \drop
\cap S$. Moreover, $I \cap \dropS = \emptyset$ by definition.
By Lemma~\ref{lemma::HT01case3},
  we have for on $\T_a$,
  \ben
  \label{eq::L1norm}
\norm{\beta_{\init, T_{0}^c}}_1 = \norm{h_{T_{0}^c}}_1 \le D_1 \lambda 
\sigma s_0 \text{ and }\;
\twonorm{h_{T_{01}}}  \le D_0 \lambda \sigma \sqrt{ s_0}.
\een
First, we find the lower bound on $\beta_{\min, A_0}$ so that $A_0 \subset I$.
Suppose there exists a threshold $t_0>0$ such that one of
the following conditions holds,
\ben
\nonumber
\beta_{\min, A_0}
- \norm{h_{A_0}}_{\infty} & \geq &  t_0
\text{ or the stronger  }
\beta_{\min, A_0} \ge \twonorm{h_{A_0}} +
t_0; \\
\label{eq::noFN}
\text{Then }\;
\forall j \in A_0,\; 
\abs{\beta_{\init, j}} & \geq &
\abs{\beta_{\min, A_0} - \norm{h_{A_0}}_{\infty} }> t_0,
\een
ensuring no FNs from $A_0$.
Next, we derive a lower bound on $t_0$.
Suppose 
\ben 
\label{eq::t0lowerbound}
t_0 > (D_0 \wedge D_1) \lambda \sigma \ge
\norm{\beta_{\init, 
    T_{01}^c} }_{\infty}
=\norm{h_{T_{01}^c} }_{\infty},
\een
we can then eliminate all variables in $T_{01}^c$ from model $I$,
where $\beta_{\init, T_{0}^c} = h_{T_{0}^c}$ and $T_1 \subset T_0^c$.
To see the second inequality in \eqref{eq::t0lowerbound}, recall the
$k$th largest value of $\abs{h_{T_0^c}}$ obeys
$\size{h_{T_{0}^c}}_{(k)} \leq  \norm{h_{T_0^c}}_1 / k$ for $k \ge 1$.
Hence by definition of $T_1$,
\bens 
\norm{\beta_{ \init, T_{01}^c} }_{\infty}
& = &
\size{h_{T_{0}^c}}_{(s_0+1)}
\le \norm{h_{T_0^c}}_1 / (s_0+1)  < D_1 
\lambda \sigma \; \text{by~\eqref{eq::L1norm}}.
\eens
Moreover, by definition of $T_1$, the largest entry in
$\abs{h_{T_{01}^c}}$
(entrywise absolute value of vector $h_{T_{01}^c}$)
is bounded by the average over the top $s_0$ largest entries of
$\abs{h_{T_0^c}}$:
$\size{h_{T_{0}^c}}_{(s_0+1)} \le \inv{s_0} \sum_{k=1}^{s_0} \size{h_{T_{0}^c}}_{(k)} = \norm{h_{T_1}}_1/s_0$,  and hence
\bens
\norm{\beta_{ \init, T_{01}^c} }_{\infty}
& = &
\norm{h_{T_{01}^c} }_{\infty}\le \norm{h_{T_{1}}}_1/s_0 \le \sqrt{s_0}
\twonorm{h_{T_{1}}}/s_0  \\
& \le &
\twonorm{h_{T_{01}}}/\sqrt{s_0} \le D_0 \lambda \sigma \; \text{ by~\eqref{eq::L1norm}}.
\eens

In summary, \eqref{eq::betaA-min-RE} and  \eqref{eq::t0lowerbound}
imply that it is feasible to set 
\ben
\nonumber
t_0 & = & (\breve{D}_1 \wedge \breve{D}_0) \lambda \sigma > 
(D_1 \wedge D_0) \lambda \sigma \; \text{ since} \\
\label{eq::betaminA0proof}
\beta_{\min, A_0}
& \ge &  \twonorm{h_{T_0}} +  t_0 > \norm{h_{T_0}}_{\infty} +  (D_0
\wedge D_1) \lambda \sigma.
\een
Then none of the variables in $T_{01}^c$ can be chosen in the thresholding 
step and hence $I \cap  T_{01}^c = \emptyset$. Hence, $\size{I} \le 
2s_0$. On event $\T_{a}$, \eqref{eq::plossintro} holds by
Lemma~\ref{prop:MSE-missing-orig}. Moreover, $\theta^2_{|I|,
  |\dropS|}$ is bounded as in Lemma~\ref{lemma:parallel}
given $|I| + |\dropS| = |I \cup \dropS| \leq \abs{S} + \abs{I \cap S^c} \le s + |I \cap T_0^c| \leq s + s_0 \leq 2s$.
Thus~\eqref{eq::threshold-general-II} follows immediately.
\qed 
\end{proof}

\subsection{Discussions}
\begin{figure}
\begin{center}
\includegraphics[width=0.8\textwidth]{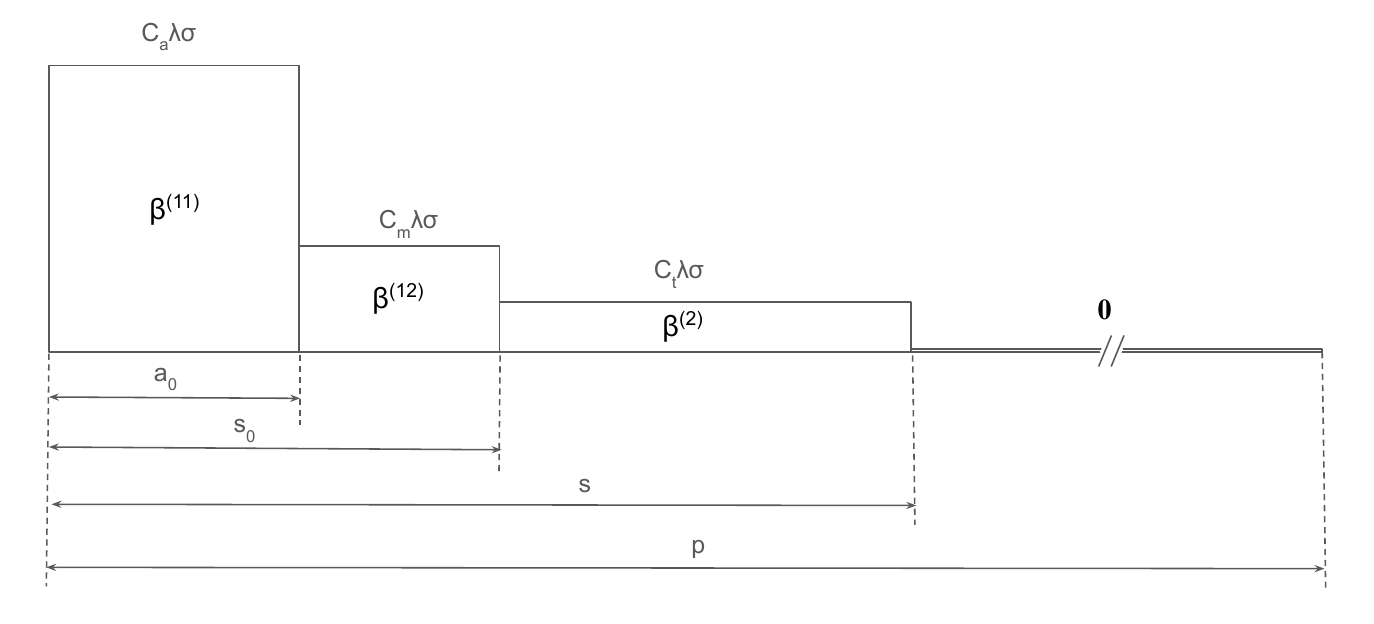} 
\caption{
In this model, the component $\beta^{(11)}$ has $a_0$ non-zero coordinates
with the same magnitude $C_a  \lambda \sigma =:\beta_{\min, A_0}$,
where $C_a \in \{1.706, 8.528\}$
and $\beta_{\min, A_0} \in \{0.2, 1\}$; the component $\beta^{(12)}$ has $s_0 - a_0$
non-zero coordinates with the same magnitude $C_m \lambda \sigma$, where $C_m =
1/{\sqrt{2}}$ for $s> s_0$ and $C_m=1$ in case $s_0 = s$;
the component $\beta^{(2)}$ has $s - s_0$ non-zero coordinates with the same
magnitude $C_t \lambda \sigma =: c_t \sigma/\sqrt{n}$. See
\eqref{eq::tailcount}. The rest are all 0s.
In the exact sparse case, namely, when $s=s_0$, all non-zero signals
are concentrated on the component $\beta^{(1)}$ without spreading
across components of $\beta^{(2)}$.}
\label{fig:tiger-beta}
\end{center}
\end{figure}

Theorem~\ref{thm::threshA0} is an immediate corollary of
Theorem~\ref{thm::RE-oracle-main}, Lemma~\ref{prop:MSE-missing-orig},
and Lemma~\ref{lemma:threshold-general-II}, 
except that we now set $\breve{s}_0 = s_0$ everywhere in Lemma~\ref{lemma:threshold-general-II}, and assume
having an upper estimate $\breve{D}_1$ (resp. $\breve{D}_0$)
of $D_1$  (resp. $D_0$), so as not to depend on an ``oracle'' telling
us an exact value.
 Given all other parameters held invariant, a larger
 $\beta_{\min,A_0}$ (e.g., a larger $C_a$ for $\beta$ in
 Figure~\ref{fig:tiger-beta}) will result in a tighter bound on
 the bias $\twonorm{\beta_{\drop}}^2$ in view
 of~\eqref{eq::off-beta-norm-bound-alt}.
 When $\beta_{\min, A_0}^2$ dominates the upper bound on the RHS of~\eqref{eq::off-beta-norm-bound-alt}, then $A_0 \cap I^c = \emptyset$,
since the loss of a single variable from $A_0$ will already saturate
this upper bound in the order of $O(\lambda^2 \sigma^2 s_0)$.
Hence when $\beta_{\min, A_0}$ is strong enough in the sense of
\eqref{eq::ideal-t0}, there will be no false negatives from the set
$A_0$, leading to the stronger and tighter bounds in~\eqref{eq::lemma-last1}
and~\eqref{eq::lemma-last2} that also control false positives and the
bias. This is validated in our experiments in
Section~\ref{sec::numeric}.

In the statement of Lemma~\ref{lemma:threshold-general-II}, we assume
the knowledge of the error bounds on various norms of $\beta_{\init} -
\beta$ and $h =\beta_{\init} - \beta_{T_0}$  implicitly
(hence the name of ``oracle''). 
In obtaining~\eqref{eq::betaA-min-RE}, we may substitute the bounds 
as derived in Lemma~\ref{lemma::HT01case3}
in~\eqref{eq::betaA-min-cond}, as 
trivially, for $A_0 \subseteq T_0$ and $h_{A_0}= \beta_{\init, A_0} -
\beta_{A_0}$, 
\ben 
\label{eq::part-norms-app}
\norm{h_{A_0}}_{\infty} & \leq &
\twonorm{h_{A_0}} \leq \twonorm{h_{T_{0}}} \leq 
D_0 \lambda \sigma \sqrt{s_0} \;  \text{ and } \\
\nonumber 
\norm{h_{T_0^c}}_1/s_0 & \leq & D_1 \lambda \sigma \leq \breve{D}_1 
\lambda \sigma \; \;  \text{ in view of  \eqref{eq::L1norm}}. 
\een
Lemma~\ref{lemma:threshold-general-II} explains our results in Theorem 1.1~\citep{Zhou09c} as well as the numerical results. We also introduce
$\breve{s}_0$ so that  the dependency of $t_0$ on the knowledge of
$s_0$ is relaxed; in particular, it can be used to express a desirable
level of sparsity for the model $I$ that one wishes to select.
In general, we may assume that $\beta_{\min, A_0}  \geq  2 D_0 \lambda 
\sigma \sqrt{s_0}$ so as to have a large range of effective 
thresholding parameters; cf. Section~\ref{sec::proofsketchRE}.
Thus one can increase $t_0$ as $\beta_{\min,A_0}$ increases in order to 
reduce the number of false positives while not increasing the number 
of false negatives from the active set $A_0$.

In this case, it is also possible to remove
variables from $T_0^c$ entirely by  increasing the threshold $t_0$
while strengthening the lower bound on $\beta_{\min,A_0}$ by a
constant factor. For example, we may set $t_0 > \norm{\beta_{\init,
    T_0^c}}_{\infty}$, rather than $\norm{\beta_{\init,    T_{01}^c}}_{\infty}$ as
in~\eqref{eq::t0lowerbound}.
Clearly, $\norm{\beta_{\init,    T_0^c}}_{\infty} = \norm{h_{T_1}}_{\infty} \le
\twonorm{h_{T_{1}}} =O_P(\lambda \sigma \sqrt{s_0})$ and hence the
condition $\beta_{\min, A_0} \ge  \sqrt{2} \twonorm{h_{T_{01}}} \geq   \twonorm{h_{T_0}}
+ \twonorm{h_{T_{1}}}$ ensures
\bens
\beta_{\min, A_0} - \norm{h_{A_0}}_{\infty} & \geq &
\beta_{\min, A_0} - \twonorm{h_{T_0}}  \ge \twonorm{h_{T_{1}}} >
\norm{\beta_{\init,    T_0^c}}_{\infty},
\eens
which is sufficient for~\eqref{eq::noFN} to hold so long as $t_0 \asymp
\norm{\beta_{\init,    T_0^c}}_{\infty} \le  \twonorm{h_{T_{1}}}$.
Then following the same analysis as in Theorem~\ref{thm::threshA0}, we
can recover $A_0$ while eliminating variables from $T_0^c$;
cf. \eqref{eq::betaminA0proof}.
In general, if the strong signals are close to each other in their strength,
then a small $\beta_{\min, A_0}$ implies that we are in a situation with low 
signal to noise ratio (low SNR); one needs to carefully tradeoff false positives 
with false negatives as shown in our numerical results in Section~\ref{sec::numeric} and the supplementary
Section~\ref{sec:experiments}. Such results and their formal statements
are omitted from the present paper.

Bounds on $\norm{h_{A_0}}_{\infty}$ are in general 
harder to obtain than $ \twonorm{h_{A_0}}$.
In general, we can still hope to  bound $\norm{h_{A_0}}_{\infty}$ by
$\twonorm{h_{A_0}}$ as done in Theorem~\ref{thm::threshA0}.
Having a tight bound on  $\twonorm{h_{T_0}}$ (or 
$\norm{h_{T_0}}_{\infty}$) and $\shtwonorm{h_{T_0^c}}$ naturally
helps relaxing the requirement on $\beta_{\min, A_0}$ for 
Lemma~\ref{lemma:threshold-general-II}, while as shown in
Lemma~\ref{lemma:threshold-general}, such tight upper bounds will help
us control the model size $\abs{I}$ and the bias
$\norm{\beta_{\drop}}$ and therefore achieve a tight bound on the
$\ell_2$ loss in the statement of Lemma~\ref{prop:MSE-missing}.
We refer to~\cite{Wai09b} for discussion and further pointers into
this literature on information theoretic limits on sparse recovery.

\subsection{Model specification}
\label{sec::newmodel}
We now pause to briefly describe our experiment setup so as to further 
discuss variable selection in the context of
Theorem~\ref{thm::threshA0}.
Moreover, we generate a class of models with different sparsity $s$
and $\ell_1$ norm on  $\beta_{T_0^c}$ to 
shed light on this connection with the work of~\cite{ZH08}.
Let $\abs{\beta_j}$'s be ordered as in~\eqref{eq::beta-2-small-intro}.
Let $S =\supp(\beta)$ and $A_0, T_0$ be as defined in 
Proposition~\ref{PROP:COUNTING-S0}.
Let $A_0 = [a_0]$ and $T_0 = [s_0]$.
Then $T_0\setminus A_0 =\{a_0+1,  \ldots, s_0\}$.
  Let $\abs{S} =s$.
See Figure~\ref{fig:tiger-beta} for an illustration of the model 
specification. We divide $\beta$ into 4 components and write $\beta = \beta^{(11)} +
  \beta^{(12)} + \beta^{(2)} + \beta^{(0)} \in \R^p$,
where we
set for some $C_a > 1$, $0< C_m \le 1$, $c_t \ge 0$,
\ben
\label{eq::tigermodel}
\beta_{j}^{(11)}
& = & \pm C_a \lambda \sigma \cdot 1_{1 \leq j \leq a_0}, \quad
\beta_{j}^{(12)}  =   \pm C_m \lambda \sigma \cdot 1_{a_0 < j \leq
  s_0}, \\
\nonumber
\beta_j^{(2)} & = &   \pm c_t  (\sigma/\sqrt{n}) \cdot 1_{s_0 < j \le s}\;
\; \; \text{ and} \; \; \beta_j^{(0)}  =  0\cdot 1_{s<j \le p},
\een
where $\lambda = \sqrt{ 2 \log p/n}$ and the unspecified coordinates
in each component are again set to $0$.
The first 3 components contain non-zero coordinates.
Hence $\twonorm{\beta^{(12)}} =C_m  \lambda  \sigma \sqrt{s_0 - a_0}$
in \eqref{eq::tigermodel}.
Since $C_a >1$,  we have by~\eqref{eq::define-s0},
\ben
\label{eq::A0norm}
\sum_{j \leq a_0} \min(\beta_j^2, \lambda^2 \sigma^2) & = &    
a_0 \lambda^2 \sigma^2, \text{ since } \abs{\beta_j} > \lambda \sigma
\text{ for}\;  j \in A_0 =[a_0], \\
\label{eq::SR-range-intro}
(s_0-a_0) \lambda^2  \sigma^2(1-C_m^2)
& \ge &  \twonorm{\beta^{(2)}}^2 =\sum_{j >  s_0} \beta_j^2 =
\sum_{j > s_0} \min(\beta_j^2, \lambda^2 \sigma^2) \\
\label{eq::etalowermain}
& \ge &
(s_0 - a_0)(1-C^2_m) \lambda^2  \sigma^2 -\lambda^2 \sigma^2. 
\een
Clearly, for $\beta^{(2)}$, its $\ell_1$ norm $\norm{\beta_{T_0^c}}_1 
=\norm{\beta^{(2)}}_1$ is proportional to the signal strength 
$\twonorm{\beta^{(2)}}^2$ and is inversely proportional to its 
$\ell_{\infty}$ norm, since by \eqref{eq::SR-range-intro},
\ben
\nonumber
(\text{support}) \quad 
\abs{T_0^c \cap S} & = & s-s_0 = \twonorm{\beta^{(2)}}^2/(c_t^2 \sigma^2/n) \\
\label{eq::beta2supp}
& \le & {(s_0-a_0) (2 \log p)(1-C_m^2)}/{c_t^2}, \; \text{where} \\
\label{eq::eta1c}
(\ell_2) \quad \norm{\beta^{(2)}}_2
& = &
\twonorm{\beta_{{T_0^c}}}  \approx \lambda  \sigma \sqrt{(1-C_m^2)
  (s_0 - a_0)}, \\
\text{and} \; (\ell_1) \quad 
\norm{\beta_{T_0^c}}_1 
&= &
\nonumber
\sqrt{\abs{T_0^c \cap S} } \twonorm{\beta_{T_0^c}} =
\twonorm{\beta^{(2)}}^2/[c_t \sigma/\sqrt{n}]  \\
\label{eq::eta1d}
&\le &
\frac{\sqrt{2\log p}(1-C_m^2) }{c_t} \lambda  \sigma  (s_0  - a_0).
\een
Here \eqref{eq::eta1d} is essentially tight with a matching lower bound at 
the same order; cf. \eqref{eq::etalowermain}.
In other words, by choosing different values of 
$C_t$, we generate the model class for $\beta$ with different 
sparsity as in \eqref{eq::tigermodel} and \eqref{eq::CT}: 
\ben 
\label{eq::CT}
(s - s_0) C_t^2  & = & (s_0 - a_0)  (1 - C_m^2), \text{ where} \;  C_t 
= c_t /\sqrt{2 \log  p} \ge 0. 
\een
By setting the height of $\beta^{(2)}$, or $\norm{\beta^{(2)}}_{\infty}$ to be $\asymp \sigma/\sqrt{n}$, the size of its support $\abs{\supp(\beta^{(2)})} \asymp \log p (s_0 - a_0)$
as desired; cf. Proposition~\ref{PROP:COUNTING-S0}
and~\eqref{eq::beta2supp}. Finally \eqref{eq::eta1c} holds
by \eqref{eq::A0norm}, \eqref{eq::SR-range-intro}, and \eqref{eq::etalowermain}.
We compare the $\ell_1$ condition in~\eqref{eq::eta1d} with those
in~\cite{ZH08}; cf. \eqref{eq::eta1def}.

We now show that the upper bound in~\eqref{eq::eta1d} is essentially tight.
First, the  bound in \eqref{eq::eta1c} is tight, since
\ben
\label{eq::beta2norms}
\twonorm{\beta^{(2)}}^2 & \le &
 \norm{\beta^{(2)}}_1  \norm{\beta^{(2)}}_{\infty} 
= [c_t \sigma/\sqrt{n}] \norm{\beta^{(2)}}_1, 
\text{ where} \\
\nonumber
\sum_{j  \in T^c_0} \min(\beta_j^2, \lambda^2 \sigma^2) 
& = &
\sum_{j > s_0} \min(\beta_j^2,  \lambda^2 \sigma^2) = \twonorm{\beta^{(2)}}^2 \\
\label{eq::etalower}
& \ge &
(s_0 - a_0)(1-C^2_m) \lambda^2  \sigma^2 -\lambda^2 \sigma^2.
\een
Moreover, \eqref{eq::eta1d} is also essentially tight by definition
since by~\eqref{eq::beta2norms} and \eqref{eq::etalower},
\ben
\nonumber
\norm{\beta^{(2)}}_1
& \ge &
\twonorm{\beta^{(2)}}^2 / [c_t \sigma/\sqrt{n}] \ge ((1-C^2_m)(s_0 -
a_0) -1) \lambda^2 \sigma /[c_t /\sqrt{n}] \\
& \ge &
\label{eq::lowereta1a}
((1-C^2_m)(s_0 - a_0) -1) \lambda \sigma \frac{\sqrt{2 \log p}}{c_t} 
\; \text{ and moreover} \\
\nonumber
\label{eq::lowereta1b}
\sum_{j  \in A_0^c} \abs{\beta_j}
& = &
\norm{\beta^{(2)}}_1  + C_m \lambda \sigma (s_0 -a_0)
 = \Omega(\lambda \sigma(s_0 - a_0) {\sqrt{(\log p)}}/{c_t}). 
\een
Finally, \eqref{eq::etalower} follows
from~\eqref{eq::SR-range-intro} and the fact that
$\twonorm{\beta^{(12)}}^2  = \sum_{a_0 < j \le s_0} \min(\beta_j^2, \lambda^2 \sigma^2) =
C^2_m \lambda^2 \sigma^2 (s_0 -a_0)$.

\subsection{Conditions in~\cite{ZH08}}
\label{sec::tightness}
Without loss of generality, we order $\abs{\beta_1} \ge \abs{\beta_2} \ge 
\ldots \ge \abs{\beta_p}$  as in~\eqref{eq::beta-2-small-intro}.
For the Lasso solution $\beta_{\init}$, it is assumed that
$\lambda_n \ge 2 \lambda \sigma \sqrt{(1+a) c^*}$
in~\eqref{eq::origin} in Theorem~3 in~\cite{ZH08}.
Hence we assume $\lambda_n \asymp 2 \lambda \sigma \sqrt{(1+a) c^*}$
throughout our discussion.
Roughly speaking, we interpret
$H_q =\{1, \ldots, q\}$ as the set of large coordinates in $\beta$
that one aims to recover in the settings of~\cite{ZH08}.
For convenience, we denote by $H_0 := [p] \setminus H_q = \{q+1,
\ldots, p\}$ and $T_1^* \subset H_0$ the index set $\{j \in H_0: \abs{\beta_j}
\ge \lambda \sigma\}$.
The sparsity assumption by~\cite{ZH08} is set
with the $\ell_1$-sparsity on $\beta_{H_0}$,
\ben
\label{eq::eta1def}
\sum_{j > q}^p \abs{\beta_j} \le
\eta_1 =\tilde{O}(r_1^2 q \lambda \sigma), \; \text{ where } \eta_1
\le 2 \lambda \sigma \frac{r_1^2 q}{\sqrt{c^{*}}}
=O(\frac{q^{*}}{\sqrt{c^{*}}} \lambda \sigma);
\een
Here the $\tilde{O}$ notation may hide some constants including $c^*
>1$ and $q^* > 4 r_1^2 q$, where $q^*, c_*, c^*, q, r_1$ (and $r_2$)
are all allowed to depend on $n$; cf. \eqref{eq::eta2def} to~\eqref{eq::qstar}.
This target set $H_q$ is slightly more restrictive than the set $A_0 =\{1, \ldots, 
a_0\}$,
since in $H_0$, one may still see large signals;
cf. the supplementary~\eqref{eq::largesignal} and the proof of Theorem 3~\citep{ZH08}.
To resolve this discrepancy and to properly interpret results in~\cite{ZH08},
we first extend $H_q$ by $T_1^* \subset H_0$ and denote this extended set by $L(q)$, where
\bens 
L(q) & := & H_q \cup T_1^*  = \{j: \abs{\beta_j} \ge \lambda
\sigma\},\;
\abs{T_1^*} \le  \eta_1/(\lambda \sigma)
= \tilde{O}(2 r_1^2 q),\\
\text{and }\;  \abs{L(q)} & \le &  \frac{2 r_1^2
  q}{\sqrt{c^*}} + q <\frac{q^*}{2\sqrt{c^*}}\; \; \; \text{ in view 
  of~\eqref{eq::qstar}}.
\eens
By construction, $A_0 \subseteq L(q) \subset T_0$ by  
definition of $A_0$~\eqref{eq::betaminA0} and $T_0$.
Hence 
\ben
\label{eq::qstars0}
a_0 \le \abs{L(q)} \le s_0, \quad q^* =\Omega(s_0 \sqrt{c^*}), \text{
  so long as } a_0 \asymp s_0;
\een
Moreover,~\cite{ZH08} require the following $\eta_2$ condition in
\eqref{eq::eta2def} and impose the sparse Riesz conditions with rank
$q^*$ as in~\eqref{eq::qstar} on $X$:
\ben
\label{eq::eta2def}
&& \eta_2/\sqrt{n} :=  \max_{A \subset H_0}
\shtwonorm{\sum_{j\in A} \beta_j X_j}/\sqrt{n} =
\tilde{O}(2 r_2 \sqrt{q} \lambda \sigma), \\
\label{eq::SRC}
&& c_* \le {\twonorm{X_A v}^2}/{(n \twonorm{v}^2)} \le c^* \quad \forall A 
\; \text{with} \quad \size{A} = q^*\; \text{ and } \; v \in \R^{q^*}.
\een
Then, by Eq. (2.15) -- (2.18) and (3.1)~\citep{ZH08},
\ben
\label{eq::qstar}
M_1^* q + 1 := (2 + 4 r_1^2 + 4 \sqrt{\kappa} r_2 + 4 \kappa) q + 1
\le q^*, \; \; \text{  where} \; \; \kappa :=c^*/c_*,
\een
and $r_1$ and $r_2$ are the same as in~\eqref{eq::eta1def}
and~\eqref{eq::eta2def} respectively.
Both conditions are needed to show an upper bound on
$\twonorm{\beta_{H_0}}  = \tilde{O}(\sqrt{q^*} \sigma \lambda)$;
cf.~\eqref{eq::tailL2lossmain}. 
Let $A_1 := \supp(\beta_{\init}) \cup H_q$.
Moreover, since $A_1^c \subset H_0$,  we have by~\eqref{eq::eta1def},
\ben
\nonumber
\norm{\beta_{\init, A_1^c} - \beta_{A_1^c}}_1 & \le &
\norm{\beta_{H_0}}_1 \le \eta_1 \; \; \text{  and }  \\
\label{eq::tailL2lossmain}
\twonorm{\beta_{\init, A_1^c} - \beta_{A_1^c}}^2 & \le &
\twonorm{\beta_{H_0}}^2  = 
\tilde{O}(\frac{r_2}{\sqrt{c_*}}q^* \lambda^2 \sigma^2).
\een
~\cite{ZH08} use~\eqref{eq::qstar}, \eqref{eq::tailL2lossmain}, and the
$\eta_1$ and $\eta_2$ conditions to bound
\ben
\label{eq::L1lossmain}
\twonorm{\beta_{\init} - \beta}  = \tilde{O}_P(\sqrt{q^*} \lambda 
\sigma), \; \text{ and }\;
\norm{\beta_{\init} - \beta}_1 
=  \tilde{O}_P(q^* \lambda \sigma),   
\een
which correctly depend on $q^*$.
Recall in our example in Section \ref{sec::newmodel},
$\twonorm{\beta_{A_0^c}} \le  (s_0  - a_0)^{1/2} \lambda  \sigma$,
where $H_q \subset A_0$ and $T_0^c \subseteq  ([p] \setminus L(q) )
\subseteq A_0^c \subset H_0$, and
\ben
\label{eq::L1A0c}
\norm{\beta_{A_0^c}}_1  & \approx &
\abs{s_0 - a_0} \lambda \sigma (c_t +\sqrt{\log p})/(\sqrt{2} c_t)  \text{ since } \\
\nonumber
\norm{\beta_{T_0^c}}_1 
& \approx &
\frac{\sqrt{\log p}}{\sqrt{2}c_t} (s_0  - a_0) \lambda  \sigma,\;
\;\text{ by  \eqref{eq::eta1d} when }\; C_m = 1/\sqrt{2},  \text{ and } \\
\nonumber
\norm{\beta_{T_0 \setminus A_0}}_1 
& = & \norm{\beta^{(12)}}_1 =
\abs{s_0 - a_0} C_m \lambda \sigma =
\abs{s_0 - a_0} \lambda \sigma/\sqrt{2}.
\een
Notice that the tight upper bound in \eqref{eq::L1A0c} (with matching
lower bound) has an extra $\sqrt{\log p}$ factor compared to the
$\ell_1$-sparsity condition in~\eqref{eq::eta1def}, where it is
understood that $q^{*} \asymp s_0$ in order for the error
$\twonorm{\beta_{\init} - \beta}$ and
$\twonorm{X(\beta_{\init} - \beta)}$ in Theorem 3~\citep{ZH08} to match those in
Theorem~\ref{thm::RE-oracle}; cf. the supplementary~\eqref{eq::signalloss}.
In summary, although we impose sparsity in the $\ell_0$
sense, the actual lower bound on $\norm{\beta_{A_0^c}}_1$,
\bens
\tilde{\eta}_1 := \norm{\beta_{A_0^c}}_1 >
\norm{\beta_{T_0^c}}_1 =\Omega(\sqrt{\log p} \lambda \sigma(s_0 -
a_0))
\eens
covers convergence results not available in
~\cite{ZH08}.
More explicitly, in~\eqref{eq::eta1def},
it is necessary for the $\ell_1$ norm on $\beta_{H_0}$ to satisfy
\bens
\norm{\beta_{H_0}}_1 < \eta_1  =\tilde{O}(r_1^2 q \lambda \sigma)
=O(1)(q^* \lambda \sigma)=O(1)(s_0 \lambda \sigma),
\eens
while in order for $\eta_1 \ge \norm{\beta_{H_0}}_1 \ge
\norm{\beta_{A_0^c}}_1$ to hold, we need to set
\ben
\label{eq::qstarlower}
\; q^{*} \ge \abs{s_0 - a_0}  \frac{c_t +\sqrt{\log p}}{\sqrt{2} c_t}
\text{ so that } 
\eta_1 \asymp  q^{*} \lambda  \sigma >
\norm{\beta_{A_0^c}}_1;
\een
cf. \eqref{eq::lowereta1a} and~\eqref{eq::lowereta1b}.
On the other hand, suppose $q^{*} \asymp  \sqrt{\log p} (s_0 -a_0)$ is
allowed, then this extra $\sqrt{\log p}$ factor in
\eqref{eq::qstarlower} will inevitably appear in the upper bounds
derived in \eqref{eq::L1lossmain} and~\eqref{eq::tailL2lossmain},
resulting in worse $\ell_2$ loss
while simultaneously, the bounded
sparse eigenvalue conditions need to hold for  design matrix $X$ with rank $q^*
=\Omega( (s_0 -a_0)\sqrt{\log p})$. See the supplementary
Section~\ref{sec::ZH08proof}.

\noindent{\bf $\ell_1$ error on $h$ versus on $\delta$.}
Recall $h = \beta_{\init} - \beta^{\ext}({T_0})$ and $\delta = \beta_{\init}-
\beta$. Since $\supp(\beta_{\init}) \asymp q^{*}$ by the supplementary
\eqref{eq::supploss},  these results on the $\ell_1$ error $\norm{\delta}_1$, where $\delta =
\beta_{\init} - \beta$, and the Lasso support in~\cite{ZH08},  cf. \eqref{eq::L1lossmain}
and the supplementary~\eqref{eq::supploss}, are different from 
Theorem~\ref{thm::RE-oracle} in the present paper, as we do not aim to 
bound the Lasso support directly.
A more subtle point is that since $q^* \asymp \size{A_1}$,
we can crudely interpret $q^*$ as the size
of $\supp(\beta_{\init})$, given $q \ll q^*$; 

In our theory, indeed, the Lasso support can be much larger than 
$s_0$ while the $\ell_1$-, $\ell_2$-norm bounds on $h$ are quite tight in the
sense of Theorem~\ref{thm::RE-oracle}.
Hence the SRC assumptions~\eqref{eq::SRC}  are somewhat similar and
closely related to the sparse eigenvalue conditions in
Theorem~\ref{thm::RE-oracle}, where $\RE(s_0, 4, X)$ also holds with
$K(s_0, 4) <\infty$.
On the one hand, the SRC condition is set at a rank $q^* \asymp 
\size{\supp(\beta_{\init})}$, which may be significantly larger than 
$s_0$ in our model.
On the other hand, the $\ell_2$ and $\ell_1$ error on $h$ in 
Theorem~\ref{thm::RE-oracle} depend on $s_0$, as well as $K^2(s_0, 
4)$,  $\Lambda_{\max}(2s_0)$,  and  $\Lambda_{\min}(2s_0)$ (cf.~\eqref{eq::orthocauchy}),
and $\Lambda_{\max}(s-s_0) \asymp \Lambda_{\max}(\log
p(s_0-a_0)/c_t^2)$, without an explicit condition 
on~\eqref{eq::eta2def}.

We demonstrate the tightness of these theoretical bounds in 
Figures ~\ref{fig:lasso-esimate-L1-b2-b0-l2-b1-error}
and~\ref{fig:model-size-var-s}.
Consider~\eqref{eq::tigermodel} in particular.
In our experiments, $c_t$ is an absolute constant $\in [0.527, 1.658]$
(See Table~\ref{tab:s0-config}), and the signs of the non-zero values 
in $\beta^{(2)}$ are chosen from $\{\pm 1\}$ at random.
Hence, we have a longer tail (in the sense of a larger
$\norm{\beta^{(2)}}_1$ with many small coordinates) and we expect the
Lasso estimate to have a larger $\norm{\delta}_1$. 
See Table~\ref{tab:s0-config}, where we fix $\sigma =1$, 
$a_0=30$, $s_0=50$, and $C_m=1/\sqrt{2}$ for $s> s_0$.

As shown in the bottom two plots in
Figure~\ref{fig:lasso-esimate-L1-b2-b0-l2-b1-error},
where we set $(p, n, s_0, a_0, \gamma, C_a) = (2048, 1600, 50, 30, 0.7, 1.706)$,
we observe that
$\norm{\delta}_1$ dominates $\norm{h}_1$ consistently across all $s \in
\{130, 370, 511\}$ as the Lasso penalty $\lambda_n = f_p
\lambda \sigma$ varies.
In particular while upon rescaling, all curves
corresponding to different values of $s$
align well for the $\ell_2$-error, namely, $\twonorm{\delta}$,
$\twonorm{h_{T_0}}$ (right plot) and the $\ell_1$-norm $\norm{h}_1$,
$\norm{h_{T_0}}_1$, and $\norm{h_{T_0^c}}_1$ (left plot), the same is
not true for $\norm{\delta_{T_0^c}}_1$ and $\norm{\delta}_1$.
However, we have by the triangle inequality,
\ben
\label{eq::gap1}
\abs{\norm{\delta}_1 - \norm{h}_1} & = & 
\abs{\norm{\delta_{T_0^c}}_1 - \norm{h_{T_0^c}}_1} \le 
\norm{\beta_{T_0^c}}_1.
\een
We show the lower bound on $\norm{\beta_{T_0^c}}_1$ in~\eqref{eq::lowereta1a}.
Hence the error $\norm{\delta_{T_0^c}}_1$ increases as $s$ increases (as $c_t$
decreases) but also remains bounded, as predicted by
\eqref{eq::eta1c},~\eqref{eq::eta1d} and \eqref{eq::gap1}.
See Sections~\ref{subsec:exp-setup}.

\subsection{Variable selection in $A_0$}
\begin{figure}
\begin{center}
\begin{tabular}{ccc}
\includegraphics[width=0.39\textwidth]{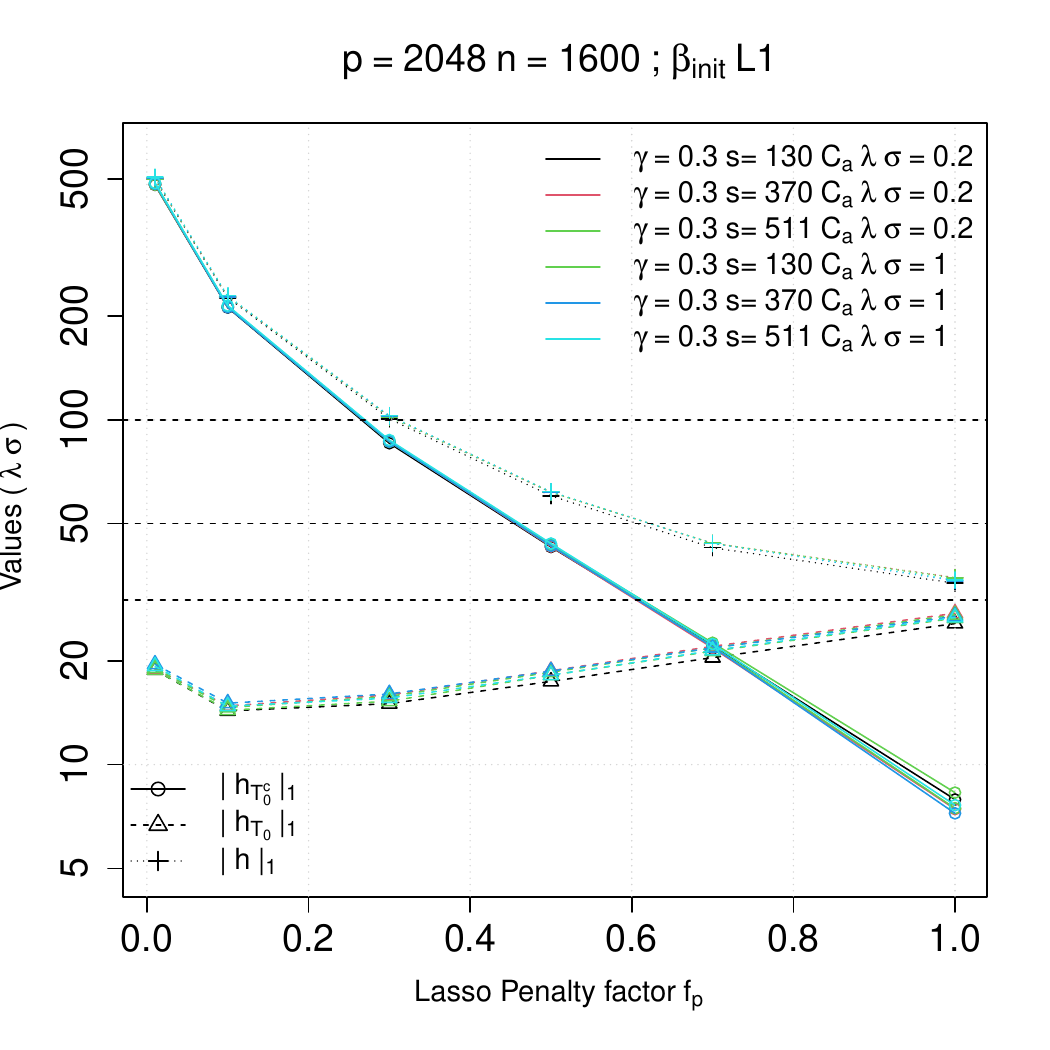} &
\includegraphics[width=0.39\textwidth]{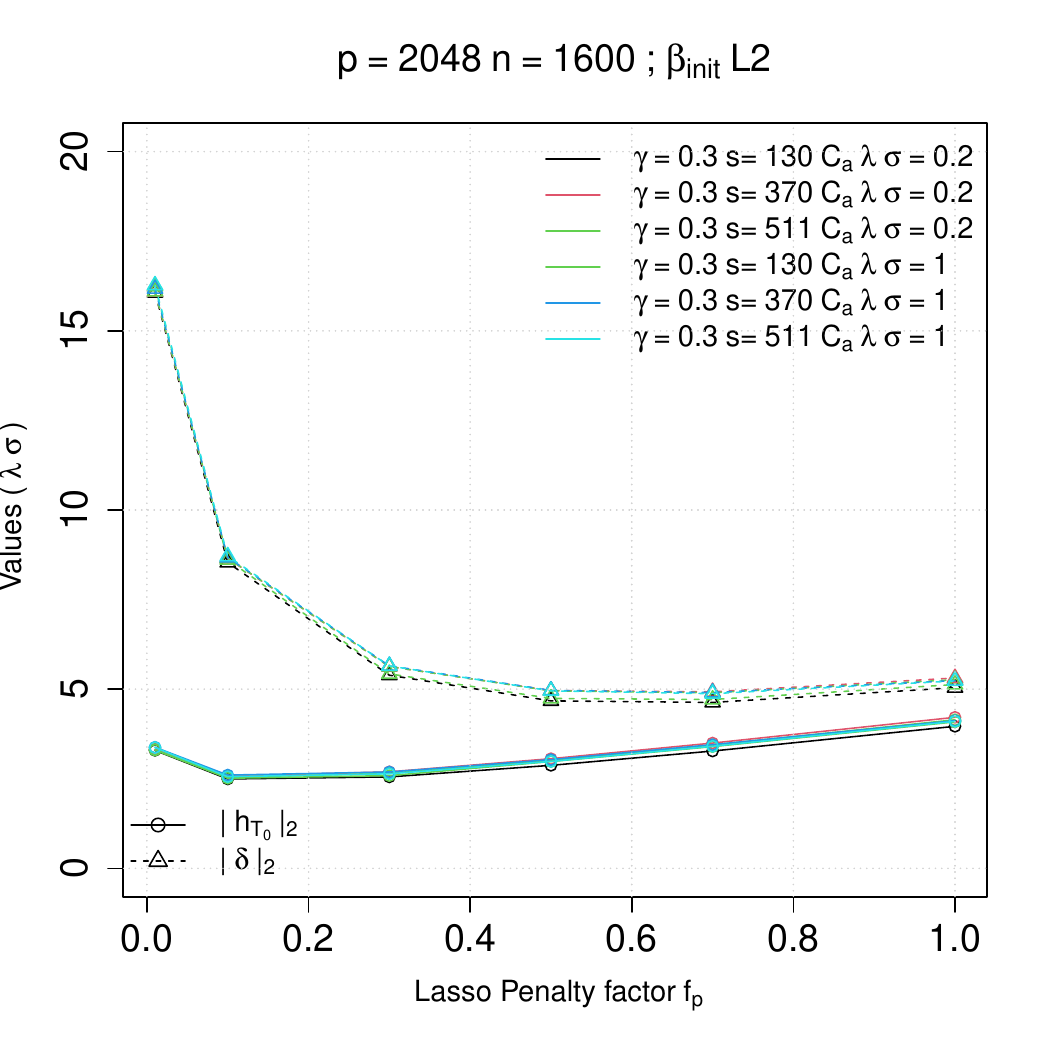}\\ 
\includegraphics[width=0.39\textwidth]{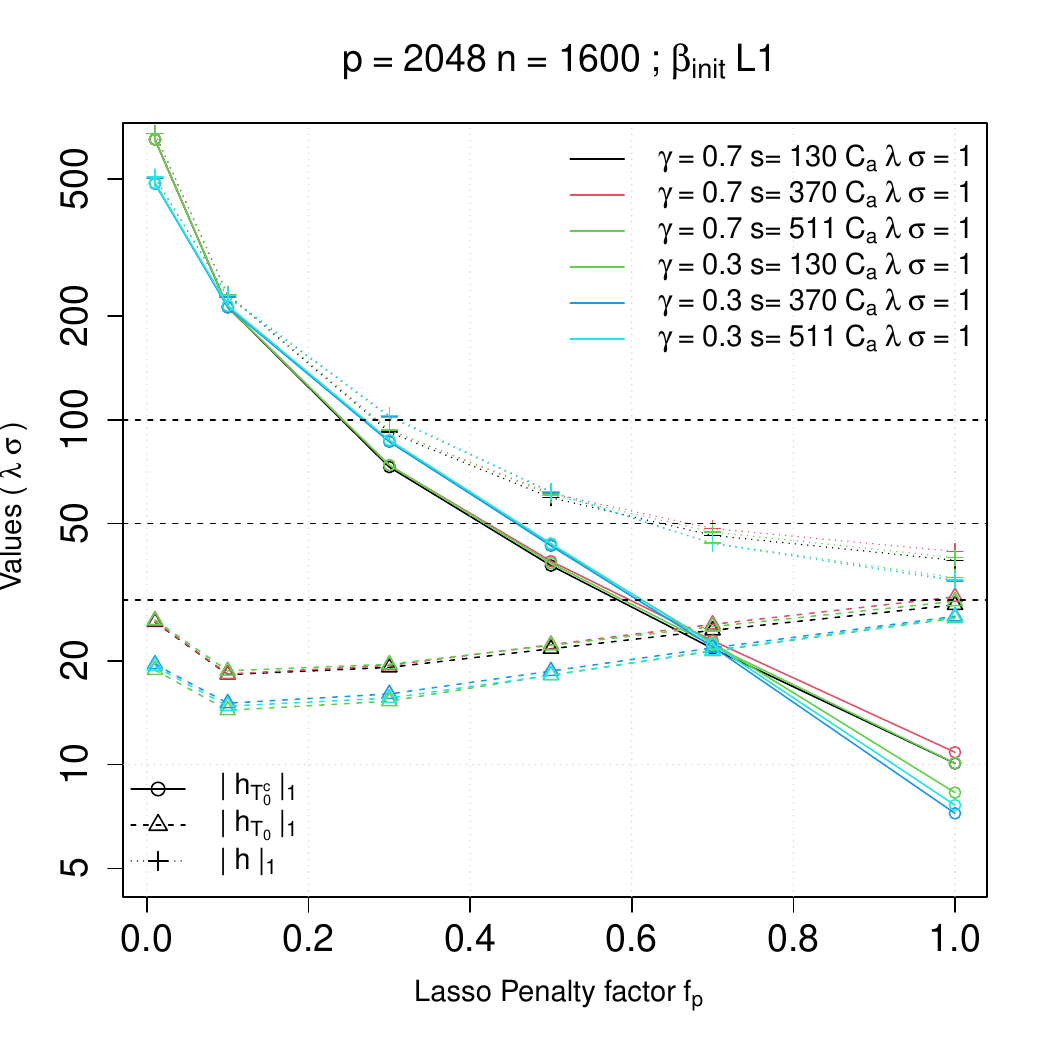} &
\includegraphics[width=0.39\textwidth]{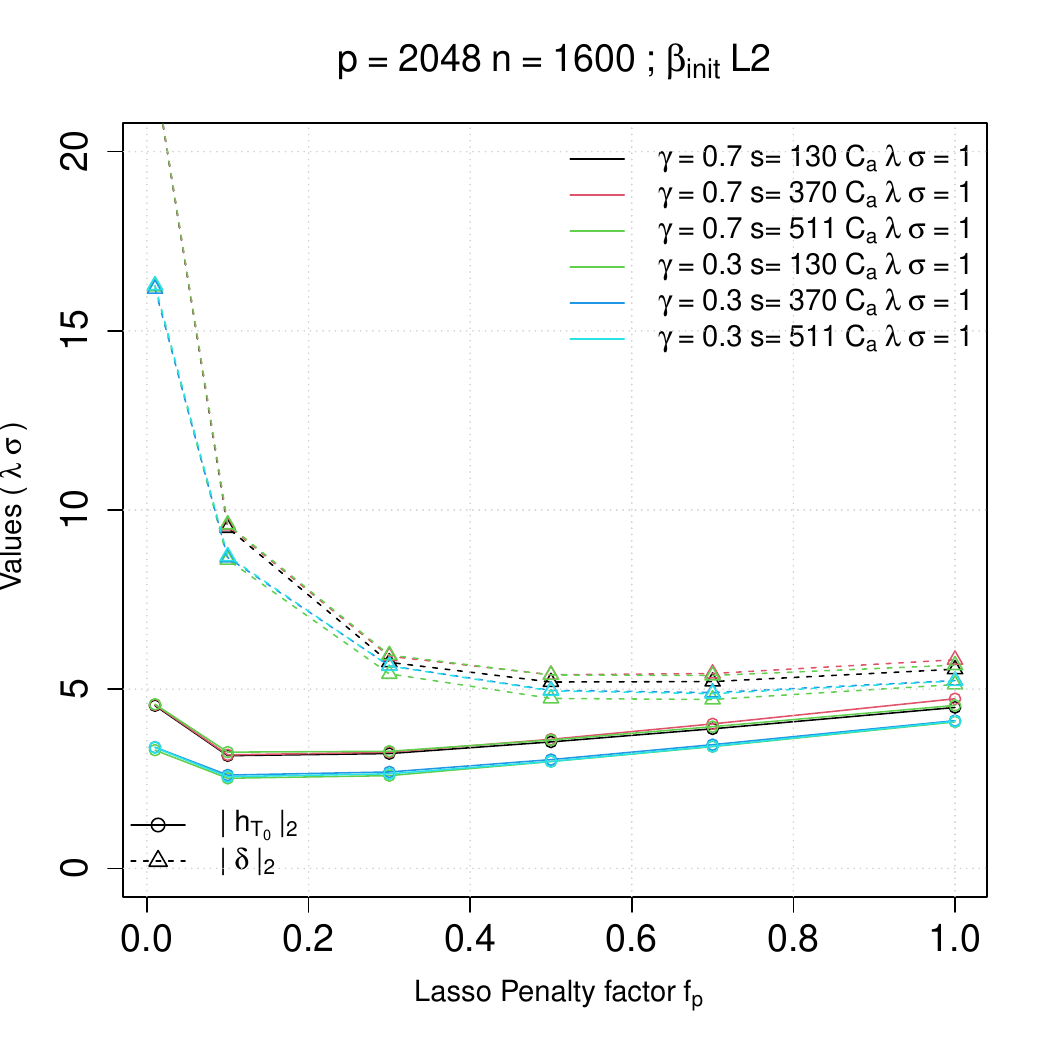} \\
\includegraphics[width=0.39\textwidth]{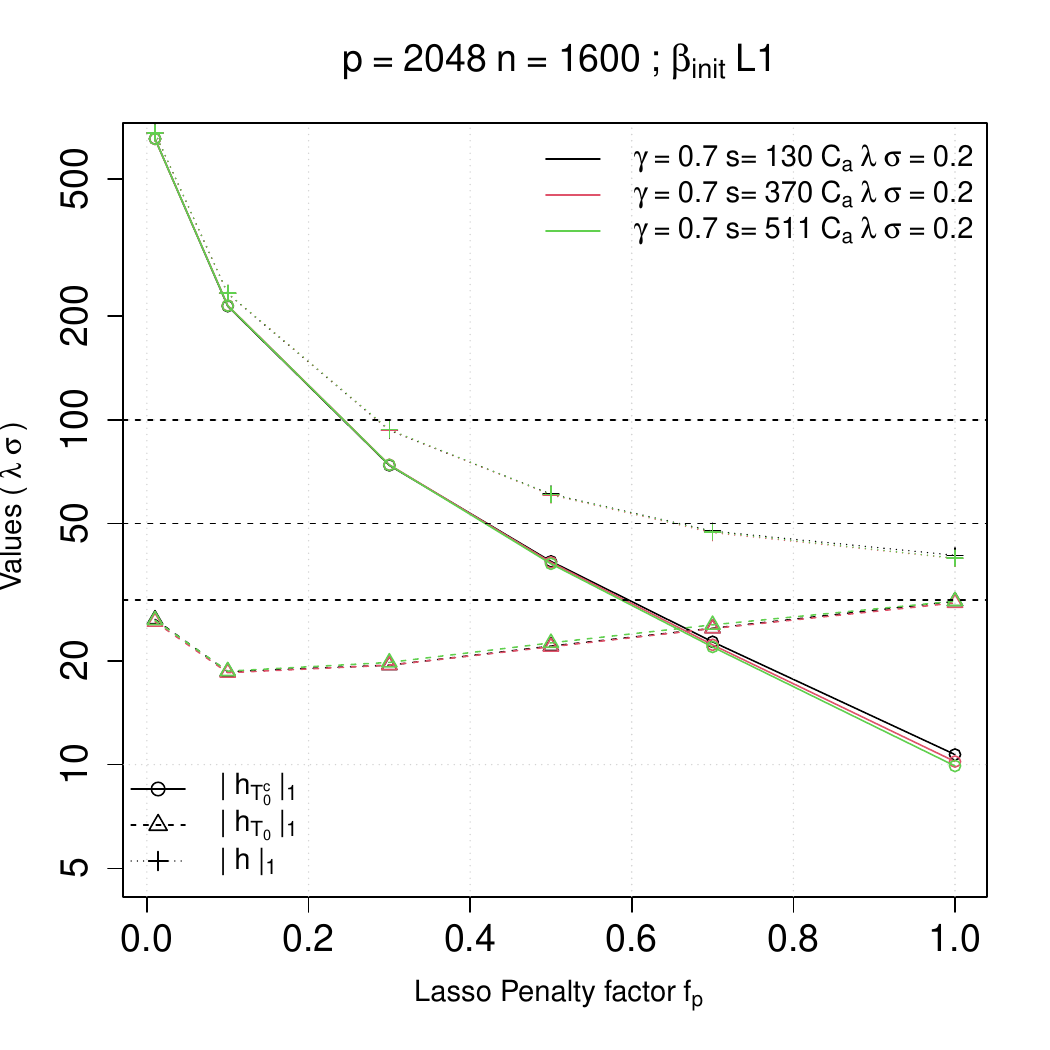} &
\includegraphics[width=0.39\textwidth]{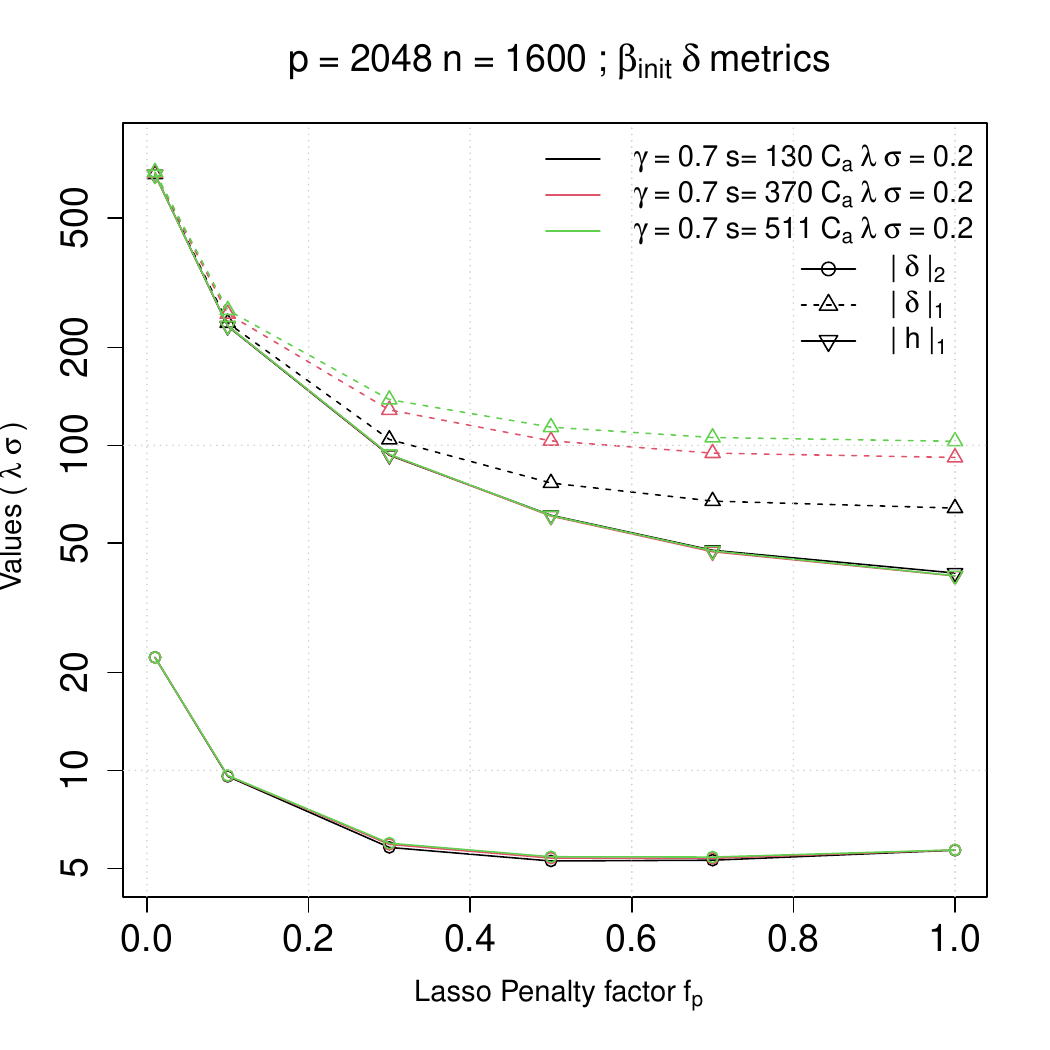} \\
\end{tabular}
\caption{
$p=2048, n=1600$. Left column: $\norm{h_{T_0^c}}_1$, $\norm{h_{T_0}}_1$, and $\norm{h}_1$ as Lasso penalty ($f_p$) increases across different sparsity $s \in \{130, 370, 511\}$. Right column: plots of  $\norm{h_{T_0}}_2$ and $\norm{\delta}_2$.
In the top panel, we fix $\gamma = 0.3$, and compare two cases of $C_a \lambda \sigma \in \{0.2, 1 \}$.  In the middle panel, we fix $C_a \lambda \sigma = 1$ and compare two cases of $\gamma \in \{0.3, 0.7\}$.  In the bottom panel, we zoom in on one case
  with $\gamma=0.7, C_a \lambda \sigma = 0.2$, and we plot
  $\norm{\delta}_1$ together with $\norm{\delta}_2$ in the bottom right panel.}
\label{fig:lasso-esimate-L1-b2-b0-l2-b1-error}
\end{center}
\end{figure}

\begin{figure}
\begin{center}
\begin{tabular}{cc}
\includegraphics[width=0.39\textwidth]{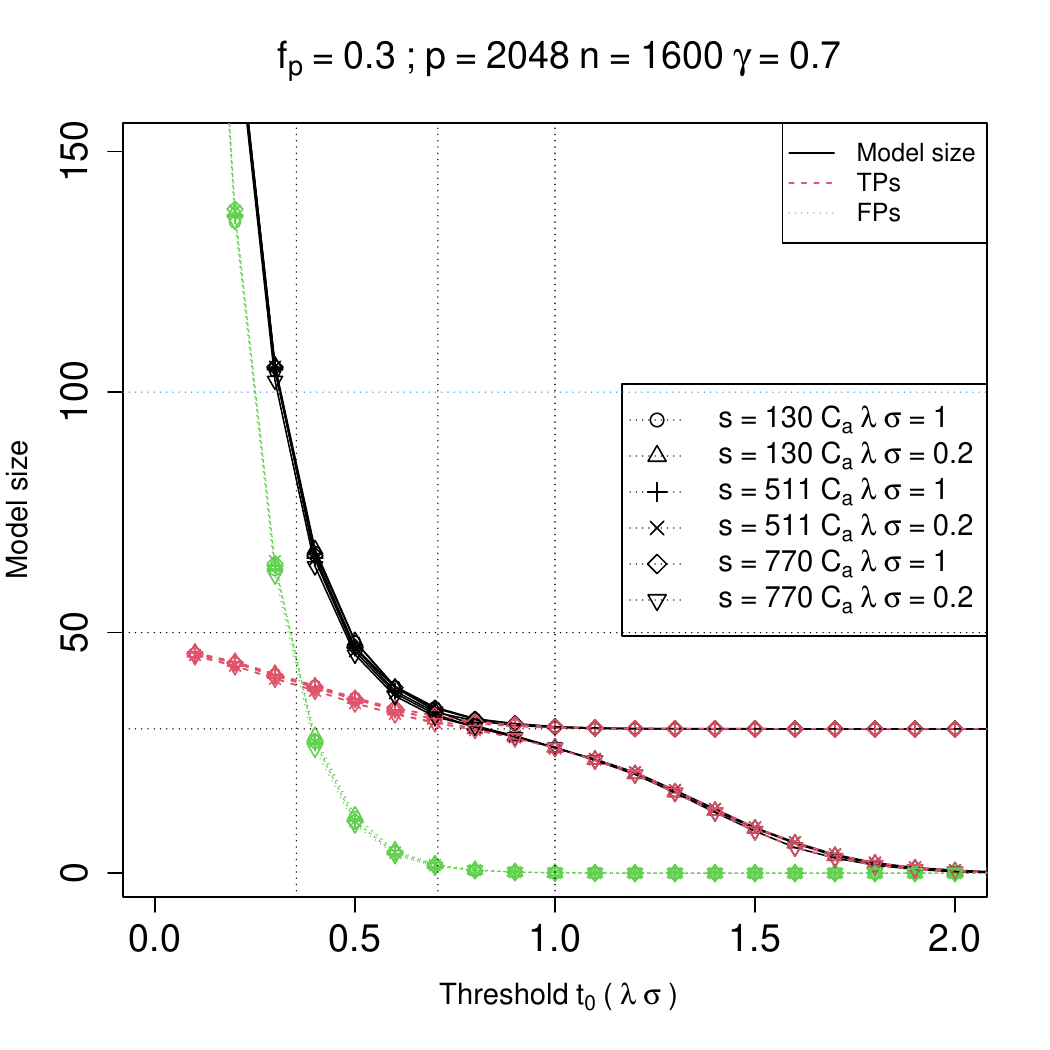} &
\includegraphics[width=0.39\textwidth]{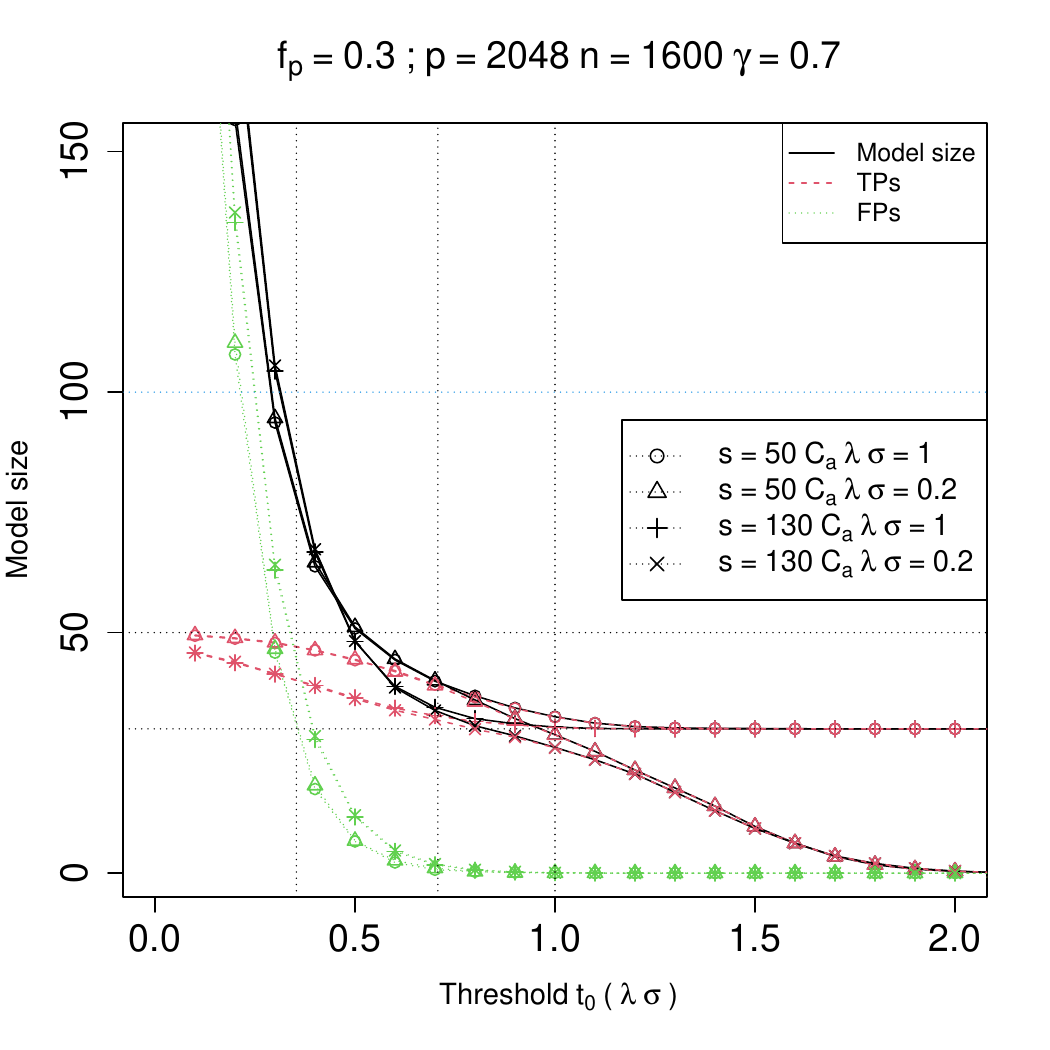}\\ 
(a) & (b) \\
\includegraphics[width=0.39\textwidth]{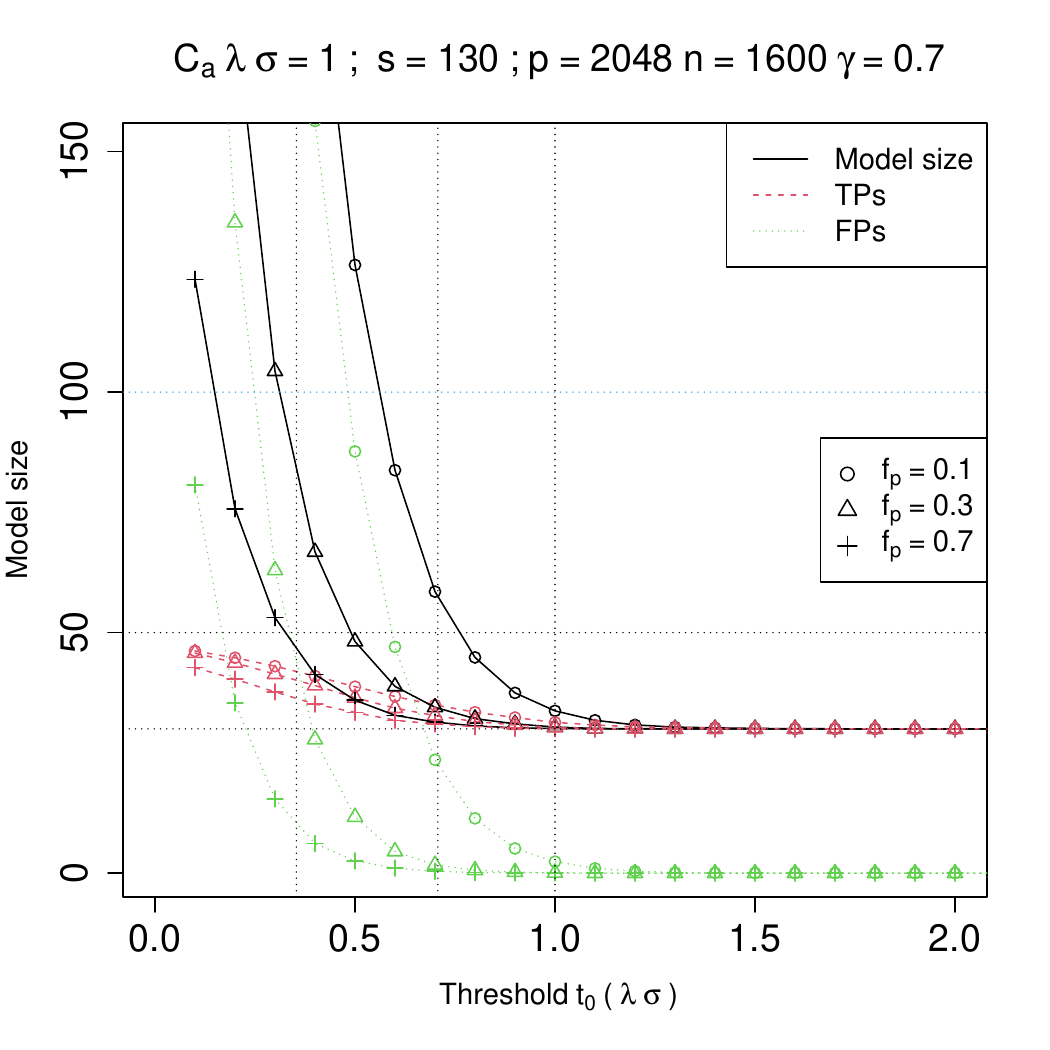}&  
\includegraphics[width=0.39\textwidth]{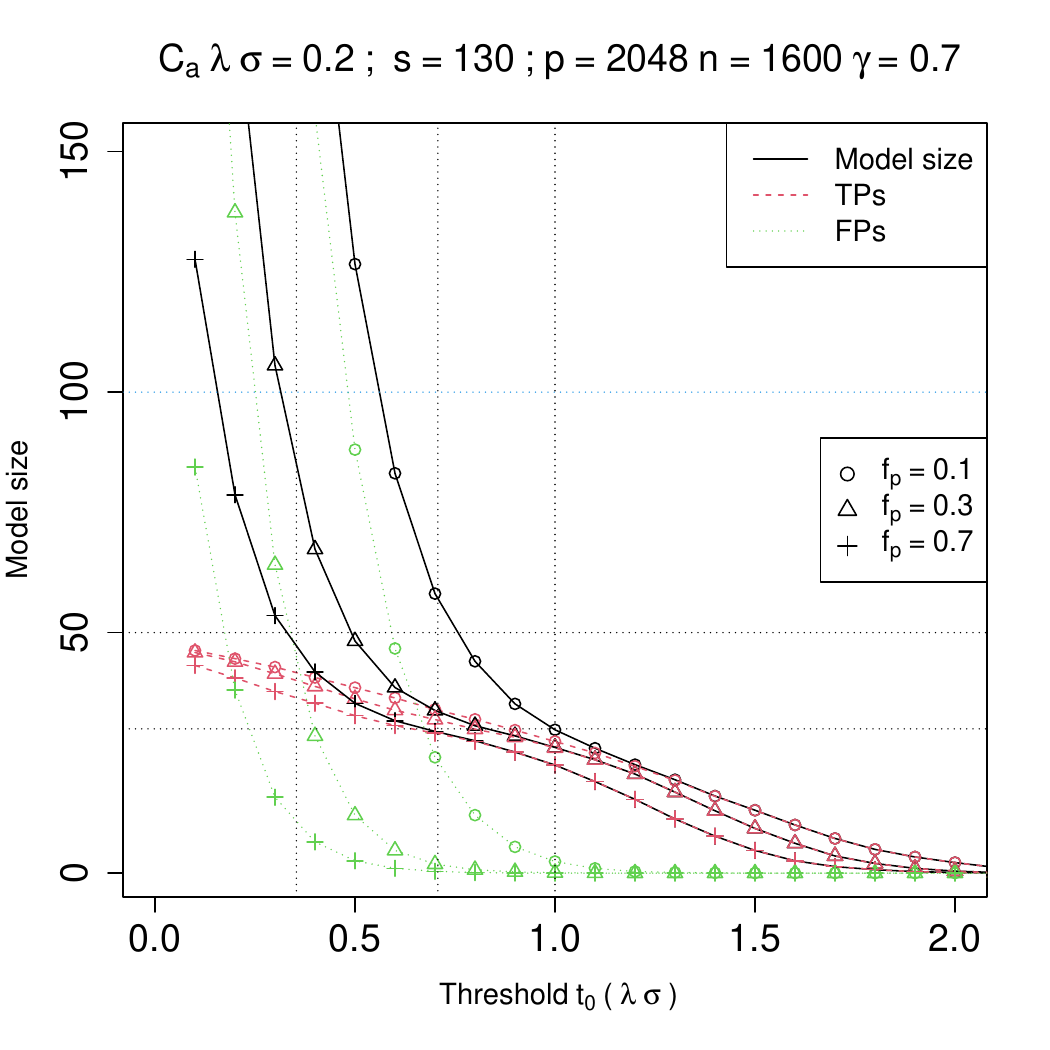}\\ 
(c) & (d) \\
\end{tabular}
\caption{
$p=2048, n=1600, \gamma=0.7$. Plots of model size ($|I|$), number of TPs and
FPs, as threshold increases.  Note $|I|=$ TPs + FPs.  In (a) and (b), Lasso
penalty factor $f_p=0.3$ is fixed, and in panel (a) $s \in \{130, 511, 710\}$,
and in panel (b) $s \in \{50, 130\}$.  In panels (c) and (d), we plot the same metrics
across different $f_p \in \{0.1, 0.3, 0.7\}$ with fixed $s=130$.  In
all panels, the 3 dotted vertical lines from left to right represent $C_m
\lambda \sigma / 2, C_m \lambda\sigma$ and $\lambda\sigma$. The model
size remains invariant and hence the diagonal  dashed lines all stay
flat for $\lambda\sigma < t_0 \le 2 \lambda \sigma$ for $\beta_{\min,
  A_0}=1$.}
\label{fig:model-size-var-s}
\end{center}
\end{figure}

Note that $\beta_{\min,A_0} = \norm{\beta^{(11)}}_\infty =C_a\lambda\sigma$, 
$\norm{\beta^{(12)}}_\infty =C_m\lambda\sigma$, and 
$\norm{\beta^{(2)}}_\infty =C_t\lambda\sigma = c_t \sigma/\sqrt{n}$.
We mention up front that no matter where we put the threshold, 
some of the signals in $\beta^{(12)}$ will be lost so long as $t_0 \asymp 
\lambda \sigma$ since $C_m \le 1$; cf. Fig.~\ref{fig:model-size-var-s}. 
In our experiments, we can consistently recover those signals in $A_0$ for $\beta_{\min,A_0} :=
\min_{j    \in A_0} |\beta_{j}| \asymp \lambda \sigma \sqrt{s_0}$ (in
case $C_a = 8.528$ for the model class in
Figure~\ref{fig:tiger-beta}), but this  is not the case when
$\beta_{\min,A_0} \asymp  \lambda \sigma$ (in case $C_a = 1.706$).
This difference occurs despite the
nearly identical $\ell_p, p=1, 2$ norm bounds on the estimation error
for the initial Lasso estimator $\beta_{\init}$ for the two models of $\beta$; 
See the top two panels in Figure~\ref{fig:lasso-esimate-L1-b2-b0-l2-b1-error},
where curves corresponding
to different values of $\beta_{\min, A_0}$ with $C_a \lambda \sigma =0.2, 1$
align well across different values of $s \in \{130, 370, 511\}$;
This is true for
$\twonorm{h_{T_0}}$, $\norm{h_{T_0}}_1$, $\norm{h_{T_0^c}}_1$, and 
$\twonorm{\delta}$ under the same design matrix $X$ of dimension $1600 \times 2048$.
 The relative effect on variable selection in component $\beta^{(12)}, 
 \beta^{(2)}$, and $\beta^{(0)}$ follows the same trend for both
 $C_a$ settings.
The relative effect on  variable selection in component $\beta^{(11)}$
is much more significant in case $\beta_{\min, A_0}  = 0.2$ ($C_a =
1.706$), since the largest error magnitude may reach the full signal
strength,  while for $\beta_{\min,  A_0}  = 1$ ($C_a = 8.528$),  the
magnitude of the  error is only a fraction $\propto 1/{\sqrt{s_0}}$ of the signal  strength in $A_0$ and hence $A_0 \subset I$ so long as $t_0 =o(\lambda \sigma \sqrt{s_0})$.

\begin{table}
    \caption{Evaluation metrics for variable selection. $I$ is the
      estimated model, and $\drop =I^c$.} 
  \begin{center}
\begin{tabular}{c | c | c}
  Metric & Definition  & Meaning \\ \hline
  TP & $I \cap T_0$    & Selected variables in $T_0$ \\
  FP & $I \cap T_0^c$  & Selected variables in $T_0^c$ \\
  TN & $\drop \cap T_0^c$  & Not selected variables in $T_0^c$  \\
  FN & $\drop \cap T_0$    & Not selected variables in $T_0$ \\ \hline
\end{tabular}
    \label{tab:tp-def}
  \end{center}
\end{table}
\section{Numerical results}
\label{sec::numeric}
In this section, we present results from numerical simulations to
validate the theoretical analysis presented in previous sections.
We consider Gaussian random matrices for the design $X$ with both
$p \times p$ identity and Toeplitz covariance.
We refer to the former as i.i.d. Gaussian ensemble, where $X_{ij} \sim N(0, 1)$ for all $i, j$, and the 
latter as Toeplitz ensemble, where the covariance matrix $T(\gamma)$
for each row vector in $\R^p$ is given by $[T(\gamma)]_{i,j} =
\gamma^{|i-j|}$.
The design matrix dimensions are either $(p=1024, n=800)$ or $(p=2048, 
n = 1600)$. 
To evaluate the impact of nominal sparsity $s$ on the 
recovery of $s_0$ components, we 
use $\beta$ as constructed as in Section \ref{sec::newmodel}. 
We have the following steps.

\begin{enumerate}
\item
Generate input $\beta \in \R^p$ as shown in Fig.~\ref{fig:tiger-beta}. $\beta$ is determined by
the parameters $(C_a, C_m, C_t)$, $\lambda = \sqrt{2 \log p/n}$,
and $\sigma$ in the noise $\epsilon$. Here,
we fix $a_0=30, s_0=50$, and $\sigma=1$.
The signs and positions of the non-zero coordinates are chosen at
random. See Table~\ref{tab:s0-config}.
\item
 Generate a {Gaussian ensemble} $X_{n \times p}$ with independent rows,
  which is then normalized to have column $\ell_2$-norm $\sqrt{n}$.
We consider two types of design: i.i.d. Gaussian ensemble, 
and Toeplitz ensemble as mentioned above.
\item
Compute $Y = X \beta + \epsilon$, where the noise $\epsilon \sim N(0, \sigma^2
I_n)$ is generated with $I_n$ being the  $n \times n$ identity matrix. 
\item
Feed $Y$ and $X$ to the Thresholded Lasso algorithm to estimate $\beta$ as described in
Section~\ref{sec:introduction}. We call the Lasso procedure $\lars(Y, X)$ ~\cite{EHJ04} to
compute the full regularization path.
We select the $\beta_{\init}$ from this output  
path with penalty parameter $\lambda_n = f_p \lambda \sigma$.
We then threshold  $\beta_{\init}$ with threshold $t_0 = f_t \lambda \sigma$,
and run OLS to obtain $\hat{\beta}^{\OLS}(I)$.
\end{enumerate}

We set $C_m =1$ for $s=s_0$ so that $\twonorm{\beta^{(12)}} = (s_0 
-a_0) \lambda \sigma$ and $\twonorm{\beta^{(2)}} = 0$.
For $s>s_0$, we fix $C_m= 1/\sqrt{2}$, 
and set $\norm{\beta^{(12)}}_2 = \norm{\beta^{(2)}}_2 =
\frac{1}{\sqrt{2}}\lambda\sigma \sqrt{s_0 - a_0}$ in 
\eqref{eq::tigermodel} and  \eqref{eq::beta2supp}.
Hence the upper bound in~\eqref{eq::SR-range-intro} becomes an equality for both scenarios.
When we lower $C_t$ ($c_t$) in~\eqref{eq::CT}, we will have a $\beta$ with a 
larger $\supp(\beta)$ with many small coefficients.
In particular, $s - s_0$ (the length of 
$\beta^{(2)}$) and $\norm{\beta^{(2)}}_1$ increase as 
$c_t$ decreases, but $\norm{\beta^{(2)}}_2$ remains the same; 
cf.  \eqref{eq::beta2supp} and \eqref{eq::eta1c}.
There are two main tuning parameters:
$\lambda_n = f_p \lambda \sigma$ and $t_0 = f_t \lambda \sigma$, where
$f_p, f_t >0$.
For each experiment, after we generate $\beta \in \R^p$ in Step 1, we
repeat Steps 2-4 100 times, and compute averages after 100 runs.  Due
to limited space, we only present results from experiments using
Toeplitz ensemble with
$\gamma \in \{0.3, 0.7\}$, but we observe similar trends for other design
matrices such as the i.i.d. Gaussian ensemble.
In the present context, we adopt a more stringent definition of metrics for variable selection
evaluation. More specifically, we define True
positives (TPs) as those variables from $I \cap T_0$,  and False positives (FPs)
refer to variables in $I \cap T_0^c$.  Note that this interpretation naturally
flags many more variables as FPs, which in the conventional notion would have
been counted as TPs. False negatives (FNs) refer to variables from
$\drop \cap T_0$ where $\drop = I^c$. True negatives (TNs) refer to variables
from $\drop \cap T_0^c$. See Table~\ref{tab:tp-def}.


\subsection{$\ell_1$ and $\ell_2$ error bounds for  $\beta_{\init}$}
\label{subsec:exp-setup} 
Denote by $\beta_{A}$ the restriction of $\beta$ to the set 
$A \subset [p]$, with all other coordinates set to zero.
We also use $\beta^{(1)} := \beta^{(11)} + \beta^{(12)} = \beta_{T_0}$ and 
$\beta^{(2)} + \beta^{(0)} = \beta_{T_0^c}$ throughout our 
discussion.
Recall $h = \beta_{\init} - \beta_{T_0}$ and $\delta = \beta_{\init}- \beta$. 
First, we investigate $\ell_1$ and $\ell_2$ error bounds for
$\delta$ and $h$, with the Lasso penalty being $\lambda_n = f_p \lambda \sigma$.
In our experiments, we vary the correlation parameter $\gamma$,
$\beta_{\min,A_0}$ and nominal sparsity $s$ under model class~\eqref{eq::CT}.

In the top and middle rows of
Fig.~\ref{fig:lasso-esimate-L1-b2-b0-l2-b1-error}, we plot
$\norm{h_{T_0^c}}_1$, $\norm{h_{T_0}}_1$, and $\norm{h}_1$ in the left column,
and $\norm{h_{T_0}}_2$ and $\norm{\delta}_2$ in the right column.
Across all plots, we observe that curves for both $\norm{h_{T_0}}_1$ and $\norm{h_{T_0}}_2$ decrease slightly first
and then increase as $f_p$ increases, due to some variables with
significant coefficients being eliminated.
We observe that as $f_p$ increases, $\norm{h_{T_0^c}}_1$ decreases quickly for
all nominal sparsity $s$ ($y$-axis is in log scale).  
Even though values in $\beta^{(2)}$ are non-zero, due to
their small magnitude, they are essentially treated the same way as the zeros
in $\beta^{(0)}$ (as they should be).
Moreover, all curves for $\norm{h_{T_0^c}}_1$ and $\norm{h}_1$
align well for the same $\gamma \in \{0.3, 0.7\}$ as predicted by Lemma~\ref{lemma::HT01case3} and
Theorem~\ref{thm::RE-oracle} since we fix $s_0$ while varying sparsity
$s \in \{130, 370, 511\}$.

\noindent{\bf Dependence on $\gamma$.}
In the top row, all curves are closely aligned, validating that
$\beta_{\min,A_0}$ (for different values of $C_a$) does not impact these error metrics significantly.
In contrast, in the middle row, we observe that 
curves for $\norm{h_{T_0}}_1$ and $\norm{h_{T_0}}_2$, 
shift downwards slightly for $\gamma = 0.3$.
This is expected since for the design matrix $X$ with 
smaller $\gamma$, the incoherence parameters appearing in 
Theorem~\ref{thm::RE-oracle} will be smaller.
The two sets of curves for $\norm{h_{T_0^c}}_1$ and $\norm{h}_1$ cross each other 
for $\gamma =0.3$ and $0.7$ with a small 
gap (middle left panel), but eventually the set with a larger $\gamma$
stays on the top.
The gap is potentially caused by the non-linear interactions between
$\gamma$ and the penalty parameters throughout the entire Lasso path.

\noindent{\bf $\ell_1$ and $\ell_2$ error for $\delta$.}
In the right column, for the $\ell_2$ error for estimating $\beta$
with $\beta_{\init}$, we observe the typical V-shaped curves as
$f_p$ increases from $0$ to $1$, since $\norm{\delta}_2$ reaches a minimum and then increases
again as the penalty $\lambda_n$ increases.
In the bottom right panel of Fig.~\ref{fig:lasso-esimate-L1-b2-b0-l2-b1-error}, we
plot $\norm{\delta}_1$ and $\norm{\delta}_2$, and also $\norm{h}_1$
(middle solid curves) for reference.
All solid curves in the bottom corresponding to the same pair of $(s_0, \gamma)$ again 
align well as predicted by Theorem~\ref{thm::RE-oracle}  (bottom right 
panel). In the same panel, we observe the three dashed curves corresponding to $\norm{\delta}_1$ 
are clearly separated under different sparsity $s \in \{130, 370, 511\}$, 
and stacked in descending order as $s$ decreases, that is,
when $\beta$ becomes more sparse in~\eqref{eq::gap1}.
We know $\norm{\delta}_1 = \norm{\delta_{T_0}}_1 +
\norm{\delta_{T_0^c}}_1$, where $\norm{\delta_{T_0}}_1 =
\norm{h_{T_0}}_1$.
In contrast, as shown in the left panels, $\norm{h_{T_0}}_1$,
$\norm{h_{T_0^c}}_1$ and hence $\norm{h}_1$ all align well across
different $s$ for the same $\gamma$.
Hence, the difference in $\norm{\delta}_1$ is due to the component 
$\norm{\delta_{T_0^c}}_1$, which depends on the sparsity.
This phenomenon is  expected and explained in
Section~\ref{sec::tightness}.
The influence of $\gamma$ on $\norm{\delta}_2$ follows a
similar trend as that for  $\norm{h_{T_0}}_2$.
Hence, all error curves corresponding to a larger $\gamma$
($0.7$) consistently dominate those with a smaller $\gamma = 0.3$ for 
$\norm{h_{T_0}}_1$, $\norm{h_{T_0}}_2$, and $\norm{\delta}_2$ in the middle panel.

\subsection{Variable selection with thresholding}
\label{sec::VST}
In Figure~\ref{fig:model-size-var-s}, we plot in all panels the
final model size $\abs{I}$ (top right solid black curves), the number of TPs
in model $I$ ($\size{I \cap T_0}$, middle diagonal red dashed lines/curves), and
the number of FPs from $T_0^c$, ($\size{I \cap T_0^c}$, left bottom
dotted green curves), as functions of threshold $t_0$ for
$\beta_{\min, A_0} =C_a \lambda \sigma \in \{0.2, 1\}$ and $C_a \in
\{1.706, 8.528\}$.
In panel (a), all curves align well across different $s \in \{130, 370,
511\}$. At $t_0 =  C_m \lambda \sigma$ (middle dotted vertical line), the model
size is only slightly above $\abs{A_0} =30$ for both cases of $C_a$.
When the top curves (plotting $\size{I}$) touch upon or cross over the horizontal 
line of $y=30$ at $t_0 = t'$, where $t' \in (\lambda \sigma/2,  \lambda 
\sigma)$, the model contains $A_0$ exactly.
For $\beta_{\min, A_0} = 1$, all coordinates in $\beta^{(11)}$ remain in 
model $I$ so long as $t_0 \le 2 \lambda \sigma$.
For $\beta_{\min, A_0} =0.2$ ($C_a  = 1.706$), the model size
$\abs{I}$ continues to shrink till $0$, as $t_0$ increases.

In Fig. ~\ref{fig:model-size-var-s}
panel (b),  we compare the exact $s_0$-sparse  case ($s=50$) with the
almost $s_0$-sparse case ($s=130$). Under the same $t_0$, we observe
that the exact sparse case recovers more non-zero components from
$T_0$ (higher red dashed lines  with more TPs and fewer FNs for $t_0 < 1.5\lambda\sigma$) and less from $T_0^c$ (lower green lines with fewer FPs for $t_0 < \lambda\sigma$), since for $s=s_0$, $\beta_j^{(12)}$ is larger in magnitude with $C_m =1$ and $\beta_{T_0^c} =0$. 
Panels (c) and (d) show that as the Lasso penalty ($f_p$) increases,
model sizes further decrease, and hence $\abs{I \cap T_0}$ and $\abs{I
  \cap T_0^c}$ both decrease. Hence all three sets of curves ($|I|$, TPs, FPs) shift downward as
$f_p$ increases since Lasso is able to remove some less significant variables as an initial
estimator. However, FPs remain at a high level without thresholding or when the
threshold is small. This is true for both $C_a$ settings.

\noindent{\bf False negatives.}
Recall FNs = $s_0$ - TPs.  The primary distinction between the two settings of $C_a$ lies
in the FNs from $\beta^{(11)}$, since we will lose some variables
from $\beta^{(12)}$ inevitably for both choices of $C_a$, no matter where we
put $t_0$, as $\beta^{(12)}_{\init, j}$ may fall within the range of
$\pm c \lambda \sigma$ for any $c \in (0, 1]$.
The entries of  the Lasso estimate $\beta_{\init}^{(12)}$ of
$\beta^{(12)}$ are indeed spread across the interval of $[-\lambda
\sigma, \lambda \sigma]$ since $\abs{\beta_{\init, j}^{(12)}} \le
\lambda \sigma$. This is indicated in Fig.~\ref{fig:model-size-var-s}
by the negative slope in the red diagonal lines, where coordinates in
$\beta^{(12)}$ are regularly cut as $t_0$ increases.
Larger $\beta_{\min, A_0}$ means variables in $A_0$
will be kept over a longer range of $t_0$, as the remaining TPs are
all from $A_0$ after a certain threshold.
Indeed, for $\beta_{\min, A_0}=1$ ($C_a = 8.528$), we observe in
Fig.~\ref{fig:model-size-var-s} panel (c), the curves for TPs (red
dashed slanted lines) have a changing slope at around $t_0 = t'$,
where $t' \in (\lambda \sigma/2, \lambda \sigma]$ for $f_p = 0.3$, and
then flattens out along the horizontal line of $y=30$ until $t_0 = 4
\lambda \sigma$.
For $\beta_{\min, A_0}  =0.2$ ($C_a = 1.706$, panel (d)),
the dashed diagonal line intersects the horizontal line of $y=30$,
and continues with the downward trend until it reaches
$y=0$ while losing all true variables.

\noindent{\bf False positives.}
We observe in Fig.~\ref{fig:model-size-var-s} panel (a)
that FPs drop sharply in both $C_a$
cases, and the rate is the same for all $s \in \{130, 511, 770\}$ as $t_0$
increases, whereas TPs drop with a slower rate due to their larger estimated values. 
At $t_0 =  \lambda \sigma/2$, the model size is about but slightly
below $50$ with both FNs and FPs.
By our theory, the coordinates in $\beta^{(2)}$ with small 
coefficients are at the noise level and hence are neither guaranteed 
nor necessary to be included in the model $I$. 
Roughly speaking, the largest magnitude of the Lasso estimate
 $\beta_{\init}^{(2)}$ and $\beta_{\init}^{(0)}$ as well
 as their $\ell_{\infty}$ norm error are nearly all bounded by $\lambda \sigma$ in
 absolute values, and hence $t_0 \asymp \lambda \sigma$ is effective
 in controlling the number of variables selected from $T_0^c$ (False
 Positives, cf. Table~\ref{tab:tp-def}), as shown on the left bottom
 dotted green curves on all panels in Figure~\ref{fig:model-size-var-s}.
On each curve, the number of FPs drops quickly as
$t_0$ increases and transitions to a horizontal line at $t_0 \approx
0.8 \lambda \sigma$, as predicted by the lower bound on $t_0$
in~\eqref{eq::t0lowerbound}.

\begin{figure}
\begin{center}
\begin{tabular}{cc}
  \includegraphics[width=0.39\textwidth]{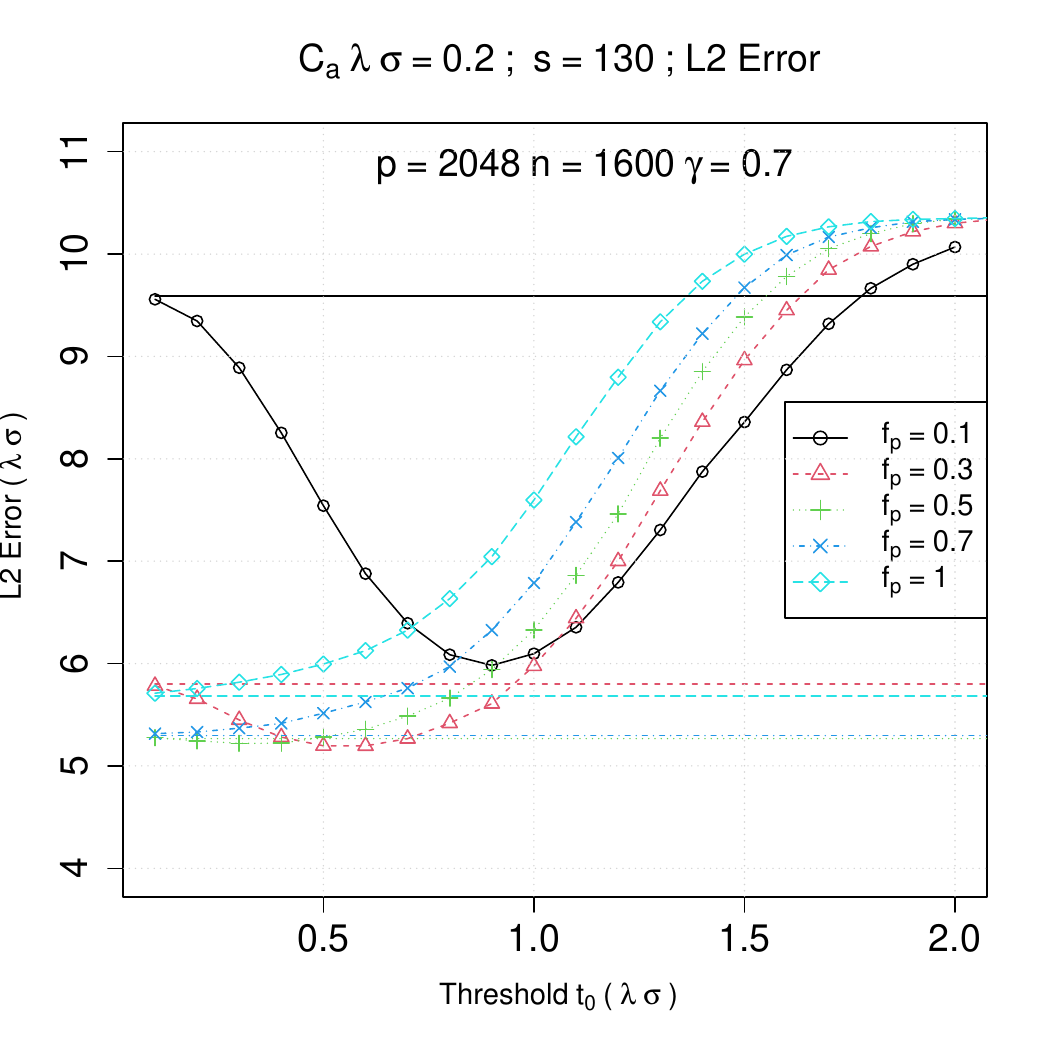} &
  \includegraphics[width=0.39\textwidth]{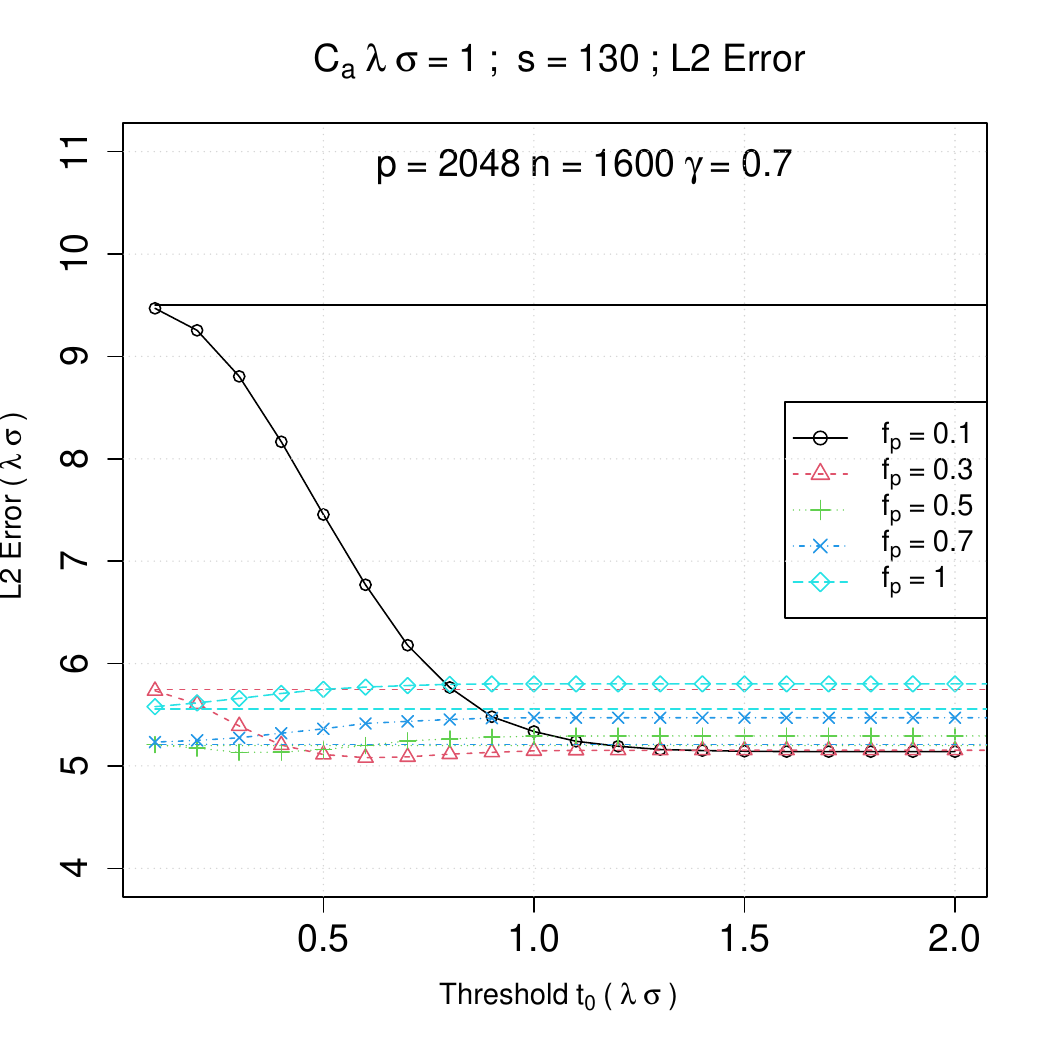}\\
\end{tabular}
  \caption{$p=2048, n=1600$, $s =130$. Plots of $\norm{\hat{\beta}^{\OLS}(I)  - \beta}_2$,
for 
  $C_a \lambda \sigma \in \{0.2, 1\}$  and $\gamma \in \{0.3, 0.7\}$. The horizontal lines
  correspond to the $\ell_2$-norm error of Lasso estimate
  $\beta_\init$, namely, $\norm{\delta}_2$.}
\label{fig:threshold-l2-errors}
\end{center}
\end{figure}

\noindent{\bf $\ell_2$-norm error of $\hat{\beta}^{\OLS}(I)$.}
In this section, we show that the previously stated tradeoffs between FPs and
TPs in Thresholded Lasso does not come at the cost of an increased
$\ell_2$ error of the final estimator.
We plot in Fig.~\ref{fig:threshold-l2-errors}, the
$\ell_2$-norm error of the final estimate, defined as
$\norm{\hat{\beta}^{\OLS}(I)  - \beta}_2$, we observe that for all the
cases of Lasso penalty ($f_p$ values), for a wide range of $t_0$,
$\norm{\hat{\beta}^{\OLS}(I)  - \beta}_2$ stays at the same level or
below $\norm{\delta}_2$. This is as predicted by Theorems \ref{thm::RE-oracle-main} and~\ref{thm::RE-oracle}.
For lower $f_p$ values, such as $0.1, 0.3$, 
$\norm{\hat{\beta}^{\OLS}(I)  - \beta}_2$ further decreases as
threshold $t_0$ increases.
This is because small coordinates in $\beta_{\init}$ with values below 
$t_0$ will be gone with thresholding, 
while variables in $A_0$ remain intact due to their larger magnitudes.
When $f_p$ increases from $0.3$ to $1$, 
$\norm{\delta}_2$ has a V-shaped curve,
due to the loss of significant coordinates in $T_0$.
This is
shown in Fig.~\ref{fig:threshold-l2-errors}, where the horizontal lines 
corresponding to $f_p = 1$ is higher than the ones correspond to $f_p=0.5$ or 
$0.7$; See also Fig.~\ref{fig:lasso-esimate-L1-b2-b0-l2-b1-error}.

\begin{landscape}
\begin{table}
\begin{tabular}{ c|c c|c c|c c c| c c c| c c c}
  \hline
  s &  $C_a$ & $C_a \lambda \sigma$ & $C_m$ & $C_m \lambda \sigma$ & $C_t$ & $C_t \lambda \sigma$ &$c_t$ &
       $\norm{\beta^{(11)}}_1$ & $\norm{\beta^{(12)}}_1$ & $\norm{\beta^{(2)}}_1$ &
       $\norm{\beta^{(11)}}_2$ & $\norm{\beta^{(12)}}_2$ & $\norm{\beta^{(2)}}_2$ \\ \hline \hline
  \multicolumn{13}{c}{$p=1024, n=800, \lambda\sigma=0.158$} \\ \hline
  50   & 6.325 & 1.0 & 1.000  & 0.158 &0.000 & 0.000  &  0.000 & 30 & 3.162 & 0.000  & 5.477 & 0.707 &0.000 \\
  130  & 6.325 & 1.0 & 0.707  & 0.112 &0.354 & 0.056  &  1.581 & 30 & 2.236 & 4.472  & 5.477 & 0.500 &0.500 \\
  370  & 6.325 & 1.0 & 0.707  & 0.112 &0.177 & 0.028  &  0.791 & 30 & 2.236 & 8.944  & 5.477 & 0.500 &0.500 \\
  511  & 6.325 & 1.0 & 0.707  & 0.112 &0.147 & 0.023  &  0.659 & 30 & 2.236 &10.738  & 5.477 & 0.500 &0.500 \\
  770  & 6.325 & 1.0 & 0.707  & 0.112 &0.118 & 0.019  &  0.527 & 30 & 2.236 &13.416  & 5.477 & 0.500 &0.500 \\ \hline
                                                                                                    
  50   & 1.265 &  0.2 & 1.000 & 0.158 &0.000 & 0.000  &   0.000 & 6 & 3.162 & 0.000  & 1.095 & 0.707 &0.000 \\
  130  & 1.265 &  0.2 & 0.707 & 0.112 &0.354 & 0.056  &   1.581&  6 & 2.236 & 4.472  & 1.095 & 0.500 &0.500 \\
  370  & 1.265 &  0.2 & 0.707 & 0.112 &0.177 & 0.028  &   0.791&  6 & 2.236 & 8.944  & 1.095 & 0.500 &0.500 \\
  511  & 1.265 & 0.2 &  0.707 & 0.112 &0.147 & 0.023  &   0.659&  6 & 2.236 &10.738  & 1.095 & 0.500 &0.500 \\
  770  & 1.265 & 0.2 &  0.707 & 0.112 &0.118 & 0.019  &   0.527&  6 & 2.236 &13.416  & 1.095 & 0.500 &0.500 \\ \hline \hline
  \multicolumn{12}{c}{$p=2048, n=1600, \lambda\sigma=0.117$} \\ \hline
  50  &8.528  & 1.0 &1.000  & 0.117 &0.000  & 0.000  &   0.000 &30 & 2.345 & 0.000  & 5.477 & 0.524 &0.000 \\
  130  &8.528  & 1.0 &0.707  & 0.083 &0.354  & 0.041  &  1.658 & 30 & 1.658 & 3.317  & 5.477 & 0.371 &0.371 \\
  370  &8.528  & 1.0 &0.707  & 0.083 &0.177  & 0.021  &  0.829 & 30 & 1.658 & 6.633  & 5.477 & 0.371 &0.371 \\
  511  &8.528  & 1.0 &0.707  & 0.083 &0.147  & 0.017  &  0.691 & 30 & 1.658 & 7.963  & 5.477 & 0.371 &0.371 \\
  770  &8.528  & 1.0 &0.707  & 0.083 &0.118  & 0.014  &  0.553 & 30 & 1.658 & 9.950  & 5.477 & 0.371 &0.371 \\ \hline
                                                                                                              
  50  &1.706  & 0.2 &1.000  & 0.117 &0.000  & 0.000  &   0.000 & 6 & 2.345 & 0.000  & 1.095 & 0.524 &0.000 \\
  130  &1.706  & 0.2 &0.707  & 0.083 &0.354  & 0.041  &  1.658  & 6 & 1.658 & 3.317  & 1.095 & 0.371 &0.371 \\
  370  &1.706  & 0.2 &0.707  & 0.083 &0.177  & 0.021  &  0.829  & 6 & 1.658 & 6.633  & 1.095 & 0.371 &0.371 \\
  511  &1.706  & 0.2 &0.707  & 0.083 &0.147  & 0.017  &  0.691  & 6 & 1.658 & 7.963  & 1.095 & 0.371 &0.371 \\
  770  &1.706  & 0.2 &0.707  & 0.083 &0.118  & 0.014  &  0.553  & 6 & 1.658 & 9.950  & 1.095 & 0.371 &0.371 \\ \hline
\end{tabular}
\caption{In this table, we list the actual values of $s$ generated
through \eqref{eq::CT} for $\beta$ configurations with $a_0=30, s_0=50.$
We also list magnitudes of each $\beta$ component
and their $\ell_1$ and $\ell_2$ norms. Here $c_t = C_t \sqrt{2 \log p}$.
The component $\beta^{(11)} = \beta_{A_0}$ has $a_0=30$ non-zero coordinates with
magnitude $C_a \lambda \sigma $, where $C_a > 1$.
The component $\beta^{(12)} = \beta_{T_0 \setminus A_0}$ has $s_0 -
a_0=20$ non-zero coordinates with magnitude $C_m \lambda \sigma$, where
$C_m = 1/{\sqrt{2}}$ for $s> s_0$. The component $\beta^{(2)} =
\beta_{S \setminus T_0}$ consists of $s - s_0$ non-zero coordinates of
magnitude $c_t \sigma/\sqrt{n}$, where $c_t$ is an absolute constant.}
\label{tab:s0-config}
\end{table}
\end{landscape}

In practice, we recommend using cross-validation to set the Lasso 
penalty $\lambda_n$, which typically ranges between $0.3\lambda\sigma$
and $0.7\lambda\sigma$ for this example, and then apply thresholding.

\section{Proofs}

\subsection{Proof of Theorem~\ref{thm::RE-oracle-main}}
\label{sec:proof-TH-main}
\begin{proofof2}
It holds by definition of $S_{\drop}$ that $I \cap S_{\drop} =
\emptyset$.
One can check via the proof of Theorem~\ref{thm::RE-oracle}
that~\eqref{eq::part-norms} holds for $D_0', D_1$ as defined
in \eqref{eq::D0-prime} and \eqref{eq::D1-define} respectively.
Hence  by Lemma~\ref{lemma:threshold-RE},  we have on event $\T_{\alpha}$
$|I| \leq 2s_0$ and $|I \cup S_{\drop}| \leq s + s_0 \leq 2s$ for $C_4 \geq D_1$; See also~\eqref{eq::IcupS}.
Moreover, $\twonorm{\beta_{\drop}} \leq \sqrt{(D'_0 + C_4)^2 + 1} 
\lambda \sigma \sqrt{s_0}.$
We have by Lemma~\ref{prop:MSE-missing-orig}, on event $\T_{\alpha}$,
for $\lambda = \sqrt{2 \log p/n}$, $\abs{I} < 2 s_0$, for $\hat{\beta}
= \hat{\beta}^{\OLS}(I)$,
\bens
\twonorm{\hat{\beta}^{\OLS}(I) - \beta}^2
& \leq & \twonorm{\beta_{\drop}}^2 \left(1 + \frac{2 \theta^2_{|I|,
|\dropS|}}{\Lambda_{\min}^2(|I|)}\right) + \frac{2|I| (1+a) \sigma^2
 \lambda^2}{\Lambda^2_{\min}(|I|)} \\
& \leq & \twonorm{\beta_{\drop}}^2
\big(1 + \frac{2 \theta^2_{\abs{S_{\drop}}, |I|}}{\Lambda_{\min}^2(2s_0)}\big) +
\frac{4 (1+ a) }{\Lambda_{\min}^2(2s_0)}  \sigma^2 \lambda^2 s_0 \\
& \leq & \lambda^2 \sigma^2 s_0 ((D'_0 + C_4)^2 + 1)
\big(  1 + \frac{2 \theta^2_{\abs{S_{\drop}},    |I|}}{\Lambda_{\min}^2(2s_0)}  + \frac{4}{9} \big),
\eens
where $\Lambda_{\min}(|I|) \ge  \Lambda_{\min}(2s_0)$
by~\eqref{eq::eigen-cond} since $\abs{I} < 2s_0$, and the
fact that ${4 (1+ a) }/{\Lambda_{\min}^2(2s_0)} \le \frac{4}{9} (D'_0)^2$ since
\bens
(D'_0)^2 & \ge & 9 d^2_0 K^4(s_0, 4) \ge 
{9 (1+a)}/{\Lambda^2_{\min}(2s_0)}, \quad  \text{ where} \\
K^4(s_0, 4)
& \ge &
K^4(s_0, 1) \ge 1/{(4 \Lambda^2_{\min}(2s_0))} \; \text{in view of~\eqref{eq::admissible2}}.
\eens
Now~\eqref{eq::MSE} clearly holds with 
$D_4^2 = ((D'_0 + C_4)^2 + 1) 
\big(1 + \frac{2 \theta_{|I|, |\dropS|}^2}{\Lambda_{\min}^2(2s_0)} +
\frac{4}{9}\big)$,
and \eqref{eq::D4-constant} holds by Lemma~\ref{lemma:parallel}.
\end{proofof2}

\subsection{Proof of Lemma~\ref{lemma:threshold-RE}}
\label{sec::proofLemmaTRE}
\begin{proofof2}
 Without loss of generality, we order $\abs{\beta_1} \ge \abs{\beta_2} \ge 
 \ldots \ge \abs{\beta_p}$. 
 Then $T_0 = \{1, \ldots,  s_0\}$.
 Let $T_1$ be the largest $s_0$
 positions of $\beta_{\init}$ outside of $T_0$.
Then
\bens
\size{I \cap T^c_0} & \le &
\frac{\norm{\beta_{\init,T^c_0}}_1}{(C_4 \lambda  \sigma)} =
\frac{\norm{h_{T^c_0}}_1}{(C_4 \lambda  \sigma)} \le s_0 D_1/C_4.
\eens
Thus $|I| = |I \cap T_0| + |I \cap T_0^c| \leq s_0 +  s_0 D_1/C_4$;
Now ~\eqref{eq::modelsize} holds since $T_0 \subseteq S$ and hence
\ben
\label{eq::IcupS}
|I  \cup S| = |S| + |I \cap \Sc| \leq s +  |I \cap T_0^c| \le s+ s_0
D_1 /C_4.
\een
We now bound $\twonorm{\beta_{\drop}}^2$.
Denote by
\bens 
\beta_{j}^{(1)} \; = \; \beta_{j} \cdot 1_{j \le s_0}\; \;
\text{ and}
\; \beta_{j}^{(2)} \; = \;  \beta_{j} \cdot 1_{j >s_0}. 
\eens
Let 
$\beta_{\drop} = \beta_{\drop}^{(1)} +  \beta_{\drop}^{(2)}$, where 
$\beta_{\drop}^{(1)} := (\beta_j)_{j \in T_0 \cap \drop}$ consists of
coefficients that are significant relative to $\lambda \sigma$, but
are dropped as $\abs{\beta_{j, \init}} < t_0$, and $\beta_{\drop}^{(2)}$
consists of those strictly below $\lambda \sigma$ in magnitude that are
dropped; cf.~\eqref{eq::beta-2-small-intro}.
Hence 
\ben
\label{eq::drop-decompose-RE}
\twonorm{\beta_{\drop}}^2 = \twonorm{\beta_{\drop}^{(1)}}^2 + 
\twonorm{\beta_{\drop}^{(2)}}^2.
\een
Now it is clear $\beta_{\drop}^{(2)}$ is bounded given~\eqref{eq::beta-2-small-intro}, 
indeed, we have for $\lambda = \sqrt{2 \log p/n}$, 
\ben
\label{eq::drop-bound-tail} 
\twonorm{\beta_{\drop}^{(2)}}^2 \leq \twonorm{\beta^{(2)}}^2 
 =  \sum_{j > s_0} \beta_j^2 = \sum_{j > s_0} \min(\beta_j^2, \lambda^2 \sigma^2)
\leq s_0 \lambda^2 \sigma^2,
\een
where the second equality is by~\eqref{eq::beta-2-small-intro} and the 
last inequality is by definition of $s_0$ as in~\eqref{eq::define-s0}.
Let $\drop_1 = \drop \cap T_0$, where $|\drop_1| <
s_0$.
Now by the triangle inequality,
\ben
\nonumber
\twonorm{\beta_{\drop}^{(1)}} =  \twonorm{\beta_{\drop_1}}
& \leq & 
\nonumber
\twonorm{\beta_{\drop_1, \init}} + 
\twonorm{\beta_{\drop_1, \init} - \beta_{\drop_1}} \\
& \leq &
\label{eq::drop-bound-RE}
t_0 \sqrt{\size{\drop_1}} +  \twonorm{h_{T_0}} \leq 
(C_4 + D_0') \lambda \sigma \sqrt{s_0},
\een
where we used the fact that 
$\twonorm{\beta_{\drop_1, \init} - \beta_{\drop_1}}=
\twonorm{h_{\drop_1}} \leq \twonorm{h_{T_0}} \le D'_0 \lambda
\sigma \sqrt{s_0}$  by~\eqref{eq::part-norms}. Hence~\eqref{eq::off-beta-norm-bound} holds
given~\eqref{eq::drop-decompose-RE},~\eqref{eq::drop-bound-tail}, and~\eqref{eq::drop-bound-RE}.
\end{proofof2}

\subsection{Proof of Proposition~\ref{PROP:COUNTING-S0}}
\label{sec::proofprop1}
\begin{proofof2}
Recall that $|\beta_j| \leq \lambda \sigma$ for all $j > a_0$ as 
defined in~\eqref{eq::betaminA0}; hence
for $\lambda = \sqrt{2 \log p/n}$, we have by~\eqref{eq::SR-range},
$\sum_{i> a_0}^p \min(\beta_i^2, \lambda^2 \sigma^2)= 
\sum_{i> a_0}^s \beta_i^2 \leq  (s_0 - a_0) \lambda^2 \sigma^2$;
hence
\bens
\size{ \{j \in A_0^c: |\beta_j| \geq  \sqrt{\log p/(c' n)} \sigma \}}
& \leq &  2 c' (s_0  - a_0) \text{ where }  \size{T_0 \setminus A_0} = s_0 - a_0.
\eens
Now given that $\abs{\beta_i} \geq \abs{\beta_j}$ for all $i \in T_0, j \in T_0^c$,
the proposition holds.
\end{proofof2}

\subsection{Proof of  Lemma~\ref{prop:MSE-missing} }
\label{sec::OLSproof-orig}
\begin{proofof2}
Recall that the random variable $\twonorm{Q}^2 \sim \chi_m^2$ is distributed
according to the chi-square distribution 
where $\twonorm{Q}^2 = \sum_{i=1}^m Q_i^2$ with $Q_i \sim N (0, 1)$ that are
independent and normally distributed. By~\cite{Joh01},
\ben
\label{eq::upper}
\prob{\frac{\chi^2_m}{m} - 1 \le -\ve} & \le &  \exp(-m\ve^2/4) \; \;
\text{ for} \; \;  0 \le \ve \le 1,\\
\label{eq::lower}
\prob{\frac{\chi^2_m}{m} - 1 \ge \ve} & \le &  \exp(-3m\ve^2/(16))  
\; \; \text{ for} \; \;  0 \le \ve \le 1/2.
\een
Although we need to bound the bad event only on one side, we provide a 
tight bound on the norm of $\twonorm{Q}$ with \eqref{eq::upper} and \eqref{eq::lower}. 
Let $\abs{I}  = m$.
Thus we have for $Q = (Q_1, \ldots, Q_{m})$ 
where $Q_j  \sim \text{i.i.d}\;  N(0, 1)$, for $\delta < 1/2$,
\bens
\lefteqn{\prob{\abs{\inv{m}\sum_{j=1}^m Q_j^2 - 1} \ge \delta} =:
  \prob{\mathcal{Q}} =}\\
&  &
\prob{\inv{m}\sum_{j=1}^m Q_j^2 - 1 \ge \delta} + 
\prob{\inv{m}\sum_{i=1}^m  Q_j^2 - 1 \le \delta} \\
& = &
\prob{\frac{\chi^2_m}{m} - 1 \ge \delta}+
\prob{\frac{\chi^2_m}{m} - 1 \le -\delta} \le \exp(-3m\delta^2/(16))  + \exp(-m\delta^2/4).
\eens
Note that $X_{I^c} \beta_{I^c} = X_{\dropS} \beta_{\dropS}$.
We have 
\bens
\nonumber
& & 
\hat\beta_{I} = (X_{I}^T X_{I})^{-1}X_{I}^T Y 
 = (X_{I}^T X_{I})^{-1} X_{I}^T (X_I \beta_I + X_{I^c} \beta_{I^c} + \epsilon)  \\
\nonumber
&  & \; \; \; \; =
\beta_I + (X_{I}^T X_{I})^{-1} X_{I}^T X_{\dropS} \beta_{\dropS} +
(X_{I}^T X_{I})^{-1}X_{I}^T \epsilon;
\eens 
Hence by the triangle inequality,
\ben
\nonumber
\twonorm{\hat{\beta}_I - \beta_{I}}  & \le &
\twonorm{
(X_{I}^T X_{I})^{-1} X_{I}^T X_{\dropS} \beta_{\dropS} +
(X_{I}^T X_{I})^{-1}X_{I}^T \epsilon}  \\
\label{eq::appen-two-terms}
& \le &
\twonorm{(X_{I}^T X_{I})^{-1} X_{I}^T X_{\dropS} \beta_{\dropS}} +
\twonorm{(X_{I}^T X_{I})^{-1}X_{I}^T \epsilon},
\een
where the two terms are bounded below as follows.

First notice that $w/\sigma = (X_{I}^T X_{I})^{-1}X_{I}^T \epsilon/\sigma$ is a mean
zero Gaussian random vector with covariance $(X_I^T X_I)^{-1}$, since
${\epsilon_i}/{\sigma} \sim N(0, 1)$ and
\bens 
\inv{\sigma^2} \E(w w^T) & = & \E  \left(  (X_{I}^T X_{I})^{-1}
  X_{I}^T [(\epsilon/\sigma) \otimes \epsilon/\sigma] X_I (X_{I}^T X_{I})^{-1}\right)  \\
& = &  (X_{I}^T X_{I})^{-1} X_{I}^T \inv{\sigma^2} \E (\epsilon \epsilon^T) X_I (X_{I}^T X_{I})^{-1}
= (X_{I}^T X_{I})^{-1}.
\eens
Then on event $\Q^c$, which holds with probability at least $1 - 2\exp(-3m/64)$
\ben
\label{eq::noise-two-norm-bound-new}
\twonorm{(X_{I}^T X_{I})^{-1}X_{I}^T \epsilon/\sigma}^2
& = &
Q^T (X_{I}^T X_{I})^{-1} Q
\leq \Lambda_{\max}\left((X_{I}^T X_{I})^{-1}\right) \twonorm{Q}^2 \\
& \leq & \frac{3m}{2\Lambda_{\min}\left(X_{I}^T X_{I}\right)}
\le \frac{3m}{2n\Lambda_{\min}(\abs{I})},
\een 
where we used an upper bound on $\twonorm{Q}^2 \le 3m/2$,
and the fact that
\bens
\Lambda_{\max}\left((X_{I}^T X_{I})^{-1} \right) =
\inv{\Lambda_{\min}\left((X_{I}^T X_{I}) \right)} =\inv{n \Lambda_{\min}\left((X_{I}^T X_{I})/n \right)} \le
\inv{n\Lambda_{\min}(\abs{I})}.
\eens
Here $\twonorm{Q}^2 \sim  \chi^2_m$.
We now focus on bounding the first term in~\eqref{eq::appen-two-terms}.
Let $P_I$ denote the orthogonal projection onto $I$.
Clearly, $I \cap S_{\drop} = \emptyset$. Let
\bens
c = (X_{I}^T X_{I})^{-1} X_{I}^T X_{\dropS} \beta_{\dropS} \quad
\quad \text{ and } \quad X_I c = P_I  X_{\dropS} \beta_{\dropS}.
\eens
By the disjointness of $I$ and $\dropS$, 
we have for $P_I  X_{\dropS} \beta_{\dropS} := X_I c$,
\ben
\nonumber
\twonorm{P_I  X_{\dropS} \beta_{\dropS}}^2 & = & 
\ip{P_I X_{\dropS} \beta_{\dropS}, X_{\dropS} {\beta_{\dropS}}} = 
\ip{X_I c, X_{\dropS} {\beta_{\dropS}}} \\
\nonumber
& \leq &
n \theta_{|I|, |\dropS|} \twonorm{c} \twonorm{\beta_{\dropS}} \le
n \theta_{|I|, |\dropS|} \twonorm{\beta_{\dropS}} \frac{\twonorm{P_I
    X_{\dropS} \beta_{\dropS}}} {\sqrt{n \Lambda_{\min}(|I|)}}, \\
\; \text{ where } \; \;
\label{eq::c-bound0}
\twonorm{c} & \leq & \frac{\twonorm{X_I c}}{\sqrt{n \Lambda_{\min}(|I|)}} 
\leq
\frac{\twonorm{P_I  X_{\dropS} \beta_{\dropS}}}
{\sqrt{n \Lambda_{\min}(|I|)}}.
\een
Hence 
\ben
\nonumber
\twonorm{P_I  X_{\dropS} \beta_{\dropS}} & \leq &
\frac{\sqrt{n} \theta_{|I|, |\dropS|}}{\sqrt{\Lambda_{\min}(|I|)}} 
\twonorm{\beta_{\dropS}}, \quad
\text { where } \quad
\twonorm{\beta_{\dropS}} = \twonorm{\beta_{\drop}} \\
\text { and }\quad
\label{eq::c-bound}
\twonorm{c} & \leq & {\theta_{|I|, |\dropS|}
\twonorm{\beta_{\drop}}}/{\Lambda_{\min}(|I|)}.
\een
Now we have on $\T_a$, by~\eqref{eq::noise-two-norm-bound-new},
\ben
\label{eq::decompose}
\twonorm{\hat{\beta}_I - \beta_{I}} 
& \leq & 
\twonorm{(X_{I}^T X_{I})^{-1} X_{I}^T X_{\dropS} \beta_{\dropS}} +
\twonorm{(X_{I}^T X_{I})^{-1}X_{I}^T \epsilon} \\
\nonumber
& \leq & 
\frac{\theta_{|I|, |\dropS|}}{\Lambda_{\min}(|I|) }\twonorm{\beta_{\drop}}
+
\frac{\sqrt{3|I| }\sigma}{\sqrt{2n\Lambda_{\min}(|I|)}}.
\een
Now the lemma holds given that $\hat\beta_I^{\OLS} = \hat\beta_I$,
$\hat\beta_{I^c}^{\OLS}  =0$, and hence
\ben
\label{eq::ploss}
\twonorm{\hat{\beta}^{\OLS} - \beta}^2 & = &
\twonorm{\hat{\beta}^{\OLS}_I - \beta_I}^2 + \twonorm{\beta_{I^c}}^2 = 
\twonorm{\hat{\beta}_I - \beta_{I}}^2 + \twonorm{\beta_{I^c}}^2\\
& \leq &
\nonumber
\frac{2 \theta^2_{|I|, |\dropS|}}{\Lambda^2_{\min}(|I|) }\twonorm{\beta_{\drop}}^2
+ \frac{3|I| \sigma^2}{n\Lambda_{\min}(|I|)} + \twonorm{\beta_{\dropS}}^2.
\een
\end{proofof2}

\begin{remark}
  We can also derive the lower bound on event $\Q^c$ for
$\size{I} =m$
\bens
\label{eq::noise-two-norm-bound-lower}
\twonorm{(X_{I}^T X_{I})^{-1}X_{I}^T \epsilon/\sigma}^2
& = &
Q^T (X_{I}^T X_{I})^{-1} Q
\ge \Lambda_{\min}\left((X_{I}^T X_{I})^{-1}\right) \twonorm{Q}^2 \\
& \ge & \frac{m}{2\Lambda_{\max}\left(X_{I}^T X_{I}\right)}
= \frac{m}{2n\Lambda_{\max}(\abs{I})},
\eens
where we used an upper bound on $\twonorm{Q}^2 \ge m/2$, $\twonorm{Q}^2 \sim  \chi^2_m$,
and the fact
that
\bens
\Lambda_{\min}\left((X_{I}^T X_{I})^{-1} \right) =
\inv{\Lambda_{\max}\left((X_{I}^T X_{I}) \right)} =\inv{n \Lambda_{\max}\left((X_{I}^T X_{I})/n \right)}=:
\inv{n\Lambda_{\max}(\abs{I})}.
\eens
Now suppose $\abs{I} \asymp 2s_0$, then
\bens
\label{eq::noise-two-norm-bound-lower}
\twonorm{(X_{I}^T X_{I})^{-1}X_{I}^T \epsilon/\sigma}^2
& \ge &  \frac{m}{2n\Lambda_{\max}(2s_0)}.
\eens
\end{remark}

\subsection{Proof of Lemma~\ref{prop:MSE-missing-orig}}
\label{sec::OLSmissing2}
\begin{proofof2}
  The only part we make modification is the following: On event $\T_{\alpha}$,
we have by~\eqref{eq::eigen-cond}, for an arbitrary set $I$ of size
$\size{I} \le 2s$,
\ben
\twonorm{(X_{I}^T X_{I})^{-1}X_{I}^T \epsilon}  \le
\twonorm{(X_{I}^T X_{I}/n)^{-1}} 
\twonorm{X_{I}^T \epsilon/n}
\label{eq::noise-two-norm-bound}
\leq \frac{\sqrt{|I|}  \basepen}{\Lambda_{\min}(|I|)}.
\een
Now we have on $\T_a$, by~\eqref{eq::noise-two-norm-bound},~\eqref{eq::c-bound0},~\eqref{eq::c-bound}, and~\eqref{eq::decompose},
\bens
\twonorm{\hat{\beta}_I - \beta_{I}} 
& \leq & 
\twonorm{(X_{I}^T X_{I})^{-1} X_{I}^T X_{\dropS} \beta_{\dropS}} +
\twonorm{(X_{I}^T X_{I})^{-1}X_{I}^T \epsilon} \\
& \leq & 
\frac{\theta_{|I|, |\dropS|}}{\Lambda_{\min}(|I|) }\twonorm{\beta_{\drop}}
+
\frac{\sqrt{|I| }\sigma \sqrt{1+a} \lambda }{\Lambda_{\min}(|I|)}.
\eens
and thus~\eqref{eq::plossintro} holds since  $\twonorm{\hat{\beta}^{\OLS}(I) - \beta}^2 =  
\twonorm{\hat{\beta}_I - \beta_{I}}^2 + \twonorm{\beta_{I^c}}^2$ by
\eqref{eq::ploss}.
\end{proofof2}

\subsection{Proof of Lemma~\ref{lemma:threshold-general}}
\label{sec::proofofBias}
\begin{proofof2}
We write $\beta = \beta^{(11)} +  \beta^{(12)} + \beta^{(2)}$ where
\bens
\beta_{j}^{(11)} =  \beta_{j} \cdot 1_{1 \leq j \leq a_0}, \; 
\beta_{j}^{(12)} =   \beta_{j} \cdot 1_{a_0 < j \leq s_0}, \;\text{ and } \; 
\beta^{(2)} =    \beta_{j} \cdot 1_{j > s_0}.
\eens
By definition of $s_0$ as in~\eqref{eq::define-s0}, we have 
$\sum_{i=1}^p \min(\beta_i^2, \lambda^2 \sigma^2)  \leq  
s_0 \lambda^2 \sigma^2$. 
Now
\ben
\nonumber
\sum_{j \leq a_0} \min(\beta_j^2, \lambda^2 \sigma^2) & = &    
a_0 \lambda^2 \sigma^2, \text{ since } \abs{\beta_j} \ge \lambda \sigma, \; \text{ and hence} \\
\label{eq::SR-range}
\sum_{j > a_0} \min(\beta_j^2, \lambda^2 \sigma^2) & = & 
\twonorm{ \beta^{(12)} + \beta^{(2)}}^2 \leq  (s_0 - a_0) \lambda^2 \sigma^2.
\een
For $\drop_{11} = \drop \cap A_0$, we have 
$\drop_{11} \subset A_0 \subset T_0 \subset S$, and $|\drop_{11}| \leq
a_0$.
Then $\norm{\beta_{\drop_{11}, \init}}_{\infty} < t_0$.
Let  $\beta_{\drop}^{(11)} := (\beta_j)_{j \in A_0 \cap \drop}$. 
Now by the triangle inequality,
\ben
\nonumber
\twonorm{\beta_{\drop}^{(11)}}
& = & \twonorm{\beta_{\drop_{11}} }
\leq 
\twonorm{\beta_{\drop_{11}, \init}} +
\twonorm{\beta_{\drop_{11}, \init} -   \beta_{\drop_{11}} } \\
\label{eq::drop-bound-3}
& \leq &
t_0 \sqrt{\size{\drop_{11}}} +  \twonorm{h_{\drop_{11}}} 
\leq
t_0 \sqrt{a_0} +  \twonorm{h_{\drop{11}}}; \\
\nonumber
\text{By~\eqref{eq::SR-range}},
\twonorm{\beta_{\drop}}^2 & \leq &  \twonorm{\beta_{\drop}^{(11)}}^2 + 
\twonorm{ \beta^{(12)} + \beta^{(2)}}^2  \\
\nonumber
& \leq & \twonorm{\beta_{\drop_{11}} }^2
+ (s_0 - a_0) \lambda^2
\sigma^2, \; \text{from which \eqref{eq::off-beta-norm-bound-2} follows}. 
\een
Now~\eqref{eq::off-beta-norm-bound-alt} holds
when we replace the crude bound of $|\drop_{11}| \leq a_0$ with
\bens
|\drop_{11}| & \leq & 
\frac{\twonorm{h_{\drop_{11}}}^2}{|\beta_{\min,A_0} - t_0|^2} \text{
  in~\eqref{eq::drop-bound-3} to obtain} \\
\twonorm{\beta_{\drop}^{(11)}} & \leq &  t_0 \frac{\twonorm{h_{\drop_{11}}}}{\beta_{\min,A_0} - t_0}
+  \twonorm{h_{\drop_{11}}} = \twonorm{h_{\drop_{11}}}
\frac{\beta_{\min,A_0}}{\beta_{\min,A_0} - t_0}.
\eens
\end{proofof2}

\subsection{Proof of Lemma~\ref{lemma:threshold-general-II}}
\label{sec::proofofA0}
\begin{proofof2}
  By the choice of $t_0$ in~\eqref{eq::ideal-t0} and by 
\eqref{eq::betaA-min-cond},
\bens
\min_{i \in A_0} \abs{\beta_{\init, i}} \geq
\abs{\beta_{\min, A_0}  -  \norm{h_{A_0}}_{\infty}} \geq t_0,\quad \text {and} \quad
\drop_{11} = \emptyset.
\eens
Thus by~\eqref{eq::ideal-t0}, we can bound
$|I \cap T_0^c|$, by either
\bens
|I \cap T_0^c| \leq
 {\norm{\beta_{T_0^c, \init}}_1}/{t_0} \leq \breve{s}_0 \quad \text{
   or } \quad |I \cap T_0^c|  \leq 
 {\twonorm{\beta_{T_0^c, \init}}^2}/{t_0^2} \leq \breve{s}_0,
 \eens
 depending on which one is applicable. The rest follows from the proof
 of Lemma~\ref{lemma:threshold-general}.
\end{proofof2}

\section{ Proof sketch of Theorem 3 in~\cite{ZH08}}
\label{sec::ZH08proof}
Let $\hat{S} :=\supp(\beta_{\init})$.
Following~\cite{ZH08}, we use $A_1$
to denote the union of support set for $\beta_{\init}$ and $H_q$:
\bens
A_1 := \{j: \beta_{\init, j} \not=0 \; \text{or} \; j \in H_q\} =
\hat{S}(\lambda) \cup H_q.
\eens
Then clearly, $A_1^c \subset H_0 =\{q+1, \ldots, p\}$.
In fact, by assumption~\eqref{eq::qstar}, $\size{H_q \cup T_1^*} <
q^*/(2 \sqrt{c^*})$ as shown earlier.
From this, we know that $q^* =\Omega(s_0)$; cf. \eqref{eq::qstars0}.

Under the SRC~\eqref{eq::SRC} and sparsity
conditions~\eqref{eq::qstar},~\eqref{eq::eta1def},
and~\eqref{eq::eta2def}, and with suitable choices of the penalty
$\lambda_n = 2 \lambda \sigma \sqrt{(1+a)  c^*}$ for some $a >0$,
the following statements  (along with other results) hold with high
probability (w.h.p.),
\ben
\label{eq::supploss}
(A_1) \; \; \size{\supp(\beta_{\init}) \cup H_q} & = & \size{A_1} \le
M_1^* q < q^{*}; \text{cf. (2.21)  and  (3.1)}, \\
\label{eq::headL2loss}
\twonorm{\beta_{\init, A_1} - \beta_{A_1}} & \le &
O(1) \lambda \sigma \sqrt{\abs{A_1}} \le
O(1) \lambda \sigma \sqrt{q^{*}}, \\
\label{eq::signalloss}
\twonorm{X (\beta_{\init} - \beta)}/\sqrt{n}
& \le & O(1)  \sigma \lambda \sqrt{\abs{A_1}} \quad \text{cf. (3.5)}.
\een
By definition of
$A_1 := \supp(\beta_{\init}) \cup H_q$, it holds that $\beta_{\init, j} =0, \; \forall 
j \in A_1^c$. Moreover, since $A_1^c \subset H_0$,  we have by~\eqref{eq::eta1def},
\ben
\nonumber
\norm{\beta_{\init, A_1^c} - \beta_{A_1^c}}_1 & \le &
\norm{\beta_{H_0}}_1 \le \eta_1 =O(\frac{r_1^2 q}{\sqrt{c^{*}}} 2\lambda
\sigma) =\tilde{O}(q^* \lambda \sigma), \\
\label{eq::tailL2loss}
\twonorm{\beta_{\init, A_1^c} - \beta_{A_1^c}}^2 & \le &
\twonorm{\beta_{H_0}}^2  = 
\tilde{O}(\frac{r_2}{\sqrt{c_*}}q^* \lambda^2 \sigma^2).
\een
To show~\eqref{eq::tailL2loss}, we first bound the $\ell_2$ norm on
the following set $T_1^* =\{j: j \in H_0, \abs{\beta_j} \ge \lambda
\sigma\}$ by~\eqref{eq::eta2def},
\bens
\twonorm{\beta_{T_1^*}}^2 
& = & 
\sum_{j >q}^{p} \beta_j^2 I(\abs{\beta_j} \ge \lambda \sigma)
\le \twonorm{\sum_{j \in T_1^*} \beta_j x_j}^2/(n c_*) \\
& \le & \eta_2^2/ (n c_*) =\tilde{O}(q (2 r_2 \sigma
\lambda)^2/{c_*});
\eens
Now clearly, \eqref{eq::tailL2loss} holds since for  all $j \in H_0 \setminus
T_1^*$, $\abs{\beta_j} <\lambda \sigma$, and hence
\ben
\label{eq::thintail}
\twonorm{\beta_{H_0 \setminus T_1^*}}^2
& = & \twonorm{\beta_{[p] \setminus L(q)}}^2
\le \norm{\beta_{H_0 \setminus T_1^*}}_1 (\lambda \sigma)\\
& \le &
\eta_1 \lambda \sigma
= \frac{r_1^2}{\sqrt{c^*}} q 2 \sigma^2 \lambda^2 =
\tilde{O}(2r_1^2 q (\sigma \lambda)^2).
\een
Putting these two sets together, one can obtain
\eqref{eq::tailL2loss}, since
\ben
\nonumber
\twonorm{\beta_{H_0}}^2
& = & 
\twonorm{\beta_{T_1^*}}^2  + \twonorm{\beta_{H_0 \setminus T_1^*}}^2
\le \tilde{O}((\frac{r_2^2}{c_*} + r_1^2) q (2 \sigma \lambda)^2)\\
& = &
\tilde{O}(\frac{r_2}{\sqrt{c_*}} q^* \sigma^2 \lambda^2) \text{ since
 }\; \; 4 (\frac{r_2}{\sqrt{c_*}} + r_1^2) q \le q^*.
 \een
 Then one obtains \eqref{eq::L1lossmain}
by~\eqref{eq::eta1def}, \eqref{eq::headL2loss} and \eqref{eq::tailL2loss}. 
However, the largest signal in $H_0$ can be $\asymp \sqrt{q} \sigma \lambda$ via the $\eta_2$ condition~\eqref{eq::eta2def}:
\ben
\label{eq::largesignal}
\norm{\beta_{H_0}}_{\infty} \le \sqrt{\eta_2^2/ (n c_*)} =\tilde{O}(
\frac{2 r_2}{\sqrt{c_*}} \sqrt{q} \sigma \lambda).
\een
We compare the conditions and error bounds with those in
Theorem~\ref{thm::RE-oracle} in~Section~\ref{sec::tightness}.

\section{Conclusion}
\label{sec:conclude}
In this paper, we show that the thresholding method is effective in 
variable selection and accurate in statistical estimation. 
It improves the ordinary Lasso in significant ways.
For example, we allow very significant  number of non-zero elements in the 
true parameter,  for which the ordinary Lasso would have failed.
On the theoretical side, we show that if $X$ obeys the RE
condition and if the true parameter is sufficiently sparse, the Thresholded 
Lasso achieves the $\ell_2$ loss within a logarithmic factor of the 
{\it ideal mean square error} one would achieve with an oracle, 
while selecting a sufficiently sparse model $I$.
This is accomplished when threshold level is at about 
$\sqrt{2\log p/n} \sigma$, assuming that columns of $X$ have 
$\ell_2$ norm $\sqrt{n}$.
When the SNR is high, almost exact recovery of the 
non-zeros in $\beta$ is possible as shown in our theory; exact recovery 
of the support of $\beta$ is shown in our simulation study when $n$ 
is only linear in $s$ for several Gaussian and Bernoulli random ensembles.
When the SNR is relatively low, the inference task is difficult for any
estimator. In this case, we show that Thresholded Lasso tradeoffs 
Type I and II errors nicely: we recommend choosing the thresholding 
parameter conservatively. 
These findings not only validate our theoretical analysis excellently
but also indicate that in practical applications, this approach could 
be made very effective and relevant.

\section*{Acknowledgement}

I dedicate this work to Professor Peter B\"{u}hlmann, on the occasion
of his 60th birthday. His unwavering support and mentorship during my career has been
deeply appreciated. The author thanks the Editor,
an AE and two anonymous referees for their valuable comments and suggestions. I thank
Professor Bin Yu warmly for hosting my visit at UC Berkeley while I finished the first
draft of this work in Spring 2010. The author thanks Peter B\"{u}hlmann, Emmanuel Cand\`{e}s,
Sara van de Geer, John Lafferty, Po-Ling Loh, Richard Samworth, Xiaotong Shen, Martin
Wainwright, Larry Wasserman, Dana Yang and Cun-Hui Zhang for helpful discussions. Part
of this work has appeared in a conference paper with title: Thresholding Procedures for High
Dimensional Variable Selection and Statistical Estimation, in Proceedings of Advances in
Neural Information Processing Systems 22.

\bibliography{subgaussian}

\appendix

\section*{Notation}
Let $T \subseteq [p]$ be a fixed subset of indices.
Let $X_T$ be the $n \times  \abs{T}$ submatrix obtained by extracting columns of $X$ indexed by $T$.
We use $\upsilon_{T}$ to represent the subvector of
$\upsilon \in \R^p$ confined to a subset $T$ of $[p]$.
Let $\shtwonorm{\upsilon}^2 = \sum_{j=1}^p \upsilon_j^2$.
Occasionally, we use $\upsilon_T \in \R^{\abs{T}}$ to also represent
its $0$-extended version $\upsilon' \in \R^p$ such that
$\upsilon'_{T^c} = 0$ and $\upsilon'_{T} = \upsilon_T$,
for example in Lemma~\ref{lemma::T01simple} and its proof
when it is clear from the context.
Let $\beta_T$ be the restriction of $\beta \in \R^p$ to the set $T$.
We will also explicitly use $\beta^{\ext}(T) \in \R^p$ to represent this $0$-extended
version of $\beta_T \in \R^{\abs{T}}$ such that $\beta_{T}^{\ext}({T}) = \beta_{T}$
and $\beta_{T^c}^{\ext}({T}) = 0$.
Given $\hat{\beta}_I = (X_I^T X_{I})^{-1} X_{I}^T Y \in \R^{\abs{I}}$, we use
$\hat{\beta}_I$ to represent  its $0$-extended version $\hat\beta^{\OLS} \in \R^p$ such that 
$\hat\beta^{\OLS}_I = \hat{\beta}_I$ and $\hat\beta^{\OLS}_{I^c} = \hat{\beta}_{I^c} =0$.

\section{Proof of the MSE lower bound}
\label{section:append-diamond}
We show the proof of~\eqref{eq::ME-diamond} for self-containment.
Note that due to different normalization of columns of $X$, 
our expressions are slightly different from those by~\cite{CT07}. 
Hence we give a complete derivation here.
  Let $I$ be a fixed subset of indices and consider the OLS 
  estimator $\hat\beta^{\OLS}$ such that 
$\hat\beta^{\OLS}_I = \hat{\beta}_I = (X_I^T X_{I})^{-1} X_{I}^T Y$
and $\hat\beta^{\OLS}_{I^c} = \hat{\beta}_{I^c} =0$. 
Consider the least square estimator 
$\hat{\beta}_I = (X_I^T X_{I})^{-1} X_{I}^T Y$, 
where $|I| \leq s$ (and $\hat{\beta}_{I^c}=0$), and consider the ideal least-squares estimator 
$\beta^{\diamond}$ which minimizes the expected mean squared error
\begin{eqnarray}
\label{eq::beta-diamond-proof}
\beta^{\diamond} = \argmin_{I \subset \{1, \ldots, p\}, \; |I| \leq s} 
\E \twonorm{\beta- \hat{\beta}^{\OLS}(I)}^2.
\end{eqnarray}
\begin{proposition}{\textnormal{\citep{CT07}}} 
If $\Lambda_{\max}(s) < \infty$, then
\begin{eqnarray}
\label{eq::ME-diamond-append}
\E \twonorm{\beta- \beta^{\diamond}}^2 \geq
\min\left(1, 1/\Lambda_{\max}(s) \right) 
\sum_{i=1}^p \min(\beta_i^2, \sigma^2/n).
\end{eqnarray}
\end{proposition}
\begin{proof}
  Here we denote by $\beta_I$ the restriction of $\beta$ to the set $I$, 
where $I \subset [p]$, and $\beta^{\ext}({I})$ its $0$-extended 
version such that $\beta^{\ext}_{I^c}(I)= 0$ and $\beta^{\ext}_{I}(I) = \beta_I$. 
The error of this OLS  estimator, namely, $\hat\beta^{\OLS}$ is given by
\ben
\label{eq::error-decomp}
\twonorm{\hat{\beta}^{\OLS}(I) - \beta}^2 & = &
\twonorm{\hat{\beta}_I - \beta_{I}}^2 + 
\twonorm{\beta^{\ext}({I}) - \beta}^2.
\een
The first term is:
\bens
\label{eq::MSE-I-est-missing-data}
\hat\beta_{I} - \beta_I 
& := & (X_{I}^T X_{I})^{-1}X_{I}^T Y - \beta_I \\
& = & (X_{I}^T X_{I})^{-1} X_{I}^T (X_I \beta_I + X_{I^c} \beta_{I^c} + \epsilon)  - \beta_I \\
& = & 
(X_{I}^T X_{I})^{-1} X_{I}^T X_{I^c} \beta_{I^c} +
(X_{I}^T X_{I})^{-1}X_{I}^T \epsilon,
\eens
and its mean squared error is given by
\bens
\E \twonorm{\hat\beta_I - \beta_I}^2 & = & 
\twonorm{(X_{I}^T X_{I})^{-1} X_{I}^T X_{I^c} \beta_{I^c}}^2 +
\E \twonorm{(X_{I}^T X_{I})^{-1}X_{I}^T \epsilon}^2,
\eens
where
\bens
\E \twonorm{(X_{I}^T X_{I})^{-1}X_{I}^T \epsilon}^2
& = &
\frac{\sigma^2}{n} 
\trace\left(\left(\frac{X_{I}^T X_{I}}{n}\right)^{-1} \right)
\geq  
\frac{\sigma^2}{n} \frac{|I|}{\Lambda_{\max}(|I|)},
\eens
since eigenvalues of $\left(\frac{X_{I}^T X_{I}}{n}\right)^{-1}$
are in the range of 
$\left[\inv{\Lambda_{\max}\left({X_{I}^T X_{I}}/{n}\right)}, 
\inv{\Lambda_{\min}\left({X_{I}^T X_{I}}/{n}\right)} \right]$.
Thus
\ben
\E \twonorm{\hat\beta_I - \beta_I}^2 & \geq & 
\frac{\sigma^2}{n} \frac{|I|}{\Lambda_{\max}(|I|)}.
\een
Thus for all sets $I$ such that $|I| \leq s$ and for 
$\Lambda_{\max}(s) < \infty$, we have for $\hat{\beta}_{I^c} =0$ and by~\eqref{eq::error-decomp},
\bens
\E \twonorm{\hat\beta^{\OLS}({I})  - \beta}^2
& = & 
\E \twonorm{\hat\beta_{I} - \beta_I}^2 + \twonorm{\beta_{I^c}}^2 \\
& \geq &  
\frac{\sigma^2}{n} \frac{|I|}{\Lambda_{\max}(s)} + 
\twonorm{\beta_{I^c}}^2  \\
& \geq &  
\min\left(1, 1/\Lambda_{\max}(s) \right) 
\left(\sum_{j \in I^c} \beta_{j}^2 + \frac{\sigma^2}{n} |I| \right),
\eens
which gives that the ideal mean squared error is bounded below by
\bens
\E \twonorm{\beta- \beta^{\diamond}}^2 
& \geq & 
\min\left(1, 1/\Lambda_{\max}(s) \right) 
\min_{I} \left(\sum_{j \in I^c} \beta_{j}^2 + 
\frac{\sigma^2}{n} |I| \right) \\
& = & 
\min\left(1, 1/\Lambda_{\max}(s) \right) 
\sum_{i=1}^p \min(\beta_i^2, \sigma^2/n).
\eens
\qed
\end{proof}

\section{Proof of Lemma~\ref{lemma:parallel}}
\label{sec::proofparallel}
\begin{proofof2}
It is sufficient to show that~\eqref{eq::parallel} holds for
$\twonorm{c}  = \twonorm{c'} = 1$.
\ben
\label{eq::parallel}
\frac{\abs{{\ip{X_I c, X_{S_{\drop}} {c'}}}}}{n} \leq  
\frac{(\Lambda_{\max}(2s) - \Lambda_{\min}(2s))  }{2}.
\een
Indeed, by~\eqref{eq::eigen-admissible-s} and~\eqref{eq::eigen-max},
we have 
$2 \Lambda_{\min}(2s) \leq
\twonorm{X_I c +  X_{S_{\drop}} {c'}}^2/{n} \leq 2 \Lambda_{\max}(2s)$
and
$2 \Lambda_{\min}(2s) \leq 
{\twonorm{X_I c -  X_{S_{\drop}} {c'}}^2 }/{n}  \leq 2 \Lambda_{\max}(2s) $.
Hence~\eqref{eq::parallel} follows from the parallelogram identity:
$$\abs{{\ip{X_I c, X_{S_{\drop}} {c'}}}} =
\abs{\twonorm{X_I c +  X_{S_{\drop}} {c'}}^2  - \twonorm{X_I c -
    X_{S_{\drop}} {c'}}^2 }/{4}.$$
Next, suppose $\Lambda_{\min}(2s) =0$.
Then
\bens
\abs{{\ip{X_I c, X_{S_{\drop}} {c'}}}}
& = & 
\abs{\twonorm{X_I c +  X_{S_{\drop}} {c'}}^2  - \twonorm{X_I c -
    X_{S_{\drop}} {c'}}^2 }/{4} \\
& \le & 
\frac{\twonorm{X_I c +  X_{S_{\drop}} {c'}}^2}{4}  \vee\frac{
  \twonorm{X_I c -    X_{S_{\drop}} {c'}}^2 }{4} \\
& \le & \Lambda_{\max}(2s)/2.
\eens
Moreover,~\eqref{eq::sparseupper} follows from the arguments 
as in~\eqref{eq::orthocauchy}, using the Cauchy-Schwarz inequality. 
\end{proofof2}

\section{Proof of  Lemma~\ref{lemma::h01-bound-CT}}
\label{sec::proofofT01}
\begin{proofof2}
Decompose $h_{T_{01}^c}$ into $h_{T_2}$, \ldots, $h_{T_K}$ 
such that $T_2$ corresponds to locations of the $s_0$ largest coefficients of 
$h_{T_{01}^c}$ in absolute values, and $T_3$ corresponds to locations of the 
next $s_0$ largest coefficients of $h_{T_{01}^c}$ in absolute values, 
and so on. 
Let $V$ be the span of columns of $X_j$, where $j \in T_{01}$, and 
$P_V$ be the orthogonal projection onto $V$.
Decompose $P_V X h$:
\bens
P_V X h & = & P_V X h_{T_{01}} + \sum_{j \geq 2} P_V X h_{T_j}
= X h_{T_{01}} + \sum_{j \geq 2} P_V X h_{T_j} , \text{ where } \\
\forall j \ge 2, \; \; \twonorm{P_V X h_{T_j}} 
& \leq & 
\frac{\sqrt{n} \theta_{s_0, 2s_0}}{\sqrt{\Lambda_{\min}(2s_0)}}\twonorm{h_{T_j}}
\text{ and }
\sum_{j \geq 2} \twonorm{h_{T_j}} \leq  \norm{h_{T_0^c}}_1 /\sqrt{s_0};
\eens
see~\cite{CT07} for details. Thus we have
\bens
\twonorm{X h_{T_{01}}} & = & 
\twonorm{P_V X h - \sum_{j \geq 2} P_V X h_{T_j} }  \leq 
\twonorm{P_V X h} + \twonorm{\sum_{j \geq 2} P_V X h_{T_j} } \\
& \leq & 
\twonorm{X h} + \sum_{j \geq 2} \twonorm{P_V X h_{T_j} } 
\leq
\twonorm{X h} +  \frac{\sqrt{n} \theta_{s_0, 2s_0}}
{\sqrt{\Lambda_{\min}(2s_0)} \sqrt{s_0}} \norm{h_{T_0^c}}_1,
\eens
where  we used the fact that $\twonorm{P_V} \leq 1$.
Hence the lemma follows given 
$$\twonorm{h_{T_{01}}} \leq 
\inv{\sqrt{\Lambda_{\min}(2s_0)} \sqrt{n}}  \twonorm{X h_{T_{01}}}.$$
For other bounds,
the fact that the $k$th largest value of $h_{T_0^c}$ obeys
$\size{h_{T_{0}^c}}_{(k)} \leq  \norm{h_{T_0^c}}_1 / k$ has been used; 
see~\cite{CT07}.
\end{proofof2}

\section{Proof of Theorem~\ref{thm::RE-oracle}}
\label{sec::LassoProof}
Let $\beta_{T_0}$ be the restriction of $\beta$  to the set $T_0$.
Let $h = \beta_{\init} - \beta^{\ext}({T_0})$. 
Throughout this proof,  as a shorthand, we write $h = \hat{\beta} - \beta_0$, where $\hat\beta  =
\beta_{\init}$ and $\beta_0 = \beta^{\ext}({T_0})$ as defined in Lemma~\ref{lemma:threshold-RE}. 
We first show 
Lemma~\ref{lemma:T0-pre-loss},
which gives us the prediction error
using $\beta_{T_0}$.
We do not focus on obtaining the best constants in the proof of
Theorem~\ref{thm::RE-oracle}.
Recall  we define a quantity $\basepen$, 
which bounds the maximum correlation between the noise and 
covariates of $X$; For each $a \geq 0$, let
\ben
\label{eq::low-noise-supp}
{\T_a} := 
\biggl \{\e: \norm{{X^T \e}/{n}}_{\infty} \leq \basepen, 
\text{ where }
\basepen = \sigma \sqrt{1 + a} \sqrt{{2\log p}/{n}}\biggr \}.
\een
Then, we have $\prob{\T_a}  \geq 1 - (\sqrt{\pi \log p} p^a)^{-1}$
when $X$ has column $\ell_2$ norms bounded by $\sqrt{n}$.
\begin{lemma}
  \label{lemma:T0-pre-loss}
  Suppose $\beta$ is {\bf $s$-sparse}.
 Let $T_0$  denote  locations of the $s_0$ largest coefficients of  $\beta$ in absolute
 values.
 Suppose~\eqref{eq::eigen-max} holds. We have
for $\lambda = \sqrt{(2 \log p)/n}$.
\ben
\label{eq::T0-pre-loss}
\twonorm{X \beta - X \beta_{0}}/\sqrt{n}  & \leq & 
\sqrt{\Lambda_{\max}(s- s_0)} \lambda \sigma \sqrt{s_0},
\een
where recall $\beta_0 = \beta^{\ext}(T_0)$, as defined in
Lemma~\ref{lemma:T0-pre-loss}.
\end{lemma}
\begin{proof}
The lemma holds given that 
$\twonorm{\beta_{T_0^c}} \leq \lambda \sigma \sqrt{s_0}$
by~\eqref{eq::define-s0} and \eqref{eq::beta-2-small-intro}.
Indeed, we have $\abs{\beta_j} < \lambda \sigma$ for all $j \in  
T_0^c$ by definition of $T_0$ and hence
\bens
\twonorm{\beta_{T_0^c}}^2 & \le & \sum_{i=1}^p \min(\beta_i^2, \lambda^2 \sigma^2) \leq 
s_0 \lambda^2 \sigma^2.
\eens
Hence
$\twonorm{X \beta - X \beta_{0}}/\sqrt{n} 
\leq \sqrt{\Lambda_{\max}(s- s_0)} \twonorm{\beta_{T_0^c}}.$
\qed
\end{proof}

Next, we state Lemma~\ref{lemma::T01simple}. 
We prove Lemma~\ref{lemma::T01simple} in Section~\ref{sec::suppT01simpleproof}.
\begin{lemma}
 \label{lemma::T01simple}
 Suppose all conditions in Theorem~\ref{thm::RE-oracle} hold.
 Let $T_1$ be the $s_0$ largest positions of $h$ outside of $T_0$. 
 Denote by $T_{01} = T_0 \cup T_1$. 
 Then
  \ben
  \label{eq::simpleT01}
\twonorm{X h_{T_{01}}}/\sqrt{n}
& \leq & 
\twonorm{X h} /\sqrt{n}+  {\sqrt{\Lambda_{\max}(s_0)}}
\norm{h_{T_0^c}}_1/{\sqrt{s_0}}.
\een
Then by Lemma~\ref{lemma::h01-bound-CT} and~\eqref{eq::simpleT01},  we have
\bens
\nonumber 
\twonorm{h_{T_{01}}} 
& \leq &
\inv{\sqrt{\Lambda_{\min}(2s_0)}}
\big(\norm{X h}_{n} +  \ell(s_0) \norm{h_{T_0^c}}_1 /\sqrt{s_0} \big)
\\
\label{eq::definell}
\; \text{where} \; \ell(s_0) 
& := &  \frac{\theta_{s_0, 2s_0} }{\sqrt{\Lambda_{\min}(2s_0)}}
\wedge \sqrt{\Lambda_{\max}(s_0)}.
\eens
\end{lemma}

\subsection{Proof of Theorem~\ref{thm::RE-oracle}}
Throughout this proof, we assume that $\T_a$ holds.

\noindent{\bf Prelude.}
First recall the following elementary inequality.
By the optimality of $\hat{\beta} = \beta_{\init}$,
we have
\begin{eqnarray}
\label{eq::opt-lasso} 
\; \; \; \;\; 
\inv{2n} \norm{Y- X  \hat\beta}^2_2  - \inv{2n} \norm{Y- X \beta_{0}}^2_2 
 \leq 
  \lambda_{n} \norm{\beta_{T_0}}_1 - \lambda_{n} \norm{\hat\beta}_1,
  \end{eqnarray} 
  \text{ where} \\
\bens
\twonorm{Y-X \hat \beta}^2 & = & 
\twonorm{X \beta -X \hat \beta + \epsilon}^2\\
& = & 
\twonorm{ X \hat \beta - X \beta }^2
+ 2(\beta - \hat\beta )^T X^T \e + \norm{\e}^2_2, \text{ and } \\
\label{eq::beta0-explain} 
\twonorm{Y-X  \beta_0}^2
& = & 
\twonorm{X \beta -X \beta_0 + \epsilon}^2 \\
  & = & 
\twonorm{ X \beta - X \beta_0 }^2
+ 2(\beta - \beta_0)^T X^T \e + \norm{\e}^2_2.
\eens
Thus by~\eqref{eq::opt-lasso} and the triangle inequality,
we have on $\T_a$ and $\lambda_n \ge 2 \lambda_{\sigma, a, p}$ and
$\norm{\frac{X^T \e}{n}}_{\infty} \le  \lambda_{\sigma, a, p} \le \lambda_n/2$,
\ben
\label{eq::missingone}
\frac{\twonorm{ X \hat \beta - X \beta }^2}{n}
& \leq  &
\nonumber
\frac{\twonorm{ X \beta - X \beta_0 }^2}{n}
+ \frac{2 h^T  X^T \e}{n} + 2 \lambda_{n}
( \norm{\beta_{0}}_1 -  \norm{h + \beta_0}_1) \\
\nonumber
& \leq  & 
\frac{\twonorm{ X \beta - X \beta_0 }^2}{n}
+ 2 \norm{h}_1 \norm{\frac{X^T \e}{n}}_{\infty} + 
2 \lambda_{n} (\norm{h_{T_0}}_1 -  \norm{h_{T_0^c}}_1 )\\
\nonumber
& \leq  & 
\frac{\twonorm{ X \beta - X \beta_0 }^2}{n}
+ 3 \lambda_n \norm{h_{T_0}}_1 -  \lambda_n \norm{h_{T_0^c}}_1,
\een
where we have used the fact that 
$\lambda_n \geq 2 \basepen$ for $a \geq 0$.
Thus we have on $\T_a$,  
\ben
\label{eq::long-shot}
& & 
\twonorm{ X \hat \beta - X \beta }^2/n
+ \lambda_n \norm{h_{T_0^c}}_1
\leq \twonorm{ X \beta - X \beta_0 }^2/n + 3 \lambda_n  \norm{h_{T_0}}_1,
\een
This is the inequality we used in~\cite{Zhou10}.
In the updated proof, we will apply the following
Lemma~\ref{lemma::convex}; cf. Lemma A.2~\cite{BLT18}. 
\begin{lemma}{\textnormal{[Lemma A.2]}~\citep{BLT18}}
  \label{lemma::convex}
  Let $h: \R^p \longrightarrow \R$ be a convex function.
  Let $f, \xi \in \R^n$, $y = f+ \xi$ and let $X$ be any $n \times p$
  matrix.
  If $\hat\beta$ is a solution of the minimization problem
  $\min_{\beta \in \R^p} (\norm{X \beta - y}^2_n + h(\beta)$, then
  $\hat\beta$ satisfies for all $\tilde\beta\in \R^p$
  \ben
  \label{eq::A4}
\lefteqn{    \norm{X \hat\beta - f}^2_n + \norm{X (\hat\beta
    -\tilde\beta)}^2_n}\\
\nonumber
&&    \le  \norm{f- X \tilde\beta}^2_n +  \frac{2}{n} \xi^T  X (\hat{\beta} - \tilde\beta) + h(\tilde\beta) -
    h(\hat\beta) 
    \een
  \end{lemma}
We use $\hat{\beta} := \beta_{\init}$ to represent the solution
to the Lasso estimator in~\eqref{eq::origin}.
Using Lemma~\ref{lemma::convex} with $y = f + \e$, where $f = X
\beta$, we obtain eq. (20)~\citep{DHL17},
where we set $\bar\beta = \beta_0 :=\beta^{\ext}({T_0})$ and $\bar\delta =\hat\beta -\beta_0  =h$,
\ben
\nonumber
\lefteqn{ \norm{X (\hat\beta - \beta)}^2_n  + \norm{X (\hat\beta -\beta_0 )}^2_n
  \le   \norm{X (\beta -\beta_0)}^2_n} \\
\label{eq::oracleD}
&& + \frac{2}{n} \e^T  X(\hat{\beta} - \beta_0) + 2 \lambda_n
(\norm{\beta_0}_1 -\norm{\hat\beta}_1)
\een
Thus, we have the following updated inequality~\eqref{eq::new-oracle},
replacing  \eqref{eq::missingone}.

Denote by $\delta :=\hat \beta -  \beta$ and $h = \hat \beta -
\beta_0$. Thus we have by \eqref{eq::oracleD}, on  $\T_a$,
\ben
\norm{X \delta }^2_{n} + \norm{X h}^2_{n}
& \leq  &
\nonumber
\norm{ X \beta - X \beta_0 }^2_{n}
+ \frac{2 h^T  X^T \e}{n} + 2 \lambda_{n} (\norm{\beta_{0}}_1 -  \norm{h + \beta_0}_1) \\
\nonumber
& \leq  & 
\norm{ X \beta - X \beta_0 }^2_{n}
+ 2 \norm{h}_1 \norm{ \frac{X^T \e}{n}}_{\infty} + 
2 \lambda_{n} (\norm{h_{T_0}}_1 -  \norm{h_{T_0^c}}_1 )\\
\label{eq::oracle0}
& \leq  &
\norm{ X \beta - X \beta_0 }^2_{n}
+ 3 \lambda_n \norm{h_{T_0}}_1 -  \lambda_n \norm{h_{T_0^c}}_1,
\een
where we have used the fact that 
$\lambda_n \geq 2 \basepen$ for $a \geq 0$.
\eqref{eq::new-oracle} then follows  immediately  from \eqref{eq::oracle0}.
We consider three cases which do not need to be mutually exclusive.

\noindent{\bf Case 1.}
In the first case, suppose 
\bens 
\norm{X\delta}^2_n  + \norm{X h}^2_{n} \ge  \norm{X \beta -
  X \beta_0}_n^2.
\eens
Then by \eqref{eq::new-oracle} and \eqref{eq::DHDomi},
\ben
\label{eq::hcone3}
\norm{h_{T_0^c}}_1  \le 3 \norm{h_{T_0}}_1, \text{ and hence }  \; \; h \in \Cone(s_0, 3). 
\een
Hence we can use the $\RE(s_0, 3, X)$ condition to bound
$\twonorm{h_{T_0}}$ with $K \norm{X h}_{n}$, 
\ben
\nonumber
\lefteqn{
\norm{X \delta}^2_n  + \norm{X h}^2_{n} -  \norm{X (\beta-\beta_0)}_n^2 
+ \lambda_n \norm{h_{T_0^c}}_1 }\\
\nonumber
& \leq &
3 \lambda_n \norm{h_{T_0}}_1  \le
3 \lambda_n \sqrt{s_0} \twonorm{h_{T_0}} \le
3 \lambda_n \sqrt{s_0} K \norm{X h}_{n}\\
  \label{eq::newcase1}
& \leq & (3 K \lambda_n \sqrt{s_0}/2)^2 + \norm{X h}^2_{n}\;
\text{where }\; K :=K(s_0, 3);
\een
Hence \eqref{eq::newboundsim}  follows immediately
from~\eqref{eq::newcase1}, upon deleting  the term $\norm{X h}^2_{n}$
from both sides.
Moreover, we obtain  the upper bound on $\norm{h_{T_0^c}}_1$ (
$\norm{h_{T_0}}$) from~\eqref{eq::newboundsim} since
\ben
\nonumber
  \norm{X \delta}^2_n
 + \lambda_n \norm{h_{T_0^c}}_1 
& \leq &\norm{X (\beta-\beta_0)}_n^2  + (3 K \lambda_n
\sqrt{s_0}/2)^2,\\
\nonumber
\norm{h_{T_0^c}}_1  
& \leq &\norm{X (\beta-\beta_0)}_n^2/\lambda_n + (3 K /2)^2\lambda_n 
s_0\\
& \leq &
\label{eq::HT0case1b}
\big(\Lambda_{\max}(s- s_0)/d_0+ 9 d_0 K(s_0, 3)^2
/{4}\big)  \lambda \sigma s_0,
\een
following the proof of Lemma~\ref{lemma:T0-pre-loss}, where
$\lambda_n = d_0 \lambda \sigma \geq 2 \basepen$.
Moreover, we will show that by~\eqref{eq::newcase1} and~\eqref{eq::newboundsim}, 
\ben
  \label{eq::case1twonorm}
\twonorm{h_{T_0}}
& \le & K \norm{X \beta - X \beta_0}_{n}  + 3K^2
\lambda_n \sqrt{s_0} \le D'_0\lambda \sigma \sqrt{s_0},\\
\nonumber
&& \text{ where } \; D'_0 := (K \sqrt{\Lambda_{\max}(s- s_0)} + 3d_0 K^2).
\een
Then, we have by \eqref{eq::hcone3} and   \eqref{eq::case1twonorm},
\ben
\label{eq::case1TOC}
\norm{h_{T_0^c}}_1  & \le &
3 \norm{h_{T_0}}_1 \le
3\sqrt{s_0} \twonorm{h_{T_0}} \le
3 D'_0\lambda \sigma s_0.
\een
Combining \eqref{eq::case1TOC} and~\eqref{eq::case1twonorm},
we have
\bens 
\norm{h_{T_0^c}}_1 & \le &
D_{1,a} \lambda \sigma s_0, \; \text{ where} \\
\; \;
D_{1,a} & :=  &
\big(\Lambda_{\max}(s- s_0)/d_0+ 9 d_0 K(s_0, 3)^2 
/{4}\big) \wedge 3 D'_0 \; \; 
\; \text{ where} \\
3 D'_0 & \le &3 K(s_0, 3) \sqrt{\Lambda_{\max}(s- s_0)} + 9d_0 K^2(s_0, 3)) \\
& \le &
\Lambda_{\max}(s- s_0)/d_0 + 9 d_0 K(s_0, 3)^2 /{4} + 
9d_0 K^2(s_0, 3).
\eens
Thus \eqref{eq::pred-error-gen} holds for {\bf Case 1}, since
for $\delta = X \hat \beta - X \beta$, we have by~\eqref{eq::newboundsim},
\bens
\norm{X \delta}_n & = & \norm{X \hat \beta - X \beta}_{n}  \le 
\big(\norm{ X \beta - X \beta_0 }^2_{n} + (\frac{3 K \lambda_n 
  \sqrt{s_0}}{2})^2\big)^{1/2} \\
& \le & 
\norm{ X \beta - X \beta_0 }_{n} + \frac{3 K(s_0, 3) \lambda_n 
  \sqrt{s_0}}{2}.
\eens

\noindent{\bf   Upper bound on $\twonorm{h_{T_0}}$~\eqref{eq::HT0case1a}.}
Under the $\RE(s_0, 3, X)$ condition, we have by
Definition~\eqref{eq::admissible} and~\eqref{eq::newbound3},
for $h = \hat{\beta} - \beta_0$, where $\beta_0 = \beta^{\ext}({T_0})$,
\bens
\twonorm{h_{T_0}}^2
& \leq & K(s_0, 3)^2 \norm{X h}^2_{n} \\
& \leq & 
K(s_0, 3)^2 
\left(\norm{X \beta - X \beta_0}^2_{n} + 3 \lambda_n  \sqrt{s_0} \twonorm{h_{T_0}}
\right),
\eens
where recall  we have on $\T_a$,  by \eqref{eq::new-oracle},
\ben 
\label{eq::newbound3}
\norm{X h}^2_{n} \leq 
\norm{X \beta - X \beta_0}_n^2 + 3 \lambda_n \norm{h_{T_0}}_1 
\; \text{ where} \;  \norm{h_{T_0}}_1 
\le \sqrt{s_0} \twonorm{h_{T_0}}.
\een 
Thus we have for $\lambda_n = d_0 \lambda \sigma$ and $K := K(s_0, 3)$,
\ben
\label{eq::headnorm}
&& \twonorm{h_{T_0}}^2- 3 K^2 \lambda_n \sqrt{s_0}
\twonorm{h_{T_0}} \le K^2   \norm{X \beta - X \beta_0}^2_{n}, \;
\text{and} \\
\nonumber
&& \big(\twonorm{h_{T_0}}- \frac{3K^2}{2} \lambda_n \sqrt{s_0}\big)^2
 \leq 
K^2 \norm{X \beta - X \beta_0}^2_{n}  + \big(\frac{3K^2}{2} \lambda_n 
\sqrt{s_0})^2.
\een
Thus we have \eqref{eq::HT0case1a} holds for Case 1 with $K := K(s_0,
3)$, since
\bens 
\abs{\twonorm{h_{T_0}}- \frac{3K^2}{2} \lambda_n \sqrt{s_0}}
& \le &
K \norm{X \beta - X \beta_0}_{n}  + \frac{3K^2}{2} \lambda_n 
\sqrt{s_0}, \text{ which implies that} \\
\twonorm{h_{T_0}} & \le & K \norm{X \beta - X \beta_0}_{n}  + 3K^2 \lambda_n \sqrt{s_0}.
\eens
Thus we obtain
\ben
\nonumber
\twonorm{h_{T_0}}
& \le & D'_{0,a} \lambda \sigma \sqrt{s_0}, \; \text{ where } \\
\label{eq::HT0case1}
D'_{0,a} & := & K(s_0, 3) \sqrt{\Lambda_{\max}(s- s_0)} + 3K(s_0, 3)^2 d_0.
\een
Similarly,  we can derive a bound on $\norm{h}_1$
following~\eqref{eq::newcase1} directly:
\bens
\nonumber
\lefteqn{
\norm{ X \delta}^2_{n} + \norm{ X h}^2_{n}
+ \lambda_n \norm{h_{T_0^c}}_1 + \lambda_n \norm{h_{T_0}}_1 
- \norm{ X \beta - X \beta_0 }^2_n }\\
& \leq & 4 \lambda_n  \norm{h_{T_0}}_1 \le
4  \lambda_n  \sqrt{s_0} \twonorm{h_{T_0}} \\
 & \leq &
4  K \lambda_n  \sqrt{s_0} \norm{X h}_{n} \le
\norm{X h}^2_n + (2 K \lambda_n  \sqrt{s_0})^2,
\eens
under $\RE(s_0, 3, X)$ condition, since $h \in \Cone(s_0, 3)$. 
Hence
\bens
\nonumber
\norm{ X \delta}_n^2  + \lambda_n \norm{h}_1 
& \leq & \norm{ X \beta - X \beta_0 }^2_{n}
+ (2 K \lambda_n  \sqrt{s_0})^2.
\eens
Then
\ben
\nonumber
\norm{h}_1 & = & \norm{h_{T_0^c}}_1 + \norm{h_{T_0}}_1 \le
\norm{ X \beta - X \beta_0 }_{n}^2/ \lambda_n + 4 K^2 \lambda_n  s_0 \\
\label{eq::H1normC1}
& \leq & D_{2,a}
\lambda \sigma s_0, \; \text{ where } \;
 D_{2,a} \leq \Lambda_{\max}(s- s_0)/d_0 + 4 K(s_0, 3)^2 d_0. 
\een

\noindent{\bf Case 2.}
Suppose
\ben
\label{eq::DHhead}
 \norm{X\delta}^2_n  + \norm{X h}^2_{n}
  \le \lambda_n \norm{h_{T_0}}_1. 
  \een
  Let $T_1$ be the $s_0$ largest positions of $h$ outside of $T_0$.
  This is the easy case. First we show \eqref{eq::oracle2}.
 By the triangle inequality,
  \bens
  \norm{X \beta - X \beta_0}_n \le  \norm{X (\beta -
    \hat{\beta})}_n +  \norm{X (\hat\beta - \beta_0)}_n =
  \norm{X \delta}_n  + \norm{X h}_n.
  \eens
Thus we have on $\T_a$, by \eqref{eq::new-oracle},
\ben
\nonumber
\norm{X \delta}^2_n  + \norm{X h}^2_{n}
+ \lambda_n \norm{h_{T_0^c}}_1 
& \leq &
\norm{X \beta - X \beta_0}_n^2 +
3\lambda_n  \norm{h_{T_0}}_1 \\
\label{eq::b2}
& \leq &
2 (\norm{X \delta}^2_n  + \norm{X h}^2_{n})
+ 3\lambda_n  \norm{h_{T_0}}_1,
\een
since $(a+b)^2 \le 2 a^2 + 2b^2$.
Then we have by assumption \eqref{eq::DHhead} and \eqref{eq::b2},
\ben
\label{eq::oracle2}
\lambda_n \norm{h_{T_0^c}}_1 
& \leq &
\norm{X \delta}^2_n
+ \norm{X h}^2_{n} + 3\lambda_n  \norm{h_{T_0}}_1 \leq 
4 \lambda_n \norm{h_{T_0}}_1.
\een
Then $h \in \Cone(s_0, 4)$. 
Hence,  under $\RE(s_0, 4, X)$ condition,  we have by assumption
\eqref{eq::DHhead}, \eqref{eq::oracle2} and Definition~\eqref{eq::admissible}, 
\ben
\nonumber
\norm{X \delta}^2_n + 
\norm{X h}^2_{n}
& \leq &
\lambda_n  \norm{h_{T_0}}_1 \le
\lambda_n  \sqrt{s_0} \twonorm{h_{T_0}}\\
\label{eq::Xh2}
& \leq &
\lambda_n  \sqrt{s_0} K(s_0, 4) \norm{X h}_{n} \\
\label{eq::Xd2}
& \leq &
\norm{X h}^2_{n} + (\lambda_n  \sqrt{s_0} K(s_0, 4) /2)^2,
\een
from which \eqref{eq::Xdelta2} and \eqref{eq::HT0case2} immediately
follow since
\ben
\norm{X \delta}_n 
& \leq &
\lambda_n  \sqrt{s_0} K(s_0, 4)/2 \;   \text{ by~\eqref{eq::Xd2}},\\
\label{eq::Xhtwo}
\norm{X h}_{n} & \leq &
\lambda_n  \sqrt{s_0} K(s_0, 4) \;  \text{ by~\eqref{eq::Xh2}, and
  hence} \; \\
\label{eq::HT0case2}
\twonorm{h_{T_0}}
& \leq & K(s_0, 4) \norm{X h}_{n} \le \lambda_n  \sqrt{s_0} K(s_0, 4)^2.
  \een
Moreover, we have for $h \in \Cone(s_0, 4)$,
\ben
\nonumber
\norm{h_{T_0}}_1
& \leq & \sqrt{s_0} \twonorm{h_{T_0}} \le \lambda_n s_0 K(s_0, 4)^2, \\
\label{eq::TailCase2}
\norm{h_{T_0^c}}_1 
& \leq &
4 \lambda_n K(s_0, 4)^2 s_0  =: D_{1,b} \lambda \sigma s_0, \; \text{ and } \\
\label{eq::H1normC2}
\norm{h}_1 
& \leq &
{5} \lambda_n K(s_0, 4)^2 s_0 = {5 d_0}
K(s_0, 4)^2 \lambda \sigma s_0 =: D_{2,b} \lambda \sigma s_0,
 \een
 where $D_{1,b} = 4 d_0  K(s_0, 4)^2$ and $D_{2,b} = 5 d_0  K(s_0,
 4)^2$.
\noindent{\bf Case 3.}
Suppose
\ben
\label{eq::case3pre}
\lambda_n
\norm{h_{T_0}}_1 \le
\norm{X\delta}^2_n  + \norm{X h}^2_{n}
\le \norm{X \beta - X \beta_0}_n^2.
\een
Thus we have on $\T_a$,  by the triangle inequality,
\eqref{eq::case3pre}, and~\eqref{eq::newbound2},
\bens
\label{eq::newbound}
\lambda_n \norm{h_{T_0^c}}_1 + \lambda_n  \norm{h_{T_0}}_1
& \leq &
\norm{X \delta}^2_n  + \norm{X h}^2_{n}
+ \lambda_n \norm{h_{T_0^c}}_1 \\
& \leq &
\norm{X \beta - X \beta_0}_n^2 + 3\lambda_n  \norm{h_{T_0}}_1.
\eens
Then
\ben
\lambda_n \norm{h_{T_0^c}}_1 
& \leq &
\norm{X \beta - X \beta_0}_n^2 + 2 \lambda_n  \norm{h_{T_0}}_1, \;
\text{ and} \\
\label{eq::h1case3}
\lambda_n \norm{h}_1 
& \leq &
\norm{X \beta - X \beta_0}_n^2 +
3 \lambda_n  \norm{h_{T_0}}_1 \le 4 \norm{X \beta - X \beta_0}_n^2.
\een
Now by~\eqref{eq::h1case3}, we have for 
$d_0 \geq 2  \sqrt{1 + a}$ and $h = \hat\beta -\beta_0$,
\ben 
\label{eq::H1normC3}
\norm{h}_1 & \leq  & 
\norm{ X \beta - X \beta_0 }^2_{n}/(\lambda_n) + 3  \norm{h_{T_0}}_1 
\leq  4 \norm{X \beta - X \beta_0 }_n^2/(\lambda_n) \\
\nonumber 
& = & 
4 \Lambda_{\max}(s- s_0) \lambda \sigma s_0/(d_0).
\een 
\noindent{\bf Putting things together.}
We will show in the proof of Lemma~\ref{lemma::HT01case3} that for
for {\bf Case 3},  the following $\ell_2$-norm error bound for $h =
\hat\beta-\beta_0$,
\bens
\twonorm{h_{T_{0}}}   & \leq &
\twonorm{h_{T_{01}}}  \le
D \lambda \sigma \sqrt{s_0} \; \text{ for $D$
as in~\eqref{eq::defineD}};
\eens
Combining the preceding bound with \eqref{eq::HT0case2} for {\bf Case
  2} and \eqref{eq::HT0case1}
for {\bf Case 1},
we have the expression for \eqref{eq::D0-prime} and the following
$\ell_2$ error bound on  ${h_{T_{0}}}$,
\bens
\twonorm{h_{T_{0}}}   & \leq &
D_0'\lambda \sigma \sqrt{s_0} \; \; \text{ where }\\
D_0' &=& \left\{D \vee [d_0 K(s_0, 4)^2] \vee
 [(K(s_0, 3) \sqrt{\Lambda_{\max}(s- s_0)} + 3 d_0 K^2(s_0, 3))]\right\}.
\eens

In summary,
for {\bf Case 1},
we have the following bounds in the $\ell_1$ norm:
\bens
\label{eq::case1TOCdir}
\norm{h_{T_0^c}}_1
& \leq &
\norm{X (\beta-\beta_0)}_n^2 /(\lambda_n) + (3 K(s_0, 3)/2)^2 \lambda_n s_0,\\
 & \leq &
[\frac{\Lambda_{\max}(s- s_0)}{d^2_0} + 9 K(s_0, 3)^2/4]
d_0 \lambda \sigma s_0.
\eens
and
for {\bf Case 2},
we have by \eqref{eq::cone4}, 
\bens 
\norm{h_{T_0^c}}_1 
& \leq &
4 d_0 K(s_0, 4)^2 \lambda \sigma s_0  =: D_{1, b} \lambda \sigma s_0;
\eens
Finally, for {\bf Case 3}, we have by \eqref{eq::HT01Case3},
\bens
\norm{h_{T_0^c}}_1 & \le &
3 \norm{ X \beta - X \beta_0 }^2_{n}/\lambda_n \le D_{1,c}
\lambda \sigma s_0 \; \\
\text{ where}\; \;
D_{1,c} & := & 3 \Lambda_{\max}(s- s_0)/d_0.
\eens
Moreover, combining 
\eqref{eq::TailboundMain},~\eqref{eq::HT0case1a},~\eqref{eq::cone4}, 
and \eqref{eq::HT01Case3}, we have for $D_1$ as in~\eqref{eq::D1-define}:
\bens
\norm{h_{T_0^c}}_1 
& \leq &
\lambda \sigma s_0
\left\{\big(D_{1, a}
  \wedge 3 D'_0 ]
\vee 4 d_0 K(s_0, 4)^2 \vee \frac{3 \Lambda_{\max}(s- s_0) }{d_0}
\right\} \\
& \leq &
d_0 \lambda \sigma s_0
\left\{4 K(s_0, 4)^2 + \frac{3 \Lambda_{\max}(s- s_0) }{d^2_0}
\right\}
\eens
Similarly, combining \eqref{eq::H1normC1}, \eqref{eq::H1normC2}
for {\bf Case   2}, and~\eqref{eq::H1normC3} for {\bf Case  3}, we obtain
the expression of $D_2$ in~\eqref{eq::D2-define}:
\bens 
\norm{h}_1 
& \leq &
\lambda \sigma s_0 \left\{\big(\frac{\Lambda_{\max}(s- s_0)}{d_0} + 4 d_0 K^2(s_0, 
 3)\big) \vee 5 d_0 K(s_0, 4)^2 \vee \frac{4 \Lambda_{\max}(s- s_0)
 }{d_0} \right\}.
\eens

\section{Proof of Lemma~\ref{lemma::HT01case3}}
\label{sec::suppHT01proof}
\begin{proofof2}
Let $h = \beta_{\init} - \beta^{\ext}({T_0})$.  
Throughout this proof,  as a shorthand, we write $h = \hat{\beta} - \beta_0$, where $\hat\beta  =
\beta_{\init}$ and $\beta_0 = \beta^{\ext}({T_0})$ as defined in Lemma~\ref{lemma:threshold-RE}.  
Under the settings of   Theorem~\ref{thm::RE-oracle}, we have on $\T_a$,
\ben 
\label{eq::newbound2}
\norm{X\delta }^2_n + \norm{X h}^2_{n}
+ \lambda_n \norm{h_{T_0^c}}_1 
\leq \norm{X (\beta - \beta_{0}) }^2_n 
+ 3 \lambda_n  \norm{h_{T_0}}_1.
\een
The rest of the proof is devoted to the $\ell_2$ error bound on
$h_{T_{01}}$ in view of  Lemma~\ref{lemma::h01-bound-CT}.
Let $T_1$ be the $s_0$ largest positions of $h$ outside of $T_0$. 
Denote by $J_0$ the locations of the $s_0$ largest  coefficients of $h$ in 
absolute values. Then $J_0 \subset T_{01}$.
\noindent{\bf Case 1}
Now by Proposition A.1 as derived by~\cite{Zhou09c},
\eqref{eq::newbound2},~\eqref{eq::HT0case1a}, and
Lemma~\ref{lemma:T0-pre-loss}, we have  for $K = K(s_0, 3)$,
\bens
\twonorm{h_{T_{01}}}
& \leq & \sqrt{2} K(s_0, 3) \norm{X h}_{n} \\
& \leq &
\sqrt{2} K
\big(\norm{X \beta - X \beta_0}_n^2 
+ 3 \lambda_n \sqrt{s_0} \twonorm{h_{T_0}} \big)^{1/2} \\
& \leq &
\sqrt{2} K(s_0, 3) 
\big(\norm{X \beta - X \beta_0}^2_n 
+ 3 \lambda_n \sqrt{s_0} K \norm{X \beta - X \beta_0}_{n}  + 9K^2
\lambda^2_n s_0) \big)^{1/2} \\
& \le &
\sqrt{2} K(s_0, 3) \big(\big(\norm{X \beta - X \beta_0}_n 
+ 3 \lambda_n \sqrt{s_0} K \big)^2 \big)^{1/2} \\
& \le &
D_{0,a}
\lambda \sigma \sqrt{s_0},
\eens
where
\bens
D_{0,a} :=
\sqrt{2} \big(\sqrt{\Lambda_{\max}(s-s_0)} K(s_0, 3) 
+ 3 d_0 K^2(s_0, 3) \big).
\eens
By Lemma~\ref{lemma::h01-bound-CT},
\bens 
\twonorm{h}^2 
& \leq &
\twonorm{h_{T_{01}}}^2 + \norm{h_{T_0^c}}^2_1/s_0  \leq
[D_{0,a}^2 + D_1^2]
\lambda^2 \sigma^2 s_0 \\
\text{ and } \; \twonorm{h}^2 
& \leq & (1 + k_0) \twonorm{h_{T_{01}}}^2, \; \text{ where} \; h \in \Cone(s_0, k_0),
\eens
when $\RE(s_0, k_0, X)$ holds.
To see this, notice that we have  by Lemma~\ref{lemma::cone}, 
where we set $k_0 = 3$, 
\bens
\twonorm{h} & \leq &
\sqrt{(1+3)} \twonorm{h_{J_0}}  \le
2 \twonorm{h_{T_{01}}} \le 2 D_{0,a} \lambda \sigma \sqrt{s_0}.
\eens

\noindent{\bf Case 2.}
Similar to Case 1, we have by~\eqref{eq::Xhtwo},
\bens
\twonorm{h_{T_{01}}}
& \le & \sqrt{2} K(s_0, 4) \norm{X h}_n  \le
\sqrt{2} K^2(s_0, 4)  \lambda_n  \sqrt{s_0}.
\eens
since $h \in \Cone(s_0, 4)$.
Hence, we have by Lemma~\ref{lemma::cone}, where we set $k_0 = 4$,
\bens
\twonorm{h} & \leq &
\sqrt{(1+4)} \twonorm{h_{J_0}}  \le
\sqrt{(1+4)} \twonorm{h_{T_{01}}} \\
& \leq &
\sqrt{5} \sqrt{2} K^2(s_0, 4)  \lambda_n  \sqrt{s_0} =
\sqrt{10} d_0 K(s_0, 4)^2 \lambda \sigma \sqrt{s_0}.
\eens

\noindent{\bf Case 3.}
Denote by 
\ben 
\label{eq::defineDlocal}
D & := & \frac{\sqrt{\Lambda_{\max}(s - s_0)}}{\sqrt{\Lambda_{\min}(2s_0)}}
\left(1+ \frac{3 \ell(s_0) \sqrt{\Lambda_{\max}(s - s_0)}}{d_0} \right). 
\een
By Lemmas~\ref{lemma::h01-bound-CT},
~\ref{lemma::T01simple},~\ref{lemma:T0-pre-loss},
Assumption~\eqref{eq::case3pre}, and~\eqref{eq::HT01Case3},
we obtain for  $h = \hat\beta - \beta_0$, 
\bens
\twonorm{h_{T_{01}}} 
& \leq &
\nonumber 
\inv{\sqrt{\Lambda_{\min}(2s_0)}}
\big(\norm{X h}_{n} +  \ell(s_0) \norm{h_{T_0^c}}_1 /\sqrt{s_0} \big)\\
\nonumber
& \leq &
\inv{\sqrt{\Lambda_{\min}(2s_0)}} \norm{X \beta - X
  \beta_0 }_{{n}} + \frac{\ell(s_0)}{\sqrt{\Lambda_{\min}(2s_0)}}
\frac{3\norm{X \beta - X \beta_0}^2_n}{\lambda_n\sqrt{s_0} } \\
\nonumber
& \leq &
\frac{\norm{X \beta - X \beta_0}_n }{\sqrt{\Lambda_{\min}(2s_0)}}
\big(1 + \frac{3\norm{X \beta - X \beta_0}_n \ell(s_0)}{\lambda \sigma
  d_0 \sqrt{s_0} } \big)\\
& \leq &
\sqrt{\frac{\Lambda_{\max}(s - s_0)}{\Lambda_{\min}(2s_0)} }
\left(1+ \frac{3 \ell(s_0) \sqrt{\Lambda_{\max}(s - s_0)}}{d_0}
\right) \lambda \sigma \sqrt{s_0}
=: D \lambda \sigma \sqrt{s_0}
\eens
where recall for {\bf Case 3},
\bens
\label{eq::simple1}
\norm{h_{T_0^c}}_1/\sqrt{s_0} &\leq & 3 \norm{X \beta - X
  \beta_0}_n^2/(\lambda_n \sqrt{s_0})
\leq  D_{1,c} \lambda \sigma \sqrt{s_0}, \\
& = &   3 \Lambda_{\max}(s- s_0) \lambda \sigma \sqrt{s_0} /d_0 \\
 \text{ and } \; \norm{X h}_{n}
& \leq &  \norm{X (\beta - \beta_0)}_{n} \le  
\sqrt{\Lambda_{\max}(s- s_0)} \lambda \sigma \sqrt{s_0}.
\eens
Combining all three cases, we obtain the expression for
\eqref{eq::D0-define-orig}:
\bens
\twonorm{h_{T_{01}}}
& \le & D_0  \lambda \sigma  \sqrt{s_0} \text{ where}\; \;
D_0  := D_{0,a} \vee [d_0 \sqrt{2} K^2(s_0, 4)]  \vee D 
\eens 
where $D$ is as defined in \eqref{eq::defineDlocal} and $D_{0,a}$ as in
\eqref{eq::TaDa}.

It remains to bound $\twonorm{h}$ for {\bf Case 3}.
Thus we have on $\T_a$ for $\lambda_n= d_0 \lambda \sigma$
and $\ell(s_0)$ as defined in~\eqref{eq::definell},
for $D$ as in~\eqref{eq::defineDlocal}, by Lemma~\ref{lemma::h01-bound-CT},
\bens 
\twonorm{h}^2 
& \leq &
\twonorm{h_{T_{01}}}^2
+ \norm{h_{T_0^c}}^2_1/s_0 \\
& \leq &
\nonumber 
\big(\inv{\sqrt{\Lambda_{\min}(2s_0)}}
\big(\norm{X h}_{n} +
\ell(s_0) \norm{h_{T_0^c}}_1 /\sqrt{s_0}
\big)\big)^2 + \norm{h_{T_0^c}}^2_1/s_0 \\
& \leq & (D^2 + D_{1,c}^2) \lambda^2 \sigma^2  s_0, \;
\text{where}  \;
D_{1,c} =  \frac{3 \Lambda_{\max}(s- s_0)}{d_0},
\eens
Finally,  for $h = \hat{\beta} - \beta_0$, where $\beta_0 = \beta^{\ext}({T_0})$,
$D_0$ as defined in~\eqref{eq::D0-define-orig},  and  
$D_1$ as defined in~\eqref{eq::D1-define}, we have by the triangle
inequality and Lemma~\ref{lemma::h01-bound-CT}, 
\bens
\twonorm{\hat\beta - \beta}
& \leq &
\twonorm{\hat \beta - \beta_{0}}
+  \twonorm{\beta - \beta_{0}}  \leq
\twonorm{h} + \lambda \sigma \sqrt{s_0} \\
& \leq & 
(\twonorm{h_{T_{01}}}^2
+ \norm{h_{T_0^c}}^2_1/s_0)^{1/2} + 
\lambda \sigma \sqrt{s_0} \\ 
& \leq &
[\sqrt{D_0^2 + D_1^2}+ 1] \lambda \sigma \sqrt{s_0},
\eens
over all cases.
Moreover, 
we have by Lemmas~\ref{lemma::h01-bound-CT} and~\ref{lemma::cone},  for both
{\bf Case 1} and {\bf Case 2}, the following slightly stronger bound:
\bens
\twonorm{\hat\beta - \beta}
& \leq &
\twonorm{h} + \lambda \sigma \sqrt{s_0} \\
& \leq & 
(\twonorm{h_{T_{01}}}^2 + \norm{h_{T_0^c}}^2_1/s_0)^{1/2} \wedge
\sqrt{1 + k_0} \twonorm{h_{T_{01}}}+ \lambda \sigma \sqrt{s_0} \\
& \leq & 
\lambda \sigma \sqrt{s_0} [(\sqrt{D_{0}^2 + D_1^2} + 1)
\wedge   (\sqrt{5} D_{0} + 1)].
\eens
\end{proofof2}

\section{Proof of  Lemma~\ref{lemma::T01simple}}
\label{sec::suppT01simpleproof}
\begin{proofof2}
  For $X h_{T_{01}} = Xh - \sum_{j \geq 2} X h_{T_j}$, we have by the
  triangle inequality and sparse eigenvalue condition,
  \bens
  \twonorm{X h_{T_{01}}}/\sqrt{n}
& \leq & 
\twonorm{X h}/\sqrt{n} + \sum_{j \geq 2} \twonorm{X h_{T_j} }
/\sqrt{n} \\
 & \le &
 \twonorm{X h} /\sqrt{n}+  \sqrt{\Lambda_{\max}(s_0)} \sum_{j \geq 2}
 \twonorm{h_{T_j}} \\
 & \le &
 \twonorm{X h} /\sqrt{n}+  
 \sqrt{\Lambda_{\max}(s_0)} \sum_{j \geq 1}
 \norm{h_{T_j}}_1 /{\sqrt{s_0}}\\
 & \le &
 \twonorm{X h} /\sqrt{n} +  \sqrt{\Lambda_{\max}(s_0)} \norm{h_{T_0^c}}_1 /{\sqrt{s_0}}.
 \eens
 Thus it follows from the proof Lemma~\ref{lemma::h01-bound-CT} that
 \ben
 \nonumber
 \twonorm{h_{T_{01}}} & \le &
 \twonorm{X h_{T_{01}}}/ \sqrt{n \Lambda_{\min}(2s_0)} \\  
 & \le &
 \label{eq::simpleT01alt}
  \inv{\sqrt{\Lambda_{\min}(2s_0)}}\left( \twonorm{X h} /\sqrt{n}+\sqrt{\Lambda_{\max}(s_0)} \norm{h_{T_0^c}}_1 /{\sqrt{s_0}}\right),
 \een
 where we replace $\frac{\theta_{s_0, 2s_0}}{\sqrt{\Lambda_{\min}(2s_0)}}$ with
 $\sqrt{\Lambda_{\max}(s_0)}$.
The lemma thus holds in view of \eqref{eq::simpleT01orig} and  \eqref{eq::simpleT01alt}.
\end{proofof2}

\section{Nearly ideal model selections under the UUP}
\label{sec:DS-threshold}
The Dantzig selector~\citep{CT07} is defined as follows:  
for some  $\lambda_n \geq 0$,
\begin{gather}
\label{eq::DS-func}
(DS) \; \; \arg\min_{\tilde{\beta} \in \R^p} \norm{\tilde{\beta}}_1
 \;\; \text{subject to} \;\; 
\norm{X^T (Y - X \tilde{\beta})/n}_{\infty} \leq \lambda_n.
\end{gather}
~\cite{Zhou09th} shows that thresholding  
of an initial Dantzig selector $\beta_{\init}$ at the level of 
$\asymp \sqrt{2 \log p/n} \sigma$ followed by OLS refitting, 
achieves the {\bf sparse oracle inequalities} under a UUP.
\begin{definition}
\textnormal{(\bf{A Uniform Uncertainly Principle})}
\label{def:CT-cond}
For some integer $1 \leq s < n/3$, assume $\delta_{2s} + \theta_{s,
  2s} < 1-\tau$ for some $\tau>0$,
which implies that $\Lambda_{\min}(2s) > \theta_{s, 2s}$.
\end{definition}

\begin{remark}
It is clear that $\delta_{2s} < 1$ implies that the sparse eigenvalues 
condition~\eqref{eq::eigen-admissible-s} and~\eqref{eq::eigen-max} hold.
Moreover, Assumption~\ref{def:CT-cond} implies 
that $\RE(s_0, k_0, X)$ as in~\eqref{eq::admissible} hold for  $s_0
\le s$ with
\ben
\label{eq::RE-UUP}
K(s_0, k_0) \leq K(s, k_0) = 
\frac{\sqrt{\Lambda_{\min}(2s)}}{\Lambda_{\min}(2s) - \theta_{s, 2s}}
 \leq 
\frac{\sqrt{\Lambda_{\min}(2s)}}{1 - \delta_{2s} -  \theta_{s, 2s}}
 \leq 
\frac{\sqrt{\Lambda_{\min}(2s)}}{\tau}
\een
for $k_0 = 1$, as $K(s, k_0)$ is nondecreasing with respect to $s$ for the same $k_0$; 
see~\cite{BRT09}.
\end{remark}

\noindent{{\bf The Gauss-Dantzig Procedure}}: Assume $\delta_{2s} + \theta_{s,2s} < 1 - \tau$, 
where $\tau > 0$:
\begin{itemize}
\item[Step 1]
\label{step::thresh}
First obtain an initial estimator $\beta_{\init}$ using the Dantzig 
selector in~\eqref{eq::DS-func} with 
$\lambda_n = (\sqrt{1+a} + \tau^{-1}) \sqrt{2 \log p/n} \sigma
=:  \lambda_{p, \tau} \sigma$, where 
$a \geq 0$; then threshold $\beta_{\init}$ with $t_0$, chosen
from the range 
$(C_1 \lambda_{p, \tau} \sigma, C_4 \lambda_{p, \tau} \sigma]$, 
for $C_1$ as defined in~\eqref{eq::DS-constants-1};
set $I :=  \left\{j \in \{1, \ldots, p\}: \beta_{j, \init} \geq t_0  \right\}.$
\item[Step 2]
Given a set $I$ as above, construct the estimator
$\hat\beta_{I} = (X_I^T X_{I})^{-1} X_{I}^T Y$ and set $\hat\beta_j =
0, \forall j \not\in I.$
\end{itemize}

\begin{theorem}\textnormal{(\bf{Variable selection under UUP})}
\label{thm:ideal-MSE-prelude}
Choose $\tau, a > 0$ and set $\lambda_n = \lambda_{p, \tau} \sigma$,
where $\lambda_{p,\tau} := (\sqrt{1+a} + \tau^{-1}) \sqrt{2 \log p/n}$,
in~\eqref{eq::DS-func}.
Suppose $\beta$ is $s$-sparse with $\delta_{2s} + \theta_{s,2s} < 1 - \tau$.
Let threshold $t_0$ be chosen from the range 
$(C_1\lambda_{p,\tau} \sigma, C_4 \lambda_{p, \tau} \sigma]$ for some
constants $C_1, C_4$ to be defined.
Then  the  Gauss-Dantzig selector $\hat{\beta}$ selects a model 
$I := \supp(\hat{\beta})$ such that
we have 
\ben
\label{eq::ideal-MSE}
|I| & \leq &  2s_0 , \quad |I \setminus S| \leq s_0 \leq s \quad
\text{ and } \quad \shtwonorm{\hat\beta - \beta}^2  \leq 
C_3^2 \lambda^2 \sigma^2 s_0
\een
with probability at least $1 - (\sqrt{\pi \log p} p^a)^{-1}$,
where $C_1$ is defined in~\eqref{eq::DS-constants-1} and 
$C_3$ depends on $a, \tau$, $\delta_{2s}$, $\theta_{s, 2s}$ and $C_4$; see~\eqref{eq::DS-constants-3}.
\end{theorem}
Our analysis for Theorem~\ref{thm:ideal-MSE-prelude} builds 
upon Proposition~\ref{prop:DS-oracle}~\citep{CT07}, which shows 
the Dantzig selector 
achieves the oracle inequality as stated in~\eqref{eq::log-MSE} under 
Assumption~\ref{def:CT-cond}. We note that, 
in Assumption~\ref{def:CT-cond}, the sparsity level is fixed at $s$
rather than $s_0$.
Hence it is stronger than the conditions we impose in
Theorem~\ref{thm::RE-oracle-main} for the Thresholded Lasso estimator.
We now show the oracle inequalities for the Dantzig selector. 
We then show in the supplementary Lemma~\ref{lemma:threshold-DS} that
thresholding at the  level of $\sigma \lambda$ as elaborated in
Step~1 in the Gauss-Dantzig Procedure selects a set $I$ of at most $2s_0$ variables, among which at most $s_0$ are from the complement of the support set $S$ as required in~\eqref{eq::ideal-MSE}.
\begin{proposition}
\textnormal{\citep{CT07}}
\label{prop:DS-oracle}
Let $Y = X \beta + \e$, for $\e$ being i.i.d. $N(0, \sigma^2)$ and
$\twonorm{X_j}^2 = n$.
Choose $\tau, a > 0$ and set $\lambda_n = (\sqrt{1+a} + \tau^{-1}) \sigma \sqrt{2 \log p/n}$ in~\eqref{eq::DS-func}. 
Then if $\beta$ is $s$-sparse with $\delta_{2s} + \theta_{s,2s} < 1 - \tau$, 
the Dantzig selector obeys with probability at least $1 - (\sqrt{\pi
  \log p} p^a)^{-1}$,
\bens
\shtwonorm{\hat\beta - \beta}^2  & \leq &   
C_2^2 (\sqrt{1+a} + \tau^{-1})^2 s_0 \lambda^2\sigma^2
\eens
\end{proposition}
From this point on we let 
$\delta := \delta_{2s}$ and $\theta := \theta_{s, 2s}$;
Analysis by~\cite{CT07} (Theorem 2) and the current paper yields the 
following constants,
\begin{eqnarray}
\label{eq::DS-constants}
C_2 & = & 2 C_0' + \frac{1 + \delta}{1 - \delta - \theta} 
\; \text { where } 
C'_0 = \frac{C_0}{1 - \delta - \theta} + 
\frac{\theta(1 + \delta)}{(1 - \delta - \theta)^2}, 
\end{eqnarray}
where 
$C_0 =  2 \sqrt{2}
\left(1 + \frac{1 - \delta^2}{1 - \delta - \theta}\right) 
+ (1 + 1/\sqrt{2})\frac{(1 + \delta)^2}{1 - \delta - \theta}$.
We now define
\begin{eqnarray}
\label{eq::DS-constants-1}
C_1 & = & C_0' + \frac{1+ \delta}{1-\delta-\theta}  \text{ and } \\
\label{eq::DS-constants-3}
  C_3^2 & = &  3 (\sqrt{1+a} + \tau^{-1})^2 ((C_0' + C_4)^2 +1 )
+ {4(1+a)}/{ \Lambda^2_{\min}(2s_0)},
\end{eqnarray}
where $C_3$ is used in~\eqref{eq::ideal-MSE} and 
has not been optimized in our  analysis.

\section{Proof of Theorem~\ref{thm:ideal-MSE-prelude}}
\label{sec:append-gauss}

Now similar to Lemma~\ref{lemma:threshold-RE},
Lemma~\ref{lemma:threshold-DS} bounds the size of $I$, as well as the 
bias that we introduce to model $I$ by thresholding.
Theorem~\ref{thm:ideal-MSE-prelude} is an immediate 
corollary of Lemmas~\ref{prop:MSE-missing} and~\ref{lemma:threshold-DS}.
The proof follows from
Lemma~\ref{lemma:threshold-DS} and Proposition~\ref{prop:DS-oracle}. We
include its proof in Section~\ref{sec::proofofUUP} for self-containment.

\begin{lemma}
\label{lemma:threshold-DS}
Choose $\tau > 0$ such that $\delta_{2s} + \theta_{s,2s} < 1 - \tau$.
Let $\beta_{\init}$ be the solution to~\eqref{eq::DS-func} with 
$\lambda_n =  \lambda_{p,\tau} \sigma 
:= (\sqrt{1+a} + \tau^{-1}) \sqrt{2 \log p/n} \sigma$.
Given some constant $C_4 \geq C_1$, for $C_1$ as 
in~\eqref{eq::DS-constants-1}, choose a thresholding parameter $t_0$ 
such that 
$C_4 \lambda_{p,\tau} \sigma \geq t_0  > C_1 \lambda_{p,\tau} \sigma$
and set $I = \{j: \size{\beta_{j, \init}} \geq t_0\}.$
Then with probability at least $1 - (\sqrt{\pi \log p} p^a)^{-1}$,
we have~\eqref{eq::ideal-MSE} and 
$\twonorm{\beta_{\drop}} \leq \sqrt{(C_0' + C_4)^2 + 1} 
\lambda_{p,\tau} \sigma \sqrt{s_0}$, where
$\drop := \{1, \ldots, p\} \setminus I$
and $C_0'$ is defined in~\eqref{eq::DS-constants}.
\end{lemma}

\begin{proofof}{Theorem~\ref{thm:ideal-MSE-prelude}}
It holds by definition of $S_{\drop}$ that $I \cap S_{\drop} = \emptyset$.
It is clear by Lemma~\ref{lemma:threshold-DS} that
$|\dropS| < s$ and $|I| \leq 2s_0$ and
$|I \cup S_{\drop}| \leq |I \cup S| \leq s + s_0 \leq 2s$;
Thus for $\hat\beta_{I} = (X_I^T X_{I})^{-1} X_{I}^T Y$,
$\hat\beta_{I^c} =0$ and  $\lambda
= \sqrt{2 \log p/n}$, we have by Lemma~\ref{prop:MSE-missing-orig}
\bens
\label{eq::re-exp}
\lefteqn{
\twonorm{\hat{\beta} - \beta}^2
 \leq 
\twonorm{\beta_{\drop}}^2
\left(1 + \frac{2 \theta^2_{s, 2s_0}}{\Lambda_{\min}^2(2s_0)}\right) +
\frac{4 s_0}{\Lambda_{\min}^2(2s_0)}\basepen^2 }\\
& \leq &
\lambda^2 \sigma^2 s_0
\left(\sqrt{1+a} + \tau^{-1})^2 ((C_0' + C_4)^2 + 1) 
\big(1 + \frac{2 \theta^2_{s, 2s_0}}{\Lambda_{\min}^2(2s_0)}\big)
+ \frac{4(1+a)}{\Lambda^2_{\min}(2s_0)}\right) \\
& \leq &
C_3^2 \lambda^2 \sigma^2 s_0
\eens
with probability at least $1 - (\sqrt{\pi \log p} p^a)^{-1}-\exp(-3m/64)$ where $m = \size{I}$.
Thus the theorem holds for $C_3$ as in~\eqref{eq::DS-constants-3}, where it holds for $\tau > 0$ that
\bens
\frac{\theta_{s, 2s_0}}{\Lambda_{\min}(2s_0)}
 \leq
\frac{\theta_{s, 2s}}{\Lambda_{\min}(2s_0)}  
 \leq
\frac{1 - \delta_{2s} - \tau}{\Lambda_{\min}(2s)} < 1
\eens
given that 
$\theta_{s, 2s} < 1 - \tau - \delta_{2s} < \Lambda_{\min}(2s)$ for $\tau > 0$.
\end{proofof}

\subsection{Proof of Lemma~\ref{lemma:threshold-DS}}
\label{sec::proofofUUP}

\begin{proofof2}
Suppose $\T_a$ holds.
Consider the set 
$I \cap T_0^c := \{j \in T_0^c: \size{\beta_{j, \init}} > t_0\}$.
It is clear by definition of $h = \beta_{\init} - \beta^{(1)}$ 
and~\eqref{eq::DS-T0c-1-bound} that
\ben
\label{eq::set-count-T0-c}
\size{I \cap T_0^c} \leq \norm{\beta_{T_0^c, \init}}_1/t_0 
= \norm{h_{T_0^c}}_1/t_0 < s_0,
\een
where $t_0 \geq C_1 \lambda_{p, \tau} \sigma$.
Thus $|I| = |I \cap T_0| + |I \cap T_0^c| \leq 2s_0$;
Now~\eqref{eq::ideal-MSE} holds given~\eqref{eq::set-count-T0-c} and
$|I  \cup S| = |S| + |I \cap \Sc| \leq s +  |I \cap T_0^c| < s+ s_0.$
We now bound $\twonorm{\beta_{\drop}}^2$.
By~\eqref{eq::DS-T01-2-bound} and~\eqref{eq::off-beta-norm-bound-2},
where $\drop_{11} \subset T_0$,
we have for $\tau < C_4 \lambda_{p, \tau} \sigma$, by the triangle inequality,
\bens
\twonorm{\beta_{\drop}}^2
& = & \twonorm{\beta_{\drop 
    \cap T_0^c}}^2   + \twonorm{\beta_{\drop \cap T_0}}^2   \le
\twonorm{\beta^{(2)}}^2 + \twonorm{\beta_{I^c \cap T_0}}^2   \\
&\leq&  s_0 \lambda^2 \sigma^2 + (t_0 \sqrt{s_0} + \twonorm{h_{T_0}})^2 \leq 
((C_4 + C_0')^2 + 1)  \lambda_{p, \tau}^2 \sigma^2 s_0.
\eens
The proof of Proposition~\ref{prop:DS-oracle}~\citep{CT07} yields 
the following on $\T_a$,
\ben
\label{eq::DS-T01-2-bound}
\twonorm{h_{T_{01}}} & \leq &  C'_0 \lambda_{p, \tau} \sigma \sqrt{s_0}, 
\text{ for }\; C'_0 \; \text{ as in}~\eqref{eq::DS-constants}, \\
\label{eq::DS-T0c-1-bound}
\norm{h_{T_0^c}}_1 & \leq & C_1 \lambda_{p, \tau} \sigma s_0, \; 
\text{ where } \; C_1 = \left(C'_0 + \frac{1 + \delta}{1 - \delta - \theta}\right), 
\; \text{ and } \\
\twonorm{h_{T_{01}^c}}
& \leq & 
\norm{h_{T_0^c}}_1/\sqrt{s_0} \; \leq \; C_1 \lambda_{p, \tau} \sigma \sqrt{s_0},
\; \text{ (cf. Lemma~\ref{lemma::h01-bound-CT}).}
\een
The rest of the proof follows that of Lemma~\ref{lemma:threshold-RE}
and hence omitted.
\end{proofof2}

\def\sleft{\hskip-5pt}
\def\lleft{\hskip-25pt}
\begin{figure}
\begin{center}
\begin{tabular}{cc}
\begin{tabular}{c}
\includegraphics[width=0.35\textwidth,angle=270]{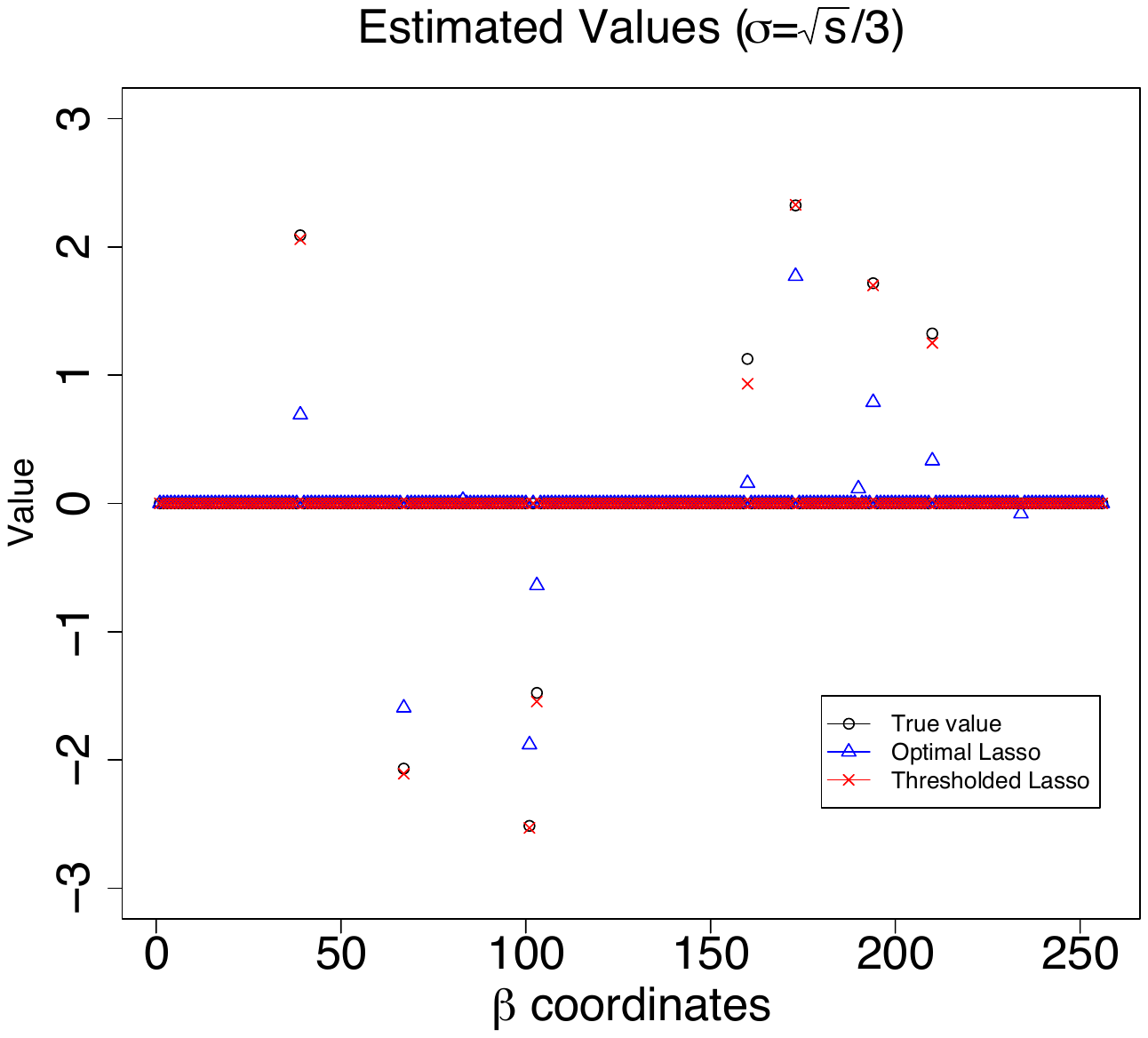}
\end{tabular}& 
\begin{tabular}{c}
\includegraphics[width=0.35\textwidth,angle=270]{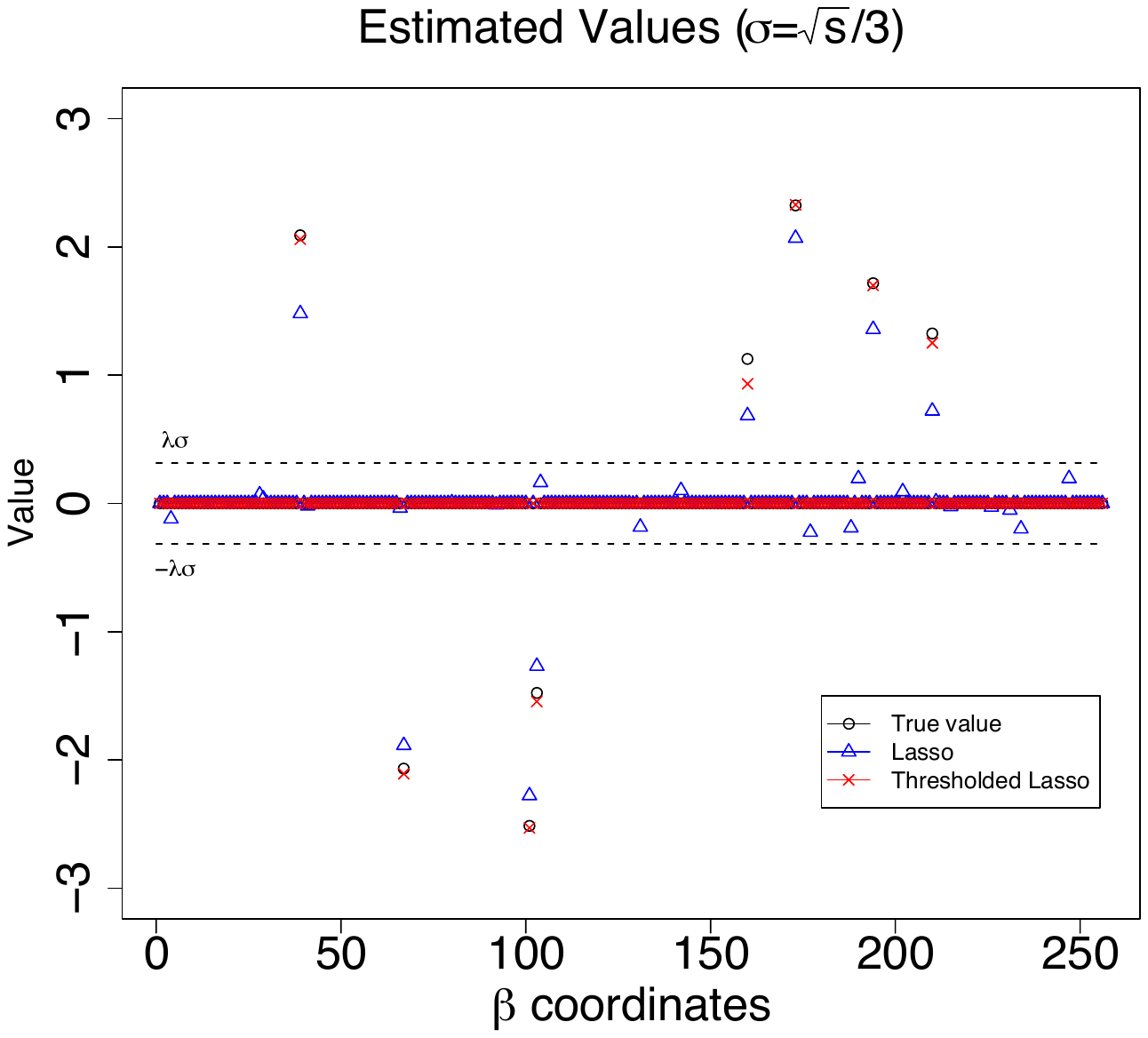} 
\end{tabular} \\
(a) & (b) \\
\end{tabular}
\caption{Illustrative example: i.i.d. Gaussian ensemble; 
$p=256$, $n=72$, $s=8$, and $\sigma = \sqrt{s}/3$.
(a) compare with the Lasso estimator $\tilde{\beta}$ which minimizes 
$\ell_2$ loss. Here $\tilde{\beta}$ has only 3 FPs, but $\rho^2$ is
large with a value of $64.73$. 
  (b) Compare with the $\beta_{\text{\rm init}}$ obtained using $\lambda_n$. 
The dotted lines show the thresholding level $t_0$. 
  The $\beta_{\text{\rm init}}$ has 15 FPs, all of which were cut after the 2nd step; 
resulting  $\rho^2= 12.73$. After refitting with OLS in the 3rd step, for the $\hat{\beta}$, $\rho^2$ is further reduced to $0.51$.}
\label{fig:example}
\end{center}
\end{figure}

\section{Numerical experiments}
\label{sec:experiments}

In this section, we present additional results from numerical simulations
designed to validate the theoretical analysis presented in this paper.

The experiment set up here is the same as the one in Section
~\ref{subsec:exp-setup}. The main difference is that all non-zero entries in
$\beta$ have large magnitudes around 1, in particular, they follow this distribution
$\beta_i = \mu_i (1 + |g_i|),$ where $\mu_i = \pm 1$ with probability 1/2 and $g_i \sim N(0,1)$.
We use $\lambda_n = 0.69 \lambda \sigma$ throughout the experiments in this
section to select a $\beta_{\init}$ as the initial estimator.  We then
threshold the $\beta_{\init}$ using a value $t_0$ typically chosen between
$0.5 \lambda \sigma$ and $\lambda \sigma$.  See each experiment for the
actual value used.  Given that columns of $X$ being normalized to have
$\ell_2$ norm $\sqrt{n}$, for each input $\beta$, we compute its
SNR as follows: $ SNR :=  \twonorm{\beta}^2 / \sigma^2.$ To evaluate
$\hat{\beta}$, we use metrics defined  in Table~\ref{tab:fpfn}; we also
compute the ratio between squared $\ell_2$ error and the ideal mean squared
error, known as the $\rho^2$; see Section~\ref{sec:exp-ell2-errors} for
details.

\subsection{Type I/II errors}
\label{subsec:type12}
We now evaluate the Thresholded Lasso estimator by comparing
Type I/II errors under different values of $t_0$ and SNR.
We consider Gaussian random matrices for the design $X$ with 
both diagonal and  Toeplitz covariance. We refer to the former
as {\it i.i.d. Gaussian ensemble} and the latter as {\it Toeplitz ensemble}.
In the Toeplitz case, the covariance is given by 
$T(\gamma)_{i,j} = \gamma^{|i-j|}$ where $0< \gamma < 1.$
We run under two noise levels: $\sigma = \sqrt{s}/3$ and $\sigma = \sqrt{s}$.
For each $\sigma$, we vary the threshold $t_0$ from 
$0.01 \lambda \sigma$ to $1.5 \lambda \sigma$. 
For each $\sigma$ and $t_0$ combination, 
we run the experiment as described in Section~\ref{subsec:exp-setup}
$200$ times with a new $\beta$ and $\epsilon$ generated each time and 
we count the number of Type I and II errors in $\hat{\beta}$.
We compute the  average at the end of 200 runs, which will correspond to
one data point on the curves in Figure~\ref{fig:type12} (a) and (b).

For both types of designs, similar behaviors are observed.
For $\sigma=\sqrt{s}/3$, FNs increase slowly; hence there is a wide
range of values from which $t_0$ can be chosen such that 
FNs and FPs are both zero. In contrast, when $\sigma=\sqrt{s}$, 
FNs increase rather quickly as $t_0$ increases due to the low SNR.
It is clear that the low SNR and high correlation combination makes
it the most challenging situation for variable selection, as predicted 
by our theoretical analysis and others.
In (c) and (d), we run additional 
experiments for the low SNR case for Toeplitz ensembles. 
The performance is improved by increasing the
sample size or lowering the correlation factor.
\begin{table}[h]
\begin{center}
\caption{Metrics for evaluating $\hat{\beta}$}
\label{tab:fpfn}
\begin{tabular}{l|l} 
\hline
Metric & Definition \\ \hline
Type I errors or False Positives (FPs) & \# of incorrectly selected non-zeros in $\hat{\beta}$ \\
Type II errors or False Negatives (FNs) & \# of non-zeros in $\beta$ that are not selected in $\hat{\beta}$ \\
True positives (TPs) & \# of correctly selected non-zeros \\ 
True Negatives (TNs) & \# of zeros in $\hat{\beta}$ that are also zero in $\beta$ \\
False Positive Rate (FPR) & $ FPR = FP / (FP + TN) = FP/(p-s) $ \\
True Positive Rate (TPR) & $ TPR = TP/ (TP+FN) = TP/ s $ \\ \hline
\end{tabular}
\end{center}
\end{table}

\subsection{$\ell_2$ loss}
\label{sec:exp-ell2-errors}
We now compare the performance of the Thresholded Lasso with 
the ordinary Lasso by examining the metric $\rho^2$ defined as follows:
$
\rho^2 = \frac{\sum_{i=1}^p (\hat{\beta}_i - \beta_i)^2}{\sum_{i=1}^p \min(\beta_i^2, \sigma^2/n)}.
$

\begin{figure}
\begin{center}
\begin{tabular}{cc}
\begin{tabular}{c}
\includegraphics[width=0.35\textwidth,angle=270]{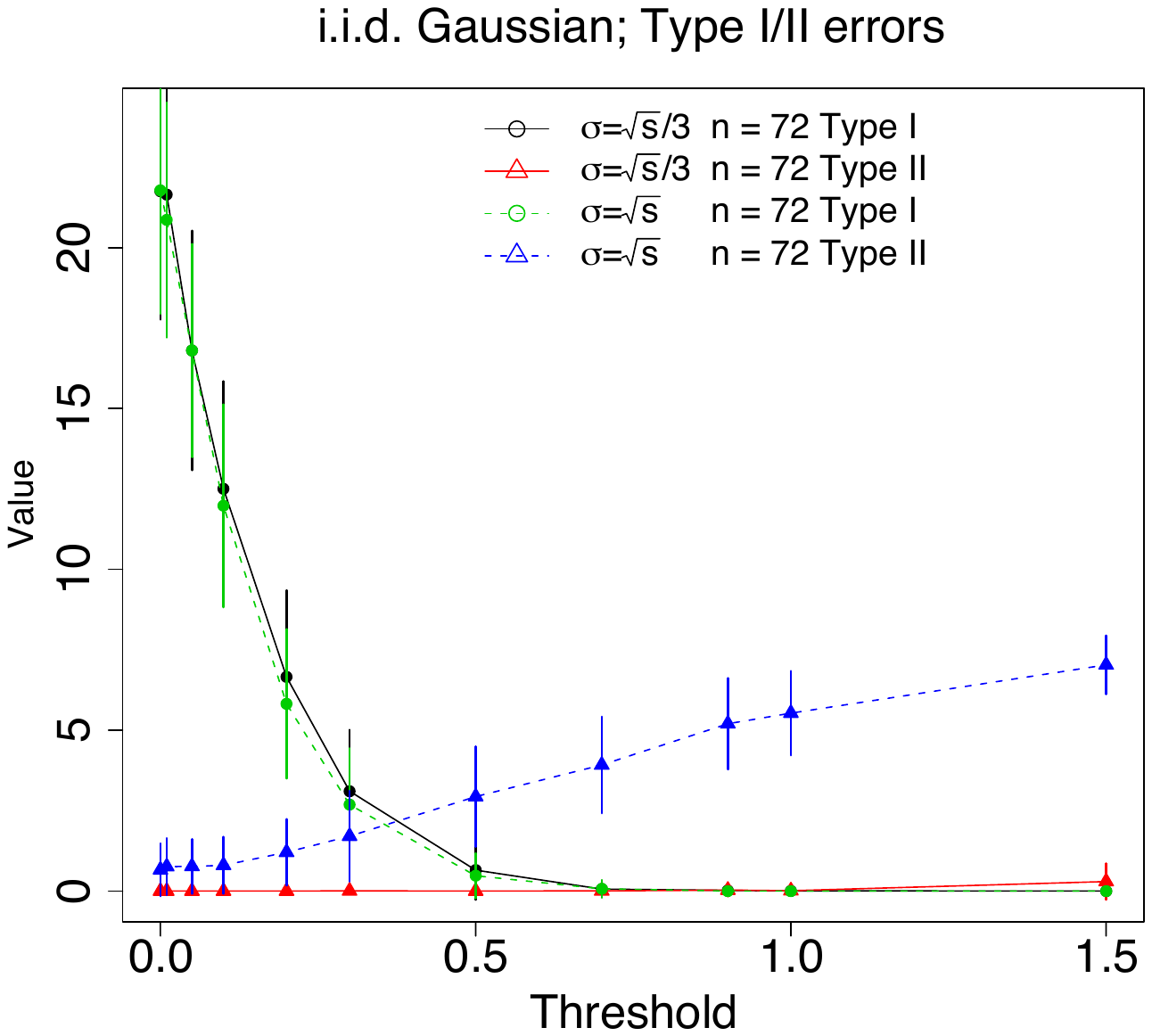}
\end{tabular}& 
\begin{tabular}{c}
\includegraphics[width=0.35\textwidth,angle=270]{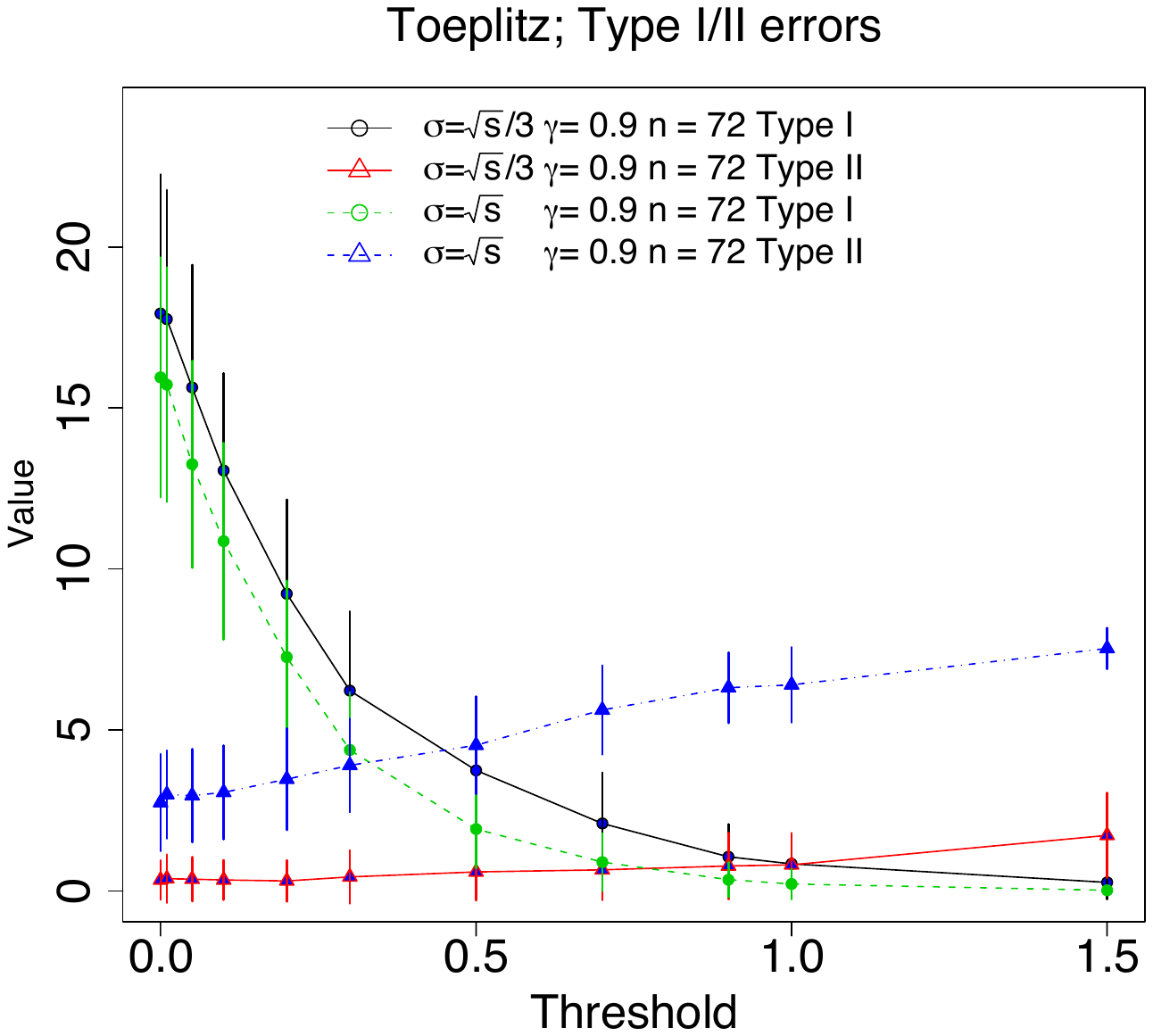} \\
\end{tabular} \\
(a) & (b) \\
\begin{tabular}{c}
\includegraphics[width=0.35\textwidth,angle=270]{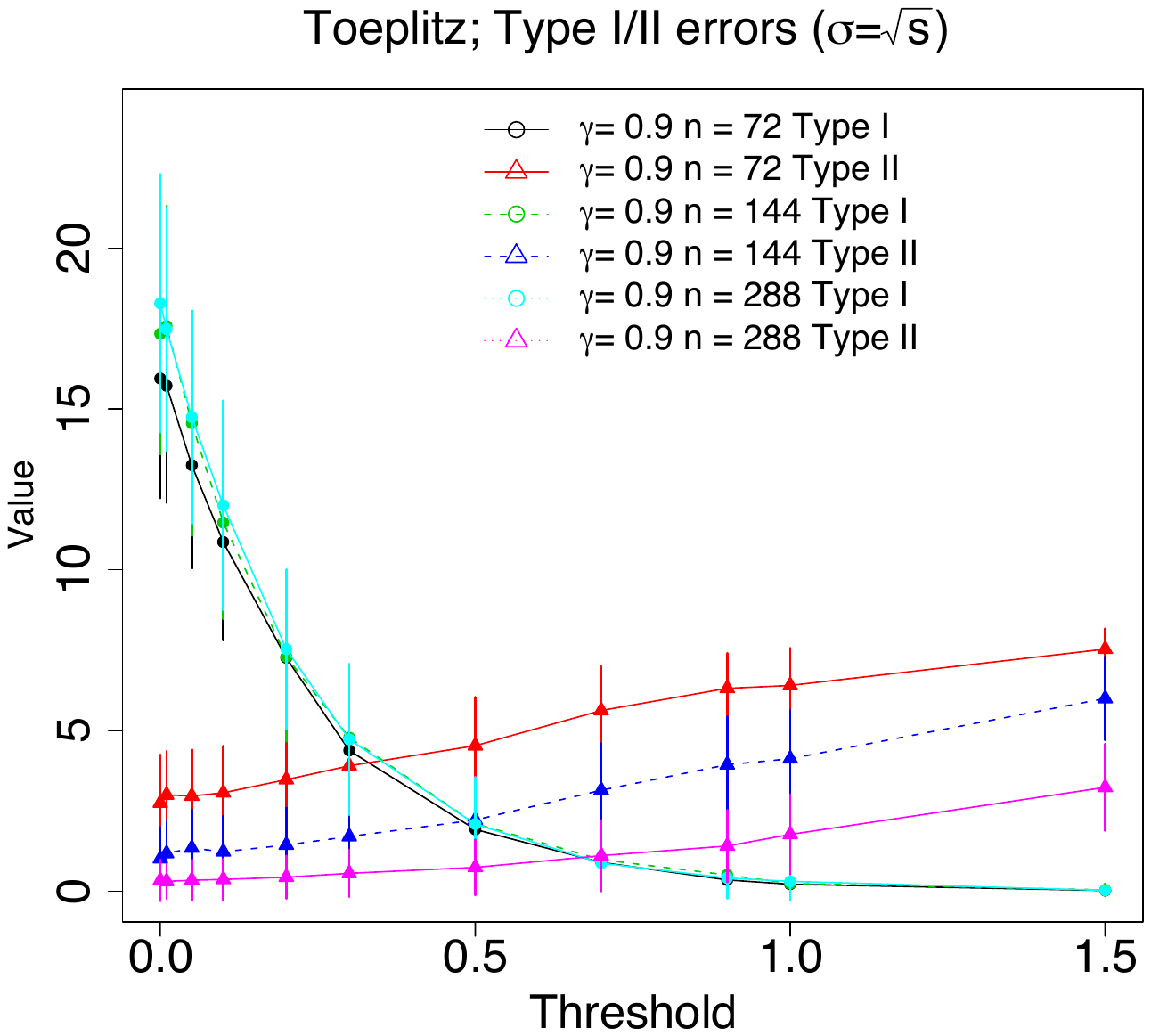}
\end{tabular}& 
\begin{tabular}{c}
\includegraphics[width=0.35\textwidth,angle=270]{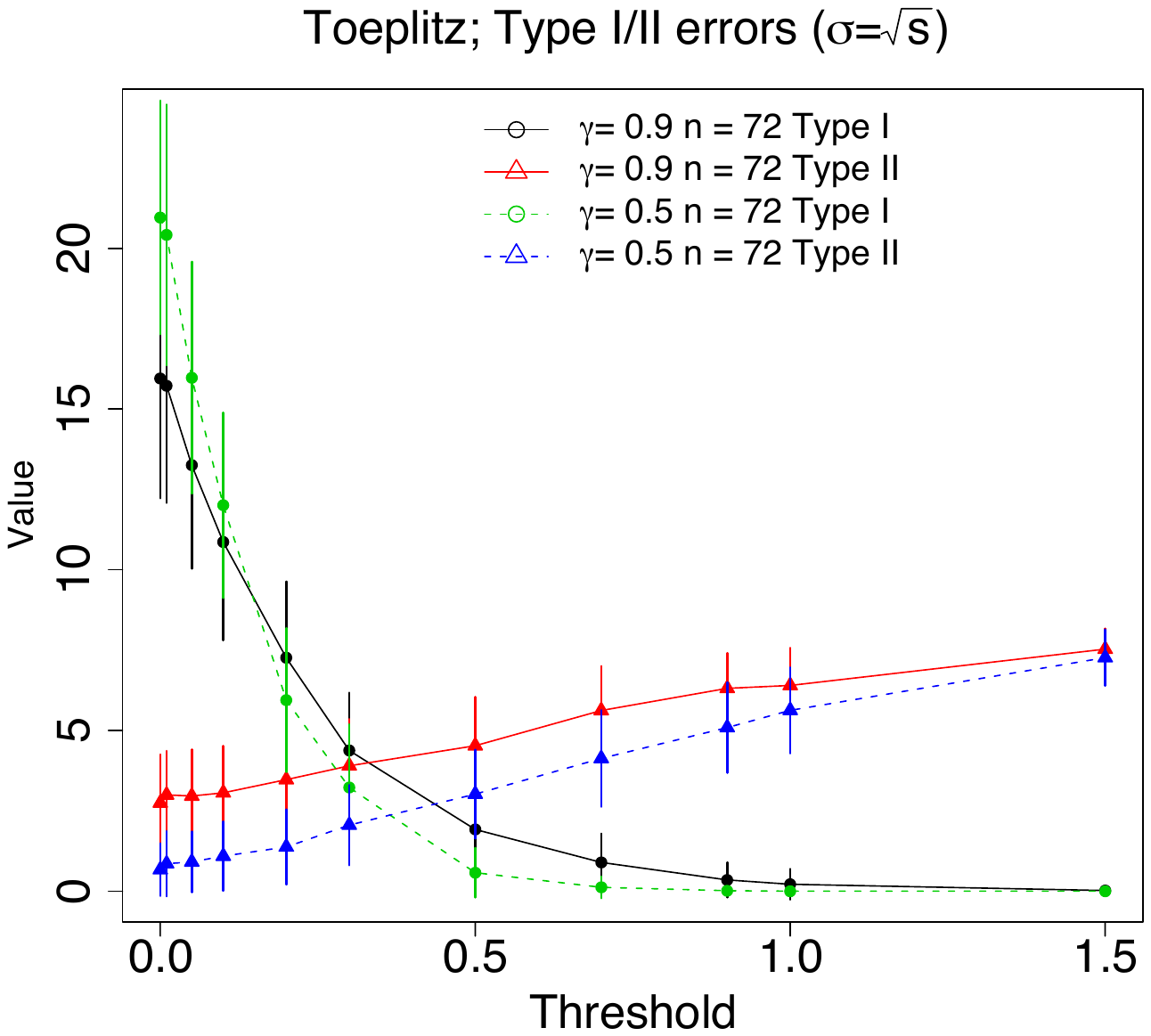} \\
\end{tabular} \\
(c) & (d) \\
\begin{tabular}{c}
\includegraphics[width=0.35\textwidth,angle=270]{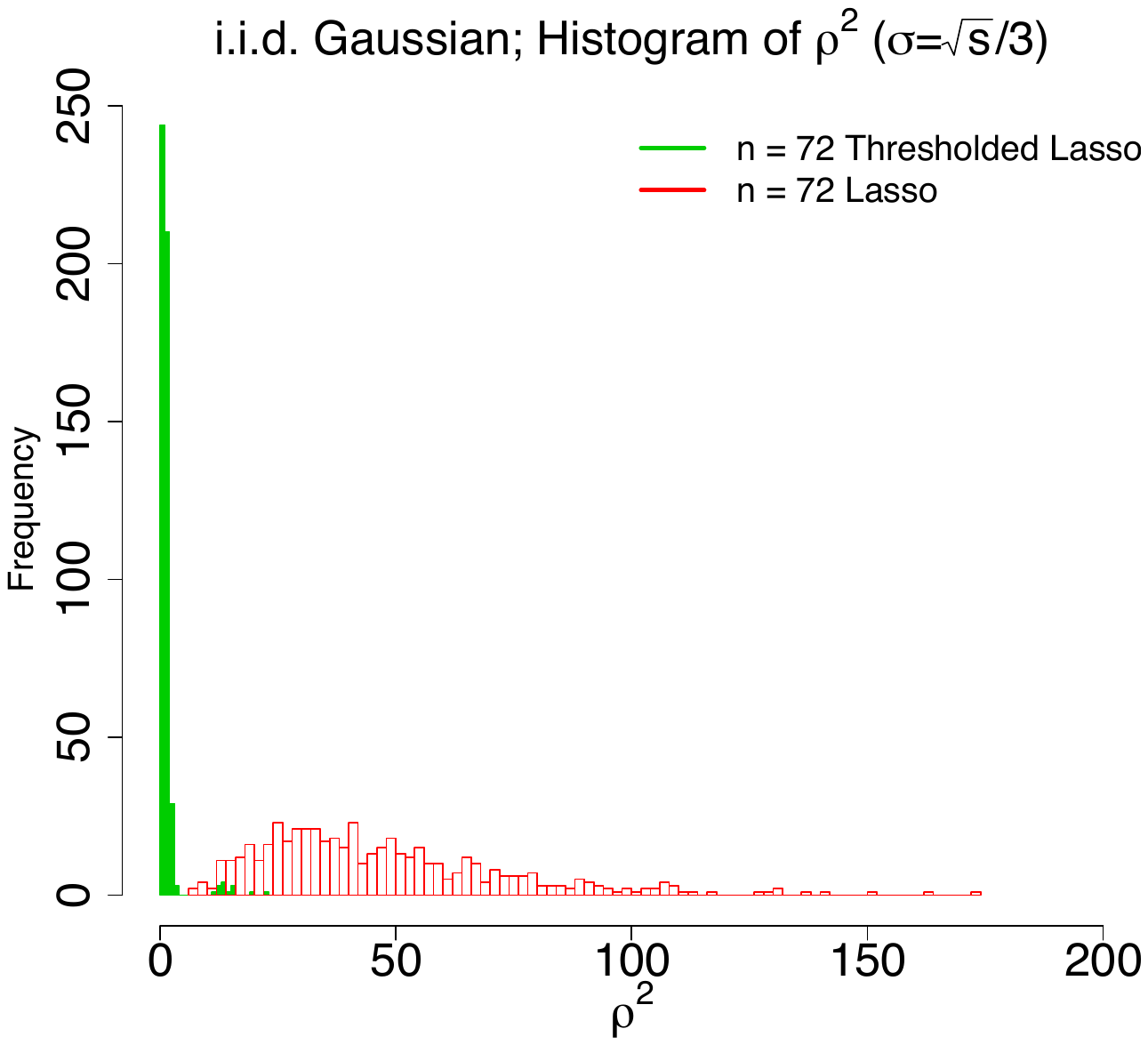}
\end{tabular}& 
\begin{tabular}{c}
\includegraphics[width=0.35\textwidth,angle=270]{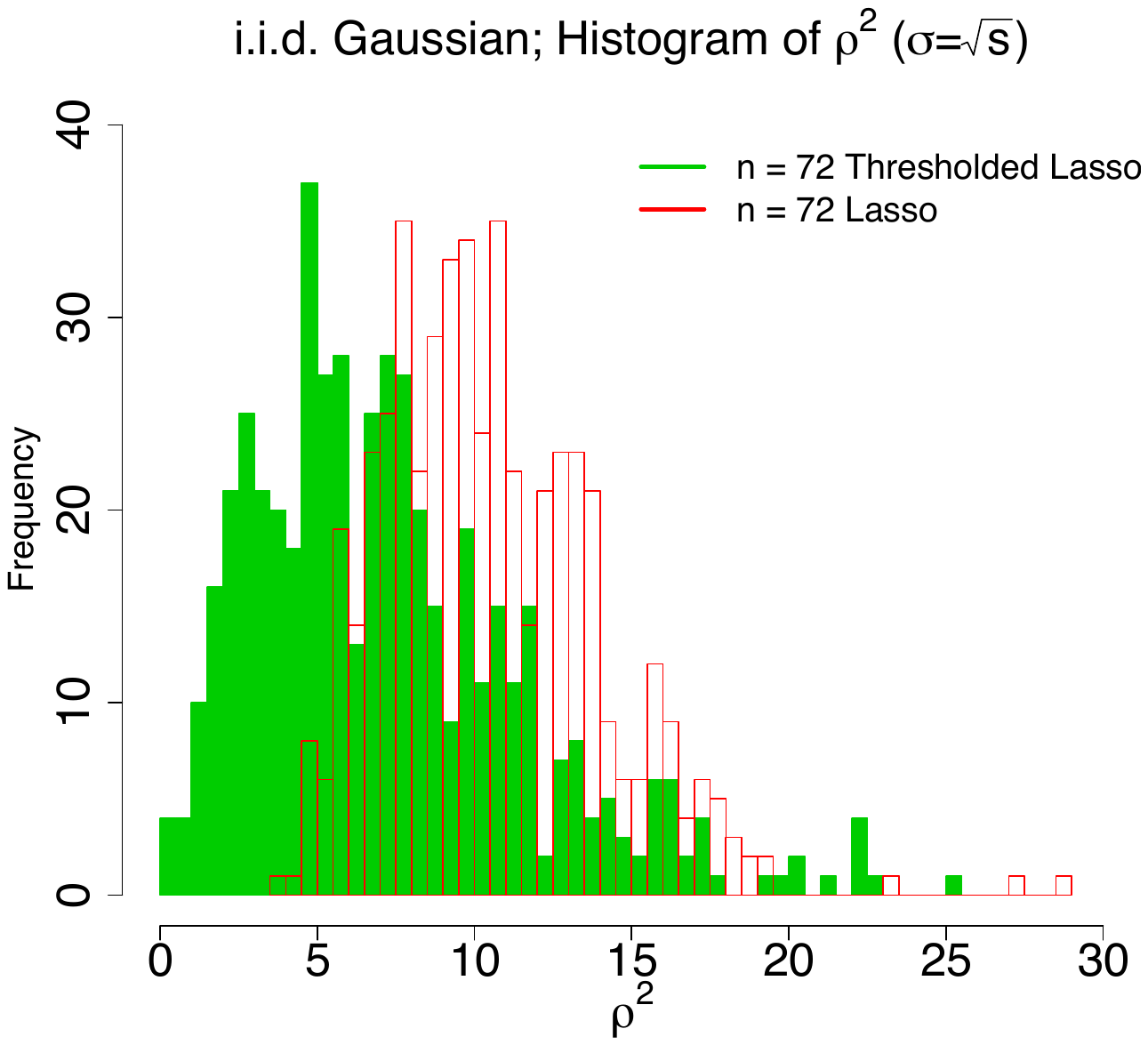} 
\end{tabular} \\
(e) & (f) \\
\end{tabular}
\caption{$p=256$ $s=8$. 
(a) (b) Type I/II errors for i.i.d. Gaussian and Toeplitz ensembles.
Each vertical bar represents $\pm 1$ std. The unit of $x$-axis is
in $\lambda \sigma$. For both types of design matrices,
FPs decrease and FNs increase as the threshold increases. 
For Toeplitz ensembles, in (c) with fixed correlation $\gamma$, 
FNs decrease with more samples, and in (d) with fixed sample size, 
FNs decrease as the correlation $\gamma$ decreases.
(e) (f) Histograms of $\rho^2$ under i.i.d Gaussian ensembles 
from 500 runs.}
\label{fig:type12}
\end{center}
\end{figure}
We first run the above experiment using i.i.d. Gaussian ensemble
under the following thresholds: $t_0 = \lambda \sigma$ for $\sigma=\sqrt{s}/3$,
and $t_0 = 0.36 \lambda \sigma$ for $\sigma=\sqrt{s}$. These are chosen based on the desire 
to have low errors of both types (as shown in Figure~\ref{fig:type12} (a)). 
Naturally, for low SNR cases,  small $t_0$ will reduce Type II errors.
In practice, we suggest using cross-validations to choose the exact 
constants in front of $\lambda \sigma$; 
See, for example, a subsequent work~\cite{ZRXB11} for details.
We plot the histograms of $\rho^2$ in Figure~\ref{fig:type12} (e) and (f).
In (e), the mean and median are $1.45$ and $1.01$ for the Thresholded 
Lasso, and  $46.97$ and $41.12$ for the Lasso. 
In (f), the corresponding values are $7.26$ and $6.60$ for 
the Thresholded Lasso and $10.50$ and $10.01$ for the Lasso.
With high SNR, the Thresholded Lasso performs extremely well;
with low SNR, the improvement of the Thresholded Lasso over the
ordinary Lasso is less prominent; this is in close correspondence with the
Gauss-Dantzig selector's behavior as shown by~\cite{CT07}.

Next we run the above experiment under different sparsity values of $s$. 
We again use i.i.d. Gaussian ensemble with 
$p=2000$, $n=400$, and $\sigma =\sqrt{s}/3$. 
The threshold is set at $t_0 = \lambda\sigma$. 
The SNR for different $s$ is fixed at around $32.36$. 
Table~\ref{tab:rho-square-snr} shows the mean of the $\rho^2$ for the 
Lasso and the Thresholded Lasso estimators. 
The Thresholded Lasso performs consistently better than the ordinary Lasso 
until about $s=80$, after which both break down.
For the Lasso, we always choose from the full regularization path
the {\em optimal} $\tilde{\beta}$ that has the minimum $\ell_2$ loss.

\begin{table}[h]
\begin{center}
\caption{$\rho^2$ under different sparsity and fixed SNR. 
Average over 100 runs for each $s$.
}
\label{tab:rho-square-snr}
\begin{tabular}{cccccccc} 
\hline
s      & 5    & 18   & 20   & 40   & 60   & 80   & 100  \\ \hline 
SNR & 34.66 &  32.99 &  32.29&   32.08 &  32.28 &  32.56 &  32.54  \\ \hline 
Lasso  & 17.42 &  22.01 &  44.89 &  52.68 &  31.88 &  29.40 &  47.63 \\ \hline 
Thresholded Lasso  & 1.02& 0.96& 1.11& 1.54& 10.32& 29.38 & 53.81 \\ \hline 
\end{tabular}
\end{center}
\end{table}

\subsection{Linear Sparsity}
We next present results demonstrating that the Thresholded Lasso 
recovers a sparse model using a small number of samples per 
non-zero component in $\beta$ when $X$ is a subgaussian ensemble.  
We run under three cases of $p = 256, 512, 1024$; for each $p$, 
we increase the sparsity $s$ by roughly equal steps from 
$s= {0.2 p}/{\log (0.2 p)}$ to $p/4$.  For each $p$ and $s$, we run with different
sample size $n$.  For each tuple $(n, p, s)$, we run an experiment similar to
the one described in Section~\ref{subsec:type12} with an i.i.d. Gaussian 
ensemble $X$ being fixed while repeating Steps $2-3$ 100 times.
In Step 2, each randomly selected non-zero coordinate of $\beta$ is 
assigned a value of $\pm 0.9$ with probability $1/2$.
After each run, we compare $\hat{\beta}$ with the true $\beta$; if all
components match in signs, we count this experiment as a success. 
At the end of the 100 runs, we compute the percentage of successful runs 
as the probability of success.  We compare with the ordinary Lasso, 
for which we search over the full regularization path of LARS and
choose the $\breve{\beta}$ that best matches $\beta$ in terms of support.

\begin{figure}
\begin{center}
\begin{tabular}{cc}
\begin{tabular}{c}
\includegraphics[width=0.34\textwidth,angle=270]{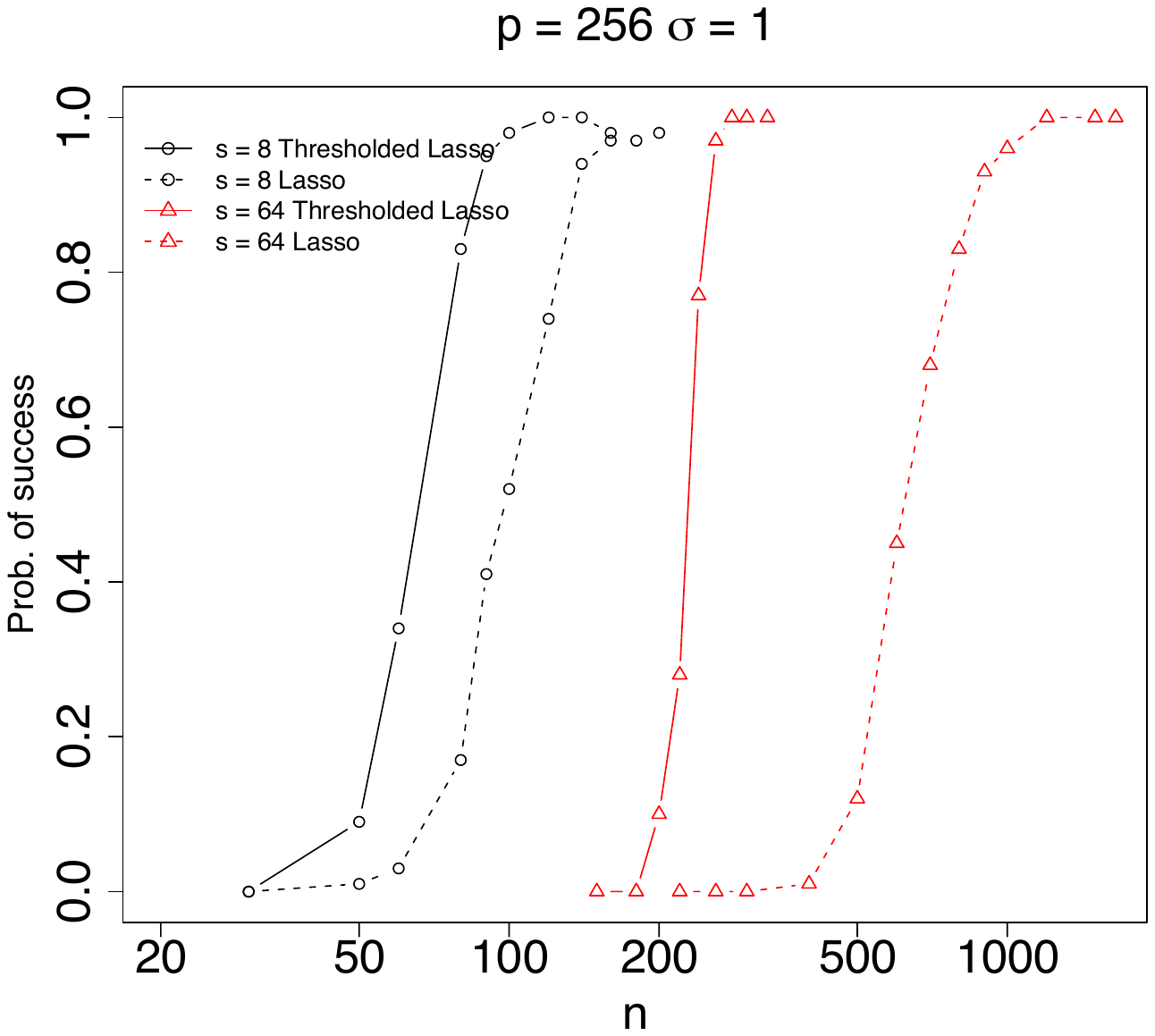}
\end{tabular}&
\begin{tabular}{c}
\includegraphics[width=0.34\textwidth,angle=270]{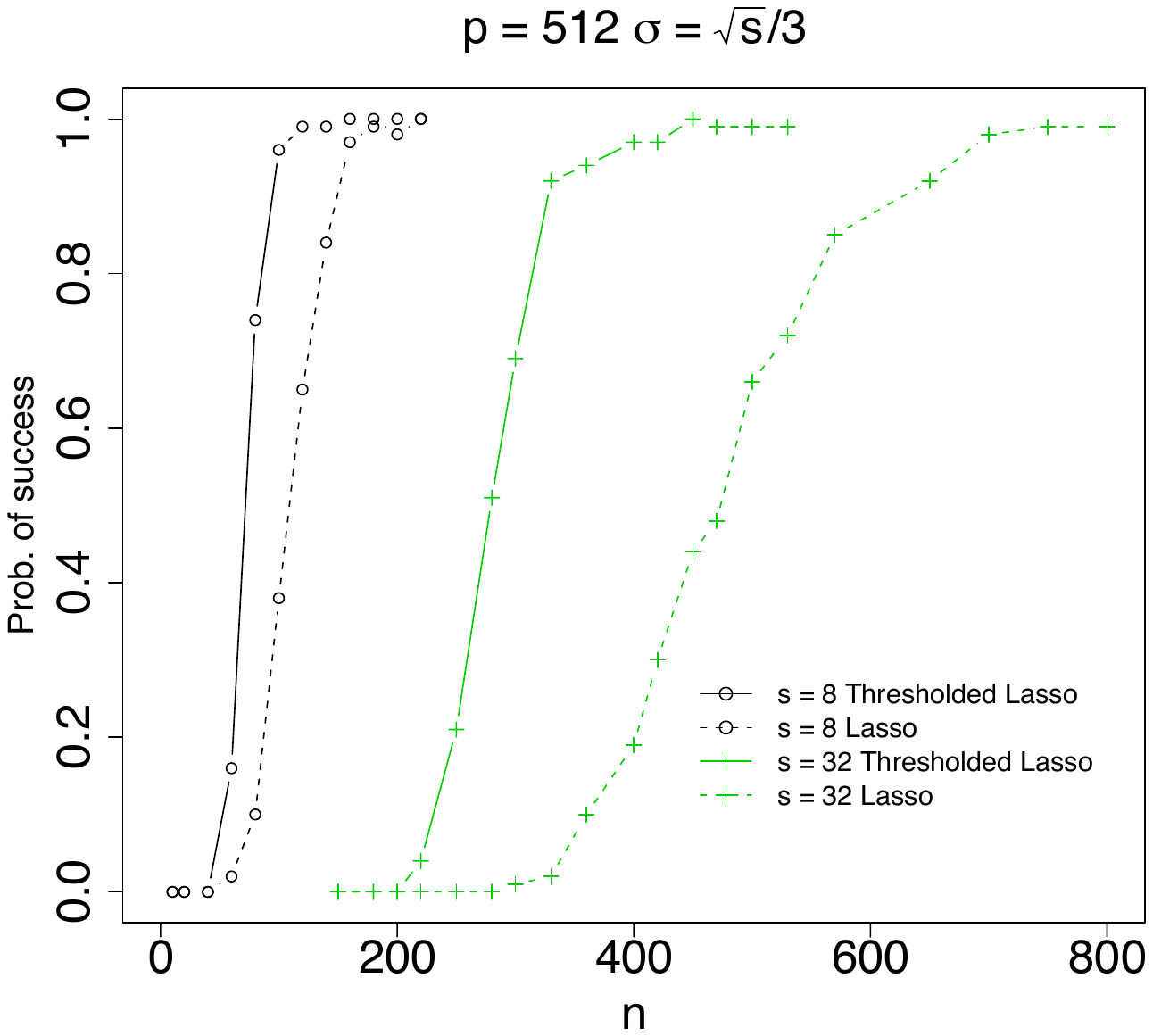}
\end{tabular} \\
(a)&(b) \\
\begin{tabular}{c}
\includegraphics[width=0.34\textwidth,angle=270]{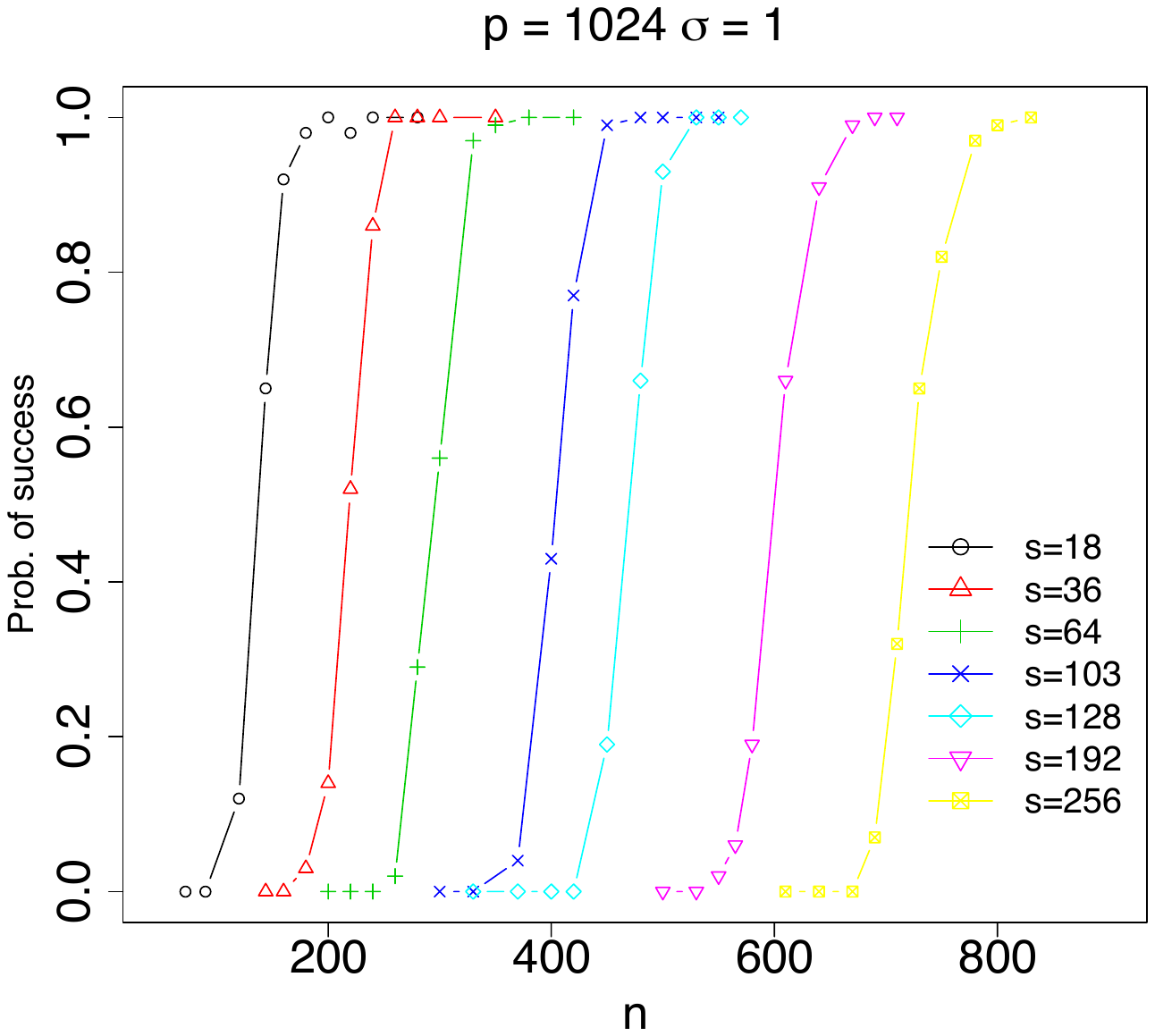}
\end{tabular}&
\begin{tabular}{c}
\includegraphics[width=0.34\textwidth,angle=270]{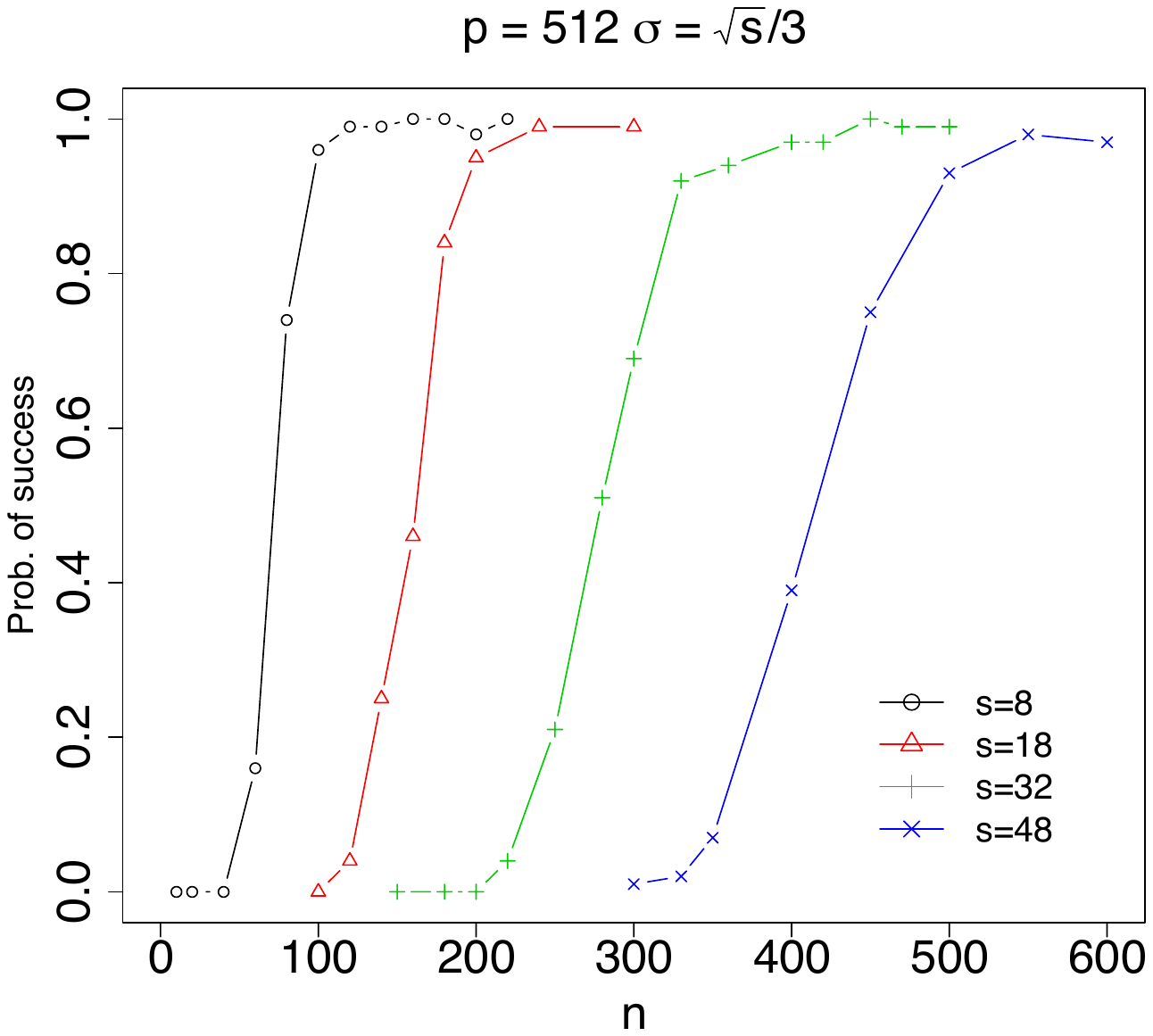}
\end{tabular} \\
(c)&(d) \\
\begin{tabular}{c}
\includegraphics[width=0.34\textwidth,angle=270]{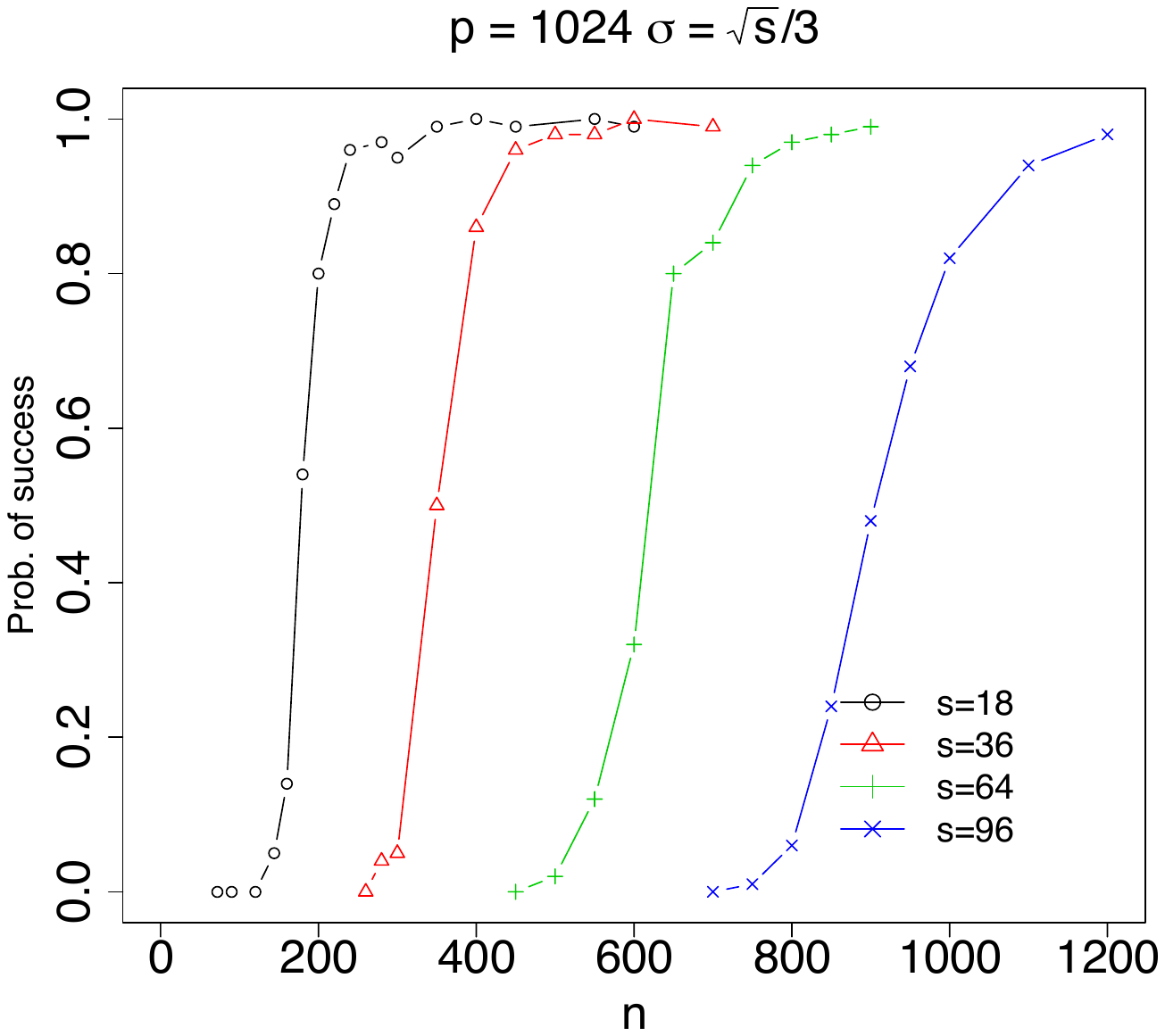}
\end{tabular}&
\begin{tabular}{c}
\includegraphics[width=0.34\textwidth,angle=270]{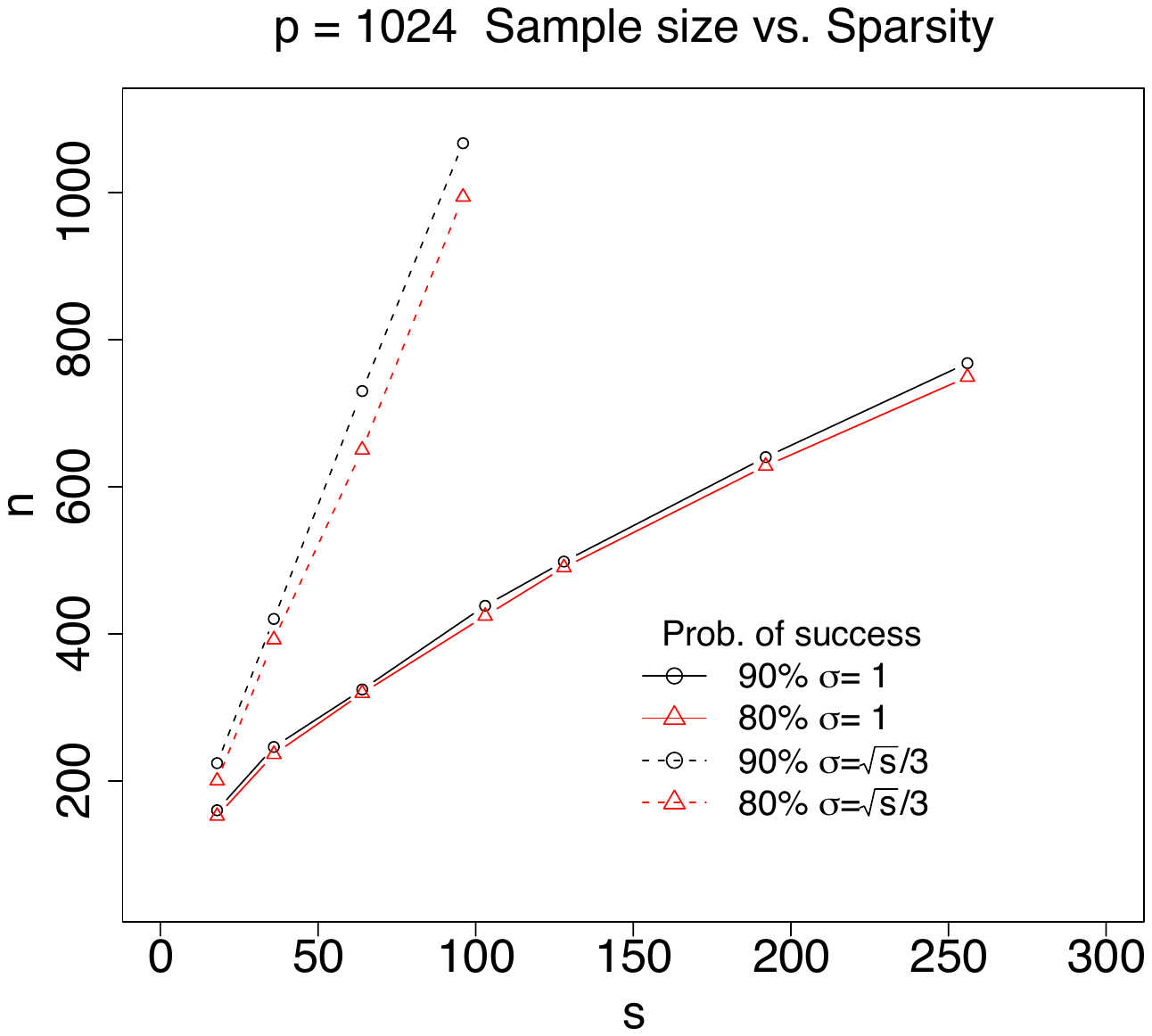}
\end{tabular} \\
(e)&(f) \\
\end{tabular}
\caption{(a) (b)
Compare the probability of success for $p = 256$ and $p=512$ under two 
noise levels. The Thresholded Lasso estimator requires much fewer samples
than the ordinary Lasso.
(c) (d) (e) show the probability of success of the Thresholded Lasso under 
different levels of sparsity and noise levels when $n$ increases for
$p =512$ and $1024$. (f) The number of samples $n$ increases almost linearly with 
$s$ for p = 1024. More samples are required to achieve the same level of success 
when $\sigma=\sqrt{s}/3$ due to the relatively low SNR.
}
\label{fig:succ-p256-four}
\end{center}
\end{figure}

We experiment with $\sigma =1$ and $\sigma = \sqrt{s}/3$.  
The results are shown in Figure~\ref{fig:succ-p256-four}. 
We observe that under both noise levels, 
the Thresholded Lasso estimator requires much fewer samples than 
the ordinary lasso in order to conduct exact recovery of the sparsity 
pattern of the true linear model when all non-zero components are 
sufficiently large.  When $\sigma$ is fixed as $s$ increases, the
SNR is increasing; the experimental results illustrate the behavior of 
sparse recovery when it is  close to the noiseless setting.  
Given the same sparsity,  more samples are required for the low SNR case 
to reach the same level of success rate. 
Similar behavior was also observed for Toeplitz and Bernoulli 
ensembles with i.i.d. $\pm 1$ entries.

\subsection{ROC comparison}
\begin{figure}
\begin{center}
\begin{tabular}{cc}
\begin{tabular}{c}
\includegraphics[width=0.34\textwidth,angle=270]{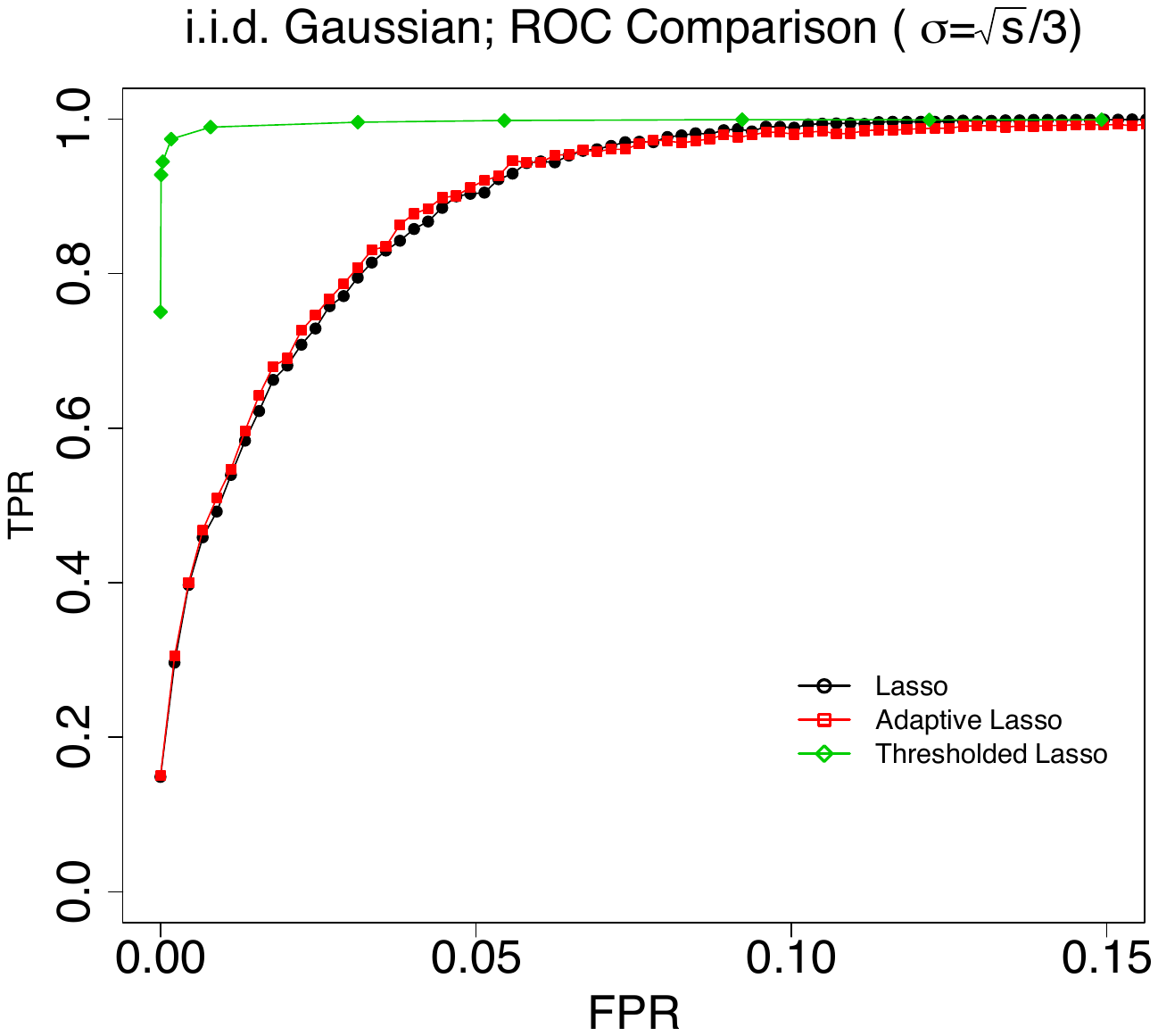}
\end{tabular}& 
\begin{tabular}{c}
\includegraphics[width=0.34\textwidth,angle=270]{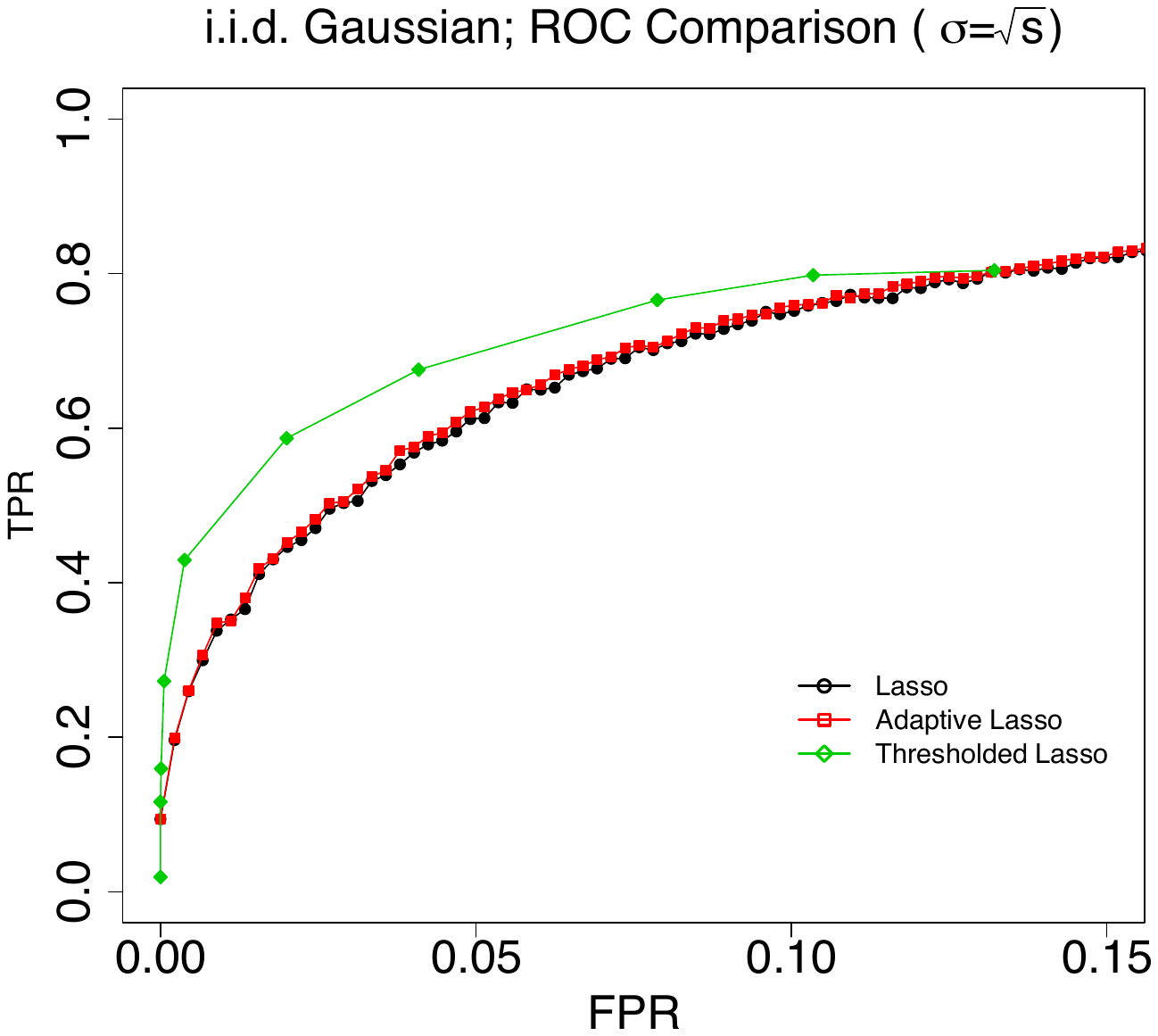} \\
\end{tabular} \\
(a) & (b) \\
\end{tabular}
\caption{$p=512$ $n=330$ $s=64$.  ROC for the Thresholded Lasso, ordinary Lasso and Adaptive
Lasso. The Thresholded Lasso clearly outperforms the ordinary Lasso and the
Adaptive Lasso for both high and low SNRs.
}
\label{fig:roc-gaussian}
\end{center}
\end{figure}

We now compare the performance of the Thresholded Lasso estimator
with the Lasso and the Adaptive Lasso by examining their ROC curves.
Our parameters are $p=512$, $n=330$, $s=64$ and we run under two cases:
$\sigma = \sqrt{s}/3$ and $\sigma = \sqrt{s}$.
In the Thresholded Lasso, we vary the threshold level from 
$0.01 \lambda \sigma$ to $1.5\lambda \sigma$. 
For each threshold, we run the experiment described in 
Section~\ref{subsec:type12} with an i.i.d. Gaussian ensemble $X$ being fixed
while repeating Steps $2-3$ 100 times. After each run, we compute the FPR
and TPR of the $\hat{\beta}$, and compute their averages after 100 runs
as the FPR and TPR for this threshold.  For the Lasso, 
we compute the FPR and TPR for each output vector along its entire
regularization path.  
For the Adaptive Lasso, we use the {\em optimal} output $\tilde{\beta}$ in terms 
of $\ell_2$ loss from the initial Lasso penalization path as the input to its 
second step, that is, we set $\beta_{\init} :=  \tilde{\beta}$
and use $w_j = 1/\beta_{\init, j}$ to compute the weights 
for penalizing those non-zero components in $\beta_{\init}$ in
the second step, while all zero components of $\beta_{\init}$ are now 
removed. We then compute the FPR and TPR for each vector
that we obtain from the second step's LARS output.
We implement the algorithms as given in ~\cite{Zou06}, the details of which 
are omitted here as its implementation has become standard.
The ROC curves are plotted in 
Figure~\ref{fig:roc-gaussian}. The Thresholded Lasso performs better than
both the ordinary Lasso and the Adaptive Lasso; its advantage is 
more apparent when the SNR is high.

\end{document}